\documentclass[psamsfonts, 12pt]{amsart}

\usepackage{amssymb,amsfonts}
\usepackage[all,arc]{xy}
\usepackage{enumerate}
\usepackage{mathrsfs}
\usepackage{hyperref}
\usepackage{indentfirst} 
\usepackage{tikz}
\usepackage{dsfont}
\usetikzlibrary{matrix,chains, arrows, decorations.markings, calc, patterns,intersections}
\usepackage{extarrows}

\usepackage{graphicx}
\usepackage{float}

\newtheorem{thm}{Theorem}[section]

\newtheorem{prop}[thm]{Proposition}
\newtheorem{lem}[thm]{Lemma}

\theoremstyle{definition}
\newtheorem{defn}[thm]{Definition}

\theoremstyle{remark}
\newtheorem{rem}[thm]{Remark}

\makeatletter
\let\c@equation\c@thm
\makeatother
\numberwithin{equation}{section}

\DeclareSymbolFont{yh}{OMX}{yhex}{m}{n}
\DeclareMathAccent{\hu}{\mathord}{yh}{"F3}

\title{Subgroups of Genus-2 Quasi-Fuchsian groups and Cocompact Kleinian Groups}

\author{Zhenghao Rao}

\usepackage{graphicx}
\usepackage{amsthm}
\usepackage{amsmath,amssymb}
\usepackage{amsfonts}
\usepackage{xcolor}
\usepackage[margin=1in]{geometry}
\usepackage[shortlabels]{enumitem}
\usepackage[bottom]{footmisc}
\usepackage{url}

\newcommand{\dmo}{\DeclareMathOperator}
\dmo{\id}{id}
\dmo{\pt}{pt}
\dmo{\ab}{ab}
\dmo{\GL}{GL}
\dmo{\SO}{SO}
\dmo{\SU}{SU}
\dmo{\SL}{SL}
\dmo{\PSL}{PSL}
\dmo{\im}{Im}
\dmo{\re}{Re}
\dmo{\Homeo}{Homeo}
\newcommand{\C}{\mathbb C}
\newcommand{\R}{\mathbb R}

\newcommand{\Z}{\mathbb Z}

\newcommand{\HH}{\mathbb H}

\newcommand{\ep}{\epsilon}
\newcommand{\de}{\delta}
\newcommand{\al}{\alpha}
\newcommand{\be}{\beta}

\newcommand{\De}{\Delta}

\newcommand{\Arg}{\mathrm {Arg}}

\newcommand{\hl}{\mathbf{{hl}}}
\newcommand{\lle}{\mathbf{{l}}}
\dmo{\Hom}{Hom}
\dmo{\Ext}{Ext}
\dmo{\Tor}{Tor}

\date{\today}

\begin{document}

\setlength{\parindent}{1.5em}

\begin{abstract}
In this paper, we want to control the geometry of some surface subgroups of a cocompact Kleinian group. More precisely, provided any genus-2 quasi-Fuchsian group $\Gamma$ and cocompact Kleinian group $G$, then for any $K>1$, we will find a surface subgroup $H$ of $G$ that is $K$-quasiconformally conjugate to a finite index subgroup $F<\Gamma$.

\end{abstract}

\maketitle

\tableofcontents

\section{Introduction}

\noindent \textbf{Main Result.}

In this paper, we want to study the geometry of cocompact hyperbolic 3-manifolds $M^3=\HH^3/G$. We say that $H<G$ is a surface subgroup if there exist a closed surface $S_g$ and a continuous map $f:S_g\mapsto M^3$ such that the induced homomophism $f_*$ is injective and $f_*(\pi_1(S_g))=G$. And we say a surface subgroup $H<G$ is $K$-quasi-Fuchsian if there is a cocompact Fuchsian group of \texorpdfstring{$\PSL(2,\R)$}{Lg} whose action on $\partial\HH^3$ is $K$-quasiconformally conjugate to the action of $H$ on $\partial\HH^3$. Here is our main result.

\begin{thm}\label{mr}
Let $\Gamma$ be a genus-2 quasi-Fuchsian group and $G$ be a cocompact Kleinian group. For any $M>1$, there is a surface subgroup $H<G$ that is $M$-quasiconformally conjugate to a finite index subgroup $F<\Gamma$. 
\end{thm}

\noindent \textbf{Related Results.}

Kahn and Markovic proved the original surface subgroup conjecture, that every closed hyperbolic 3-manifold $M^3$ contains a immersed $\pi_1$-injective surface in \cite{KM12b} and also proved that the number of genus $g$ surface subgroups grows as $g^{2g}$ in \cite{KM12a}. And in \cite{KM15}, they gave a proof of the Ehrenpreis conjecture that for any $K>1$, any two closed Riemann surfaces of genus greater than 1 have finite-degree covers which are $K$-quasiconformal. 

We say that a collection of quasi-Fuchsian surface
subgroups is ubiquitous if, for any pair of hyperbolic planes $\Pi,\Pi'$ in $\mathbb{H}^3$ with distance $d(\Pi,\Pi') > 0$, there is a surface subgroup whose boundary circle lies between $\partial\Pi$ and $\partial\Pi'$.

Cooper and Futer used different methods than those used in \cite{KM12b} to prove that the set of closed immersed quasi-Fuchsian surfaces in a complete, finite-volume hyperbolic 3-manifold is ubiquitous in \cite{CF19}. Kahn and Wright showed that for a finite volume but not compact hyperbolic 3-manifold $M^3$ and for all $K>1$, the collection of $K$-quasi-Fuchsian surface subgroups is ubiquitous in \cite{KW21}.\\

\noindent \textbf{Motivation.}

As mentioned in \cite{KW21}, one can regard \cite{KM12b}, \cite{KM15} and \cite{KW21} as special cases of the general question: whether a lattice in a Lie group $\mathcal{G}$ contains surface subgroups which are close to lying in a given subgroup of $\mathcal{G}$ isomorphic to $\PSL(2,\R)$. Here we care about a similar question but $\PSL(2,\R)$ is replaced with a discrete group. And Theorem \ref{mr} is the solution to the case where $\mathcal{G}=\PSL(2,\C)$, $\Gamma$ is cocompact and the discrete group is a genus-2 quasi-Fuchsian group. We hope that good pants homology can help us remove the restriction of genus-2, and one can refer Remark \ref{rem8.1} for more details. We also expect that the methods developed in \cite{KW21} can be applied to this case with cofinite volume $\Gamma$. Other case including higher rank Lie group like $\mathcal{G}=\SL(n,\R)$ raise special interest as well. 

One may relate our result to Hausdorff dimension of limit sets of quasi-Fuchsian groups. In \cite{Bro03}, Brock proved that there are quasi-Fuchsian groups with the Hausdorff dimension of limit sets arbitrarily close to 2. Gehring and Väisälä discussed how the Hausdorff dimension changes under a quasiconformal mapping between 2-dimensional spaces in \cite{GV73}. Therefore we can have the following result:

\begin{thm}\label{1.2}

Given a cocompact Kleinian group $G$ and a real number $1\leq\al<2$. Then for any $\ep>0$, there exists a surface subgroup $H<G$ such that
\begin{equation*}
    |\textup{H-dim}(\lambda(H))-\al|<\ep,
\end{equation*}
where $\textup{H-dim}$ denotes the Hausdorff dimension.  

\end{thm}
Bowen had a similar result but with free groups in \cite{Bow09}, which says the set of Hausdorff dimensions of limit sets of free subgroups of $\PSL(2,\C)$ is dense in $[0,2]$.

One can also ask for a given $g$, how the genus-$g$ surface subgroups are distributed in a hyperbolic 3-manifold and how the Hausdorff dimensions of genus-$g$ surface subgroups are distributed.
\\

\noindent \textbf{Sketch of the Proof.}

Roughly speaking, we divide the proof into three parts.

(1) We take any non-separating pants decomposition $\sigma_1$ of the given genus-2 quasi-Fuchsian group and we use it to find a pants decomposition $\sigma_2$ with similar cuff lengths. And then we use a specific ideal triangulation of $\sigma_2$ to construct a good pants decomposition $\sigma_3$ such that the lengths of cuffs can be as large as we want and the real part of the three shears are positive and bounded from above and below.

(2) For a pair of pants $P$, $P$ can be determined by a cuff $\gamma$ and a third connection $\alpha$ to $\gamma$. By the Margulis argument, we can count the number of $(R_i,\epsilon)$-good curves. Moreover, to count $(R_i,\epsilon)_{i=1}^3$-good pants with $\gamma$ a boundary curve, we can count the corresponding third connection to $\gamma$ and prove that these pants are distributed along $\gamma$ almost evenly. Then by the Hall Marriage theorem, we can find a permutation with good behavior determined by the shears in (1), on the set of all wanted pants with boundary cuff $\gamma$. And by using the doubling trick, we can glue these pants along $\gamma$ together and finally get a good assembly.

(3) By generalizing the idea in the appendix to \cite{KW21}, we can prove that the good assembly we have in (2) is close to a quasi-Fuchsian surface which is a finite covering space of the given surface. Therefore Theorem (\ref{mr}) is proved.
\\

\noindent \textbf{Organization.}

Section \ref{S2} introduces basics of pants, ideal triangles for (quasi-)surfaces and assemblies, where some are recalled from \cite{KW21}. In Section \ref{S3}, we present the left and right actions of $\PSL(2,\C)$ on the frame bundle of $\HH^3$.

And in Section \ref{S4}, we build the inefficiency theory for 3-dimensional hyperbolic spaces, based on the work in \cite{KM15} and \cite{LM15}. The theory is used to study a specific family of broken geodesics in a hyperbolic 3-manifold, where the segments are alternately long enough and relatively short. And we also estimate the complex length of the geodesic representative for such a broken geodesic.

Section \ref{S5} constructs a desired pants decomposition of the given genus-2 quasi-Fuchsian surface, where the cuff lengths are large enough and of the same size and the complex shears have bounded and positive real parts. We use the $spinning$ $construction$ to build a new pants decomposition from a given non-separating one on a genus-2 surface. And then we take a special ideal triangulation to estimate the lengths and shears by a result in Section \ref{S4}.

Section \ref{S6} collects and generalizes results from \cite{KW21} on construction, counting and equidistribution of good pants, and we also generalize the idea of $good$ to utilize in our case.

In Section \ref{S7}, we apply the method in the appendix to \cite{KW21} to prove that a good assembly is close to the given quasi-Fuchsian group, which is also the heart of this paper. For a good assembly $\mathcal{A}$, we construct its perfect model $\hat{\mathcal{A}}$ and a map $e$ from $\hat{\mathcal{A}}$ to $\mathcal{A}$ which maps each perfect pants to the corresponding good pants. To give a quantitative estimate on $e$, we locally lift $e$ and compare lifts of assemblies in the universal cover, where we use the matrix multiplication estimate to obtain a desired geometric estimate called $\ep$-bounded distortion up to distance $D$. Then we can extend the local estimate to a global result by some more general machinery.
\\

\noindent \textbf{Acknowledgment.}

The author would like to thank his advisor Jeremy Kahn for patient guidance and helpful comments.

\section{Right-angled hexagon and pants, ideal triangles and good assemblies}\label{S2}

\subsection{Right-angled hexagons and pants} A hyperbolic right-angled hexagon in $\HH^3$ is determined by the complex lengths of any three mutually non-adjacent sides, which can take any complex values with positive real parts. For each side $e$ of the hexagon, the real part of the complex length of $e$ denotes the length of $e$, and the imaginary part measures the rotation between two tangent vectors to the incoming and outgoing edges adjacent to $e$.

For a pair of pants in a hyperbolic 3-manifold $M^3$, it is composed of two isometric right-angled hexagons joined along three mutually non-adjacent sides. And the complex lengths of the others three sides, which form the boundary components of the pair of pants, determines this pair of pants.

For a quasi-Fuchsian group $\Gamma$ of genus $g$, let $\mathcal{C}$ be a maximal set of disjoint non-trivial, non-homotopic simple closed curves of the genus-$g$ surface; then $\mathcal{C}$ has $3g-3$ elements. We can construct an immersed surface $S_{\Gamma, \mathcal{C}}$ relative to $\mathcal{C}$ in the convex core of $\Gamma$ such that for each curve $\gamma$ in $\mathcal{C}$ there exists a geodesic freely homotopic to $\gamma$, and each component of $S_{\Gamma,\mathcal{C}}-\mathcal{C}$ (here $\mathcal{C}$ consists of geodesic representatives in the manifold) is a pair of pants with geodesic boundaries. And we say this pants decomposition of $S_{\Gamma,\mathcal{C}}$ is a pants decomposition of $\Gamma$ relative to $\mathcal{C}$. For more details, one can check \cite{Tan94}.

For a connected surface, a pants decomposition is called $\mathit{non}$-$\mathit{separating}$, if when we remove any one cuff of this pants decomposition, the remaining part of the surface is still connected. In particular, if the genus of the surface is $2$, a non-separating pants decomposition will have two identical pairs of pants.

Before we define the shears between the pairs of pants sharing a common cuff, we need the definitions of normal bundles, orthogeodesics, feet, connections and third connections, and one can see Section 2.1, 2.7 and 3.2 in \cite{KW21}, which also apply for pants in hyperbolic 3-manifolds. We also call the base point of a foot (as a unit vector) a footpoint. 

Now let $P_1,P_2$ be two pairs of pants which share a common cuff $K$, $A_1,A_2$ and $B_1,B_2$ be the remaining cuffs of $P_1$ and $P_2$ respectively, and $a_1,a_2,b_1,b_2$ be the footpoints of orthogeodesics between $A_1,A_2,B_1,B_2$ and $K$ respectively. Given the orientation of $K$, we define the $\mathit{short\ shear}$ $s$ between $P_1$ and $P_2$ be the complex distance, which is in $\C/\hl(K)\Z+2\pi i\Z$, between the orthogeodesics from $A_1,B_1$ to $K$, where $\hl(K)$ is the half-length of $K$. The fact that the complex distance between orthogeodesics from $A_1,A_2$ to $K$ is $\hl(K)$ guarantees that this definition is well-defined.

Next, let $\al,\be$ be the third connections of $K$ in $P_1$ and $P_2$ with footpoints $\al_1,\al_2,\be_1,\be_2$ satisfying $\al_1$ is between $a_1$ and $a_2$ and $\be_1$ is between $b_1$ and $b_2$ along $K$. Denote the complex distance from $\al_i$ to $\be_i$ in $\C/2\hl(K)\Z+2\pi i \Z$ by $t_i$, and we call $(t_1,t_2)$ the $\mathit{long\ shear}$ between $P_1$ and $P_2$. It is easy to see that $t_1+t_2=2s$ in $\C/2\hl(K)\Z+2\pi i\Z$.

\subsection{Ideal triangulation} An $\mathit{ideal\ triangle}$ in a hyperbolic 3-manifold $M^3$ is the image of an injective, local isometry from an ideal triangle in $\HH^3$ to $M$. Moreover, an $\mathit{ideal\ triangulation}$ of $S$ is a lamination with finitely many leaves whose complementary components are ideal triangles.

For two ideal triangles which share a common edge, we define the $\mathit{shear}$ between them along this edge as follows: drop the altitudes at the remaining vertices to the common edge, then the complex distance between two footpoints is called the shear along this edge.

For more details, one can check \cite{BBFS13} and \cite{Z17}.

\section{Right action of \texorpdfstring{$\PSL(2,\C)$}{Lg} on the frame bundle of hyperbolic 3-space}\label{S3}

In this paper, we always use the upper half plane model for $\HH^2$ and the upper half space model for $\HH^3$.

In $\HH^2$, each point can be written as $z=x+iy$ with $x,y\in\R$ and $y>0$. And it turns out that $\PSL(2,\R)$ is the isometry group of $\HH^2$, where $\PSL(2,\R)$ left acts on $\HH^2$ by
\begin{align*}
    T(z)=\frac{az+b}{cz+d},
\end{align*}
where
\begin{align*}
    T=\begin{pmatrix}a & b \\ c & d\end{pmatrix} \in \PSL(2,\R).
\end{align*}
Since the left action is isometry, we have the left action of $\PSL(2,\R)$ on the unit tangent bundle $T^{1}(\HH^2)$. Based on the left action on $T^{1}(\HH^2)$, we want to describe another action of $\PSL(2,\R)$ on $T^{1}(\HH^2)$ which is called the right action.

We first fix the base point $z_0=i\in \HH^2$ and the base unit tangent vector $\mathbf{u}_0=(i,i)\in T_{z_0}^{1}(\HH^2)$. Then there is a unique right action as follows:
\begin{enumerate}

\item for any $g\in \PSL(2,\R)$, $g\cdot \mathbf{u}_0=\mathbf{u}_0\cdot g$;

\item for any $\mathbf{v}\in T^{1}(\HH^2)$ and $g,h\in \PSL(2,\R)$, $g\cdot(\mathbf{v}\cdot h)=(g\cdot\mathbf{v})\cdot h$.

\end{enumerate}
It is easy to see that this is a transitive and faithful action. 

Since the left action is transitive, then for any $\mathbf{v}\in T^{1}(\HH^2)$, there is a $g\in \PSL(2,\R)$ such that $\mathbf{v}=g\cdot \mathbf{u}_0$. Thus for $h\in \PSL(2,\R)$, we have
\begin{align*}
    \mathbf{v}\cdot h=(g\cdot \mathbf{u}_0)\cdot h=g\cdot(\mathbf{u}_0\cdot h).
\end{align*}
Hence we can use the right action of $h$ on $\mathbf{u}_0$, which is the same as the left action of $h$ on $\mathbf{u}_0$, to describe the right action of $h$ on $\mathbf{v}=g\cdot \mathbf{u}_0$. Here are some useful examples:
\begin{enumerate}[(a)]

\item Geodesic flow. Let 
\begin{align*}
    A(t)=\begin{pmatrix} e^{t} & 0 \\ 0 & e^{-t}\end{pmatrix} \in \PSL(2,\R),
\end{align*}
then $\mathbf{u}_0\cdot A(t)=A(t)\cdot\mathbf{u}_0=(ie^{2t},ie^{2t})$. Therefore $A(t):T^{1}(\HH^2)\to T^{1}(\HH^2)$ is the geodesic flow.

\item Rotation. Let 
\begin{align*}
    B(\theta)=\begin{pmatrix} \cos(\theta/2) & \sin(\theta/2) \\ -\sin(\theta/2) & \cos(\theta/2)\end{pmatrix} \in \PSL(2,\R),
\end{align*}
then $\mathbf{u}_0\cdot B(\theta)=B(\theta)\cdot\mathbf{u}_0=(i,-\sin\theta+i\cos\theta)$. Thus the action of $B(\theta)$ on any $\mathbf{v}\in T^{1}(\HH^2)$ is the counterclockwise rotation by $\theta$.

\end{enumerate}

For an oriented dim-$n$ Riemannian manifold $M$, a point $x\in M$ together with an orthonormal basis of $T_x M$ with positive orientation is called an $n$-frame in $M$. And we denote by $\mathcal{F}M$ the set of all $n$-frames in $M$, which forms a fiber bundle over $M$ and is called the frame bundle of $M$. Now we can define the right action of $\PSL(2,\C)$ on the frame bundle $\mathcal{F}\HH^3$ of $\HH^3$ based on the corresponding left action by isometries. Take the upper half space model of $\HH^3=\{(z,t):z\in\C,t>0\}$ and let
\begin{equation*}
    \mathcal{F}\HH^3=\{\langle p,u,v\rangle:p\in\HH^3,u,v\in T_p^1 \HH^3,u\perp v\}.
\end{equation*}
We also fix the base frame 
\begin{equation*}
    \Psi_0=\langle p_0,\mathbf{u}_0,\mathbf{v}_0\rangle=\langle(0,1),(0,1),(1,0)\rangle
\end{equation*}
at $(0,1)\in\HH^3$, and then the right action of $\PSL(2,\C)$ on $\mathcal{F}\HH^3$ is defined as follows:
\begin{enumerate}

\item for any $g\in \PSL(2,\C)$, $g\cdot \Psi_0=\Psi_0\cdot g$;

\item for any $\Psi\in \mathcal{F}\HH^3$ and $g,h\in \PSL(2,\C)$, $g\cdot(\Psi\cdot h)=(g\cdot\Psi)\cdot h$.

\end{enumerate}

For $\Psi\in\mathcal{F}\HH^3$, there exists $g\in\PSL(2,\C)$ such that $\Psi=g\cdot\Psi_0$. Then for any $h\in\PSL(2,\C)$, we have
\begin{align*}
    \Psi\cdot h=(g\cdot \Psi_0)\cdot h=g\cdot(\Psi_0\cdot h).
\end{align*}
Thus the right action on $\Psi_0$ will tell us the right action on any other frames. Here are the examples of geodesic flow and rotations, which will play important roles in the future.

\begin{enumerate}[(a)]

\item Geodesic flow. Let 
\begin{align*}
    A(x)=\begin{pmatrix} e^{x} & 0 \\ 0 & e^{-x}\end{pmatrix}
\end{align*}
with $x\in\R$. Then 
\begin{equation*}
    \Psi_0\cdot A(x)=A(x)\cdot\Psi_0=\langle(0,e^{2x}),(0,e^{2x}),(e^{2x},0)\rangle.
\end{equation*}
Therefore $A(x):\mathcal{F}\HH^3\to \mathcal{F}\HH^3$ is the geodesic flow.

\item Rotation along $\mathbf{u}_0$. Let
\begin{align*}
    A(iy)=\begin{pmatrix} e^{iy} & 0 \\ 0 & e^{-iy}\end{pmatrix}
\end{align*}
with $y\in\R$. Then
\begin{equation*}
    \Psi_0\cdot A(iy)=A(iy)\cdot\Psi_0=\langle(0,1),(0,1),(e^{i2y},0)\rangle.
\end{equation*}
Hence $A(iy)$ is the rotation along $\mathbf{u}_0$.

\item Rotation along $\mathbf{u}_0\times\mathbf{v}_0$. Let 
\begin{align*}
    B(\theta)=\begin{pmatrix} \cos(\theta/2) & \sin(\theta/2) \\ -\sin(\theta/2) & \cos(\theta/2)\end{pmatrix} \in \PSL(2,\R),
\end{align*}
with $\theta\in\R$. Then
\begin{equation*}
    \Psi_0\cdot B(\theta)=B(\theta)\cdot\Psi_0=\langle(0,1),(-\sin(\theta),\cos(\theta)),(\cos(\theta),\sin(\theta))\rangle.
\end{equation*}
So the action of $B(\theta)$ is the rotation along $\mathbf{u}_0\times\mathbf{v}_0$.

\end{enumerate}

\section{The theory of inefficiency in three dimension}\label{S4}

In this section, we will build the theory of inefficiency in 3-dimensional hyperbolic manifolds, which is used to determine the wanted pants decomposition of the quasi-Fuchsian group in the next section. Some definitions and results directly comes from Section 4 in \cite{KM15} and Section 4 in \cite{LM15}.

\subsection{Terminology}

Suppose $M$ is an oriented hyperbolic 3-manifold. We introduce some concepts in the geometry of framed geodesic segments.

\begin{defn}

An \textit{oriented\ framed\ segment} in $M$ is a triple:
\begin{equation*}
    \mathfrak{s}=(s,\vec{n}_{ini},\vec{n}_{ter}),
\end{equation*}
such that $s$ is an oriented immersed compact geodesic segment, and that $\vec{n}_{ini}$ and $\vec{n}_{ter}$ are two unit normal vectors at the initial endpoint and terminal endpoint of $s$, repectively.

• The \textit{carrier\ segment} is the oriented segment $s$;

• The \textit{initial\ endpoint} $p_{ini}(\mathfrak{s})$ and the \textit{terminal\ endpoint} $p_{ter}(\mathfrak{s})$ are the initial endpoint and the terminal endpoint of $s$, respectively;

• The \textit{initial\ framing} $\vec{n}_{ini}(\mathfrak{s})$ and the \textit{terminal\ framing} $\vec{n}_{ter}(\mathfrak{s})$ are the unit normal vectors $\vec{n}_{ini}$ and $\vec{n}_{ter}$,

• The \textit{initial\ direction} $\vec{t}_{ini}(\mathfrak{s})$ and the \textit{terminal\ direction} $\vec{t}_{ter}(\mathfrak{s})$ are the unit tangent vectors in the direction of $s$ at the initial point and the terminal point, respectively.

The \textit{orientation\ reversal} of $\mathfrak{s}$ is defined to be
\begin{equation*}
    \mathfrak{\overline{s}}=(\overline{s},\vec{n}_{ter},\vec{n}_{ini}),
\end{equation*}
where $\overline{s}$ is the orientation reversal of $s$. The \textit{framing\ flipping} of $\mathfrak{s}$ is defined as
\begin{equation*}
    \mathfrak{s}^*=(s,-\vec{n}_{ini},-\vec{n}_{ter}).
\end{equation*}

\end{defn}

It follows from the definition that
\begin{equation*}
    \overline{\mathfrak{s}^*}=\overline{s}^*.
\end{equation*}

\begin{defn}
For an oriented framed segment $\mathfrak{s}$ in $M$, the length of $\mathfrak{s}$, denoted as $l(\mathfrak{s})\in[0,+\infty)$, is the length of the unframed segment $s$ carrying $\mathfrak{s}$, and the phase of $\mathfrak{s}$, denoted as $\phi({\mathfrak{s}})\in \R/2\pi\Z$, is the angle from the initial framing $\vec{n}_{ini}$ to the parallel transpotation of $\vec{n}_{ter}$ to the initial endpoint of $\mathfrak{s}$ via $s$, signed with respect to the normal orientation induced from $\vec{t}_{ini}$ and the orientation of $M$. And we define the \textit{complex\ length} of $\mathfrak{s}$ as 
\begin{equation*}
    \lle(\mathfrak{s})=l(\mathfrak{s})+i\phi(\mathfrak{s}).
\end{equation*}
For an oriented closed geodesic curve $c$, we can also talk about its \textit{length,\ phase}, or \textit{complex\ length}, by taking an arbitrary unit normal vector $\vec{n}$ at a point $p\in c$, and regarding $c$ as a framed segment obtained by cutting $c$ at $p$ and assigned with framing $\vec{n}$ at both endpoints.
\end{defn}

It is  easy to see that length and phase are invariant under orientation reversal and under framing flipping.

\begin{defn}
Let $0\leq\de<\frac{\pi}{3}$, $L>d>0$, and $0<\theta<\pi$ be constants.
\begin{enumerate}
\item An oriented framed segments $\mathfrak{s}$ is said to be $\delta$-\textit{consecutive} to another oriented framed segment $\mathfrak{s}'$ if the terminal endpoint of $\mathfrak{s}$ is the initial endpoint of $\mathfrak{s}'$, and if the terminal framing of $\mathfrak{s}$ is $\de$-close to the initial framing of $\mathfrak{s}'$. We say $\mathfrak{s}$ is \textit{consecutive} to $\mathfrak{s}'$ if $\de$ is 0. When $\mathfrak{s}$ is $\delta$-consecutive to $\mathfrak{s}'$, the \textit{bending\ angle} between $\mathfrak{s}$ and $\mathfrak{s}'$ is the angle between the terminal direction of $\mathfrak{s}$ and the initial direction of $\mathfrak{s}'$, which is valued in $[0,\pi]$ and denoted by $\angle(\mathfrak{s},\mathfrak{s}')$.

\item A $\de$-\textit{consecutive\ chain} of oriented framed segments is a finite sequence $\mathfrak{s}_1\mathfrak{s}_2\cdots\mathfrak{s}_m$ such that each $\mathfrak{s}_i$ is $\de$-consecutive to $\mathfrak{s}_{i+1}$. It is a $\de$-\textit{consecutive\ cycle} if furthermore $\mathfrak{s}_m$ is $\de$-consecutive to $\mathfrak{s}_1$. A $\de$-consecutive cycle $\mathfrak{s}_1\cdots\mathfrak{s}_m$ is called $(L,\theta)$-tame, if the length of each $\mathfrak{s}_i$ is greater than $2L$ and the bending angle between $\mathfrak{s}_i$ and $\mathfrak{s}_{i+1}$ is less than $\theta$.

\item For a $\de$-consecutive chain $\mathfrak{s}_1\cdots\mathfrak{s}_m$, we define its \textit{reduced\ concatenation}, denoted by
\begin{equation*}
    \langle\mathfrak{s}_1\mathfrak{s}_2\cdots\mathfrak{s}_m\rangle,
\end{equation*}
to be the oriented framed segments as follows. The carrier segment of $\langle\mathfrak{s}_1\mathfrak{s}_2\cdots\mathfrak{s}_m\rangle$ is the geodesic arc which is homotopic to the concatenation of the carrier segments of $\mathfrak{s}_i$'s, relative to the initial endpoint of $\mathfrak{s}_1$ and the terminal endpoint of $\mathfrak{s}_m$; the initial framing of $\langle\mathfrak{s}_1\mathfrak{s}_2\cdots\mathfrak{s}_m\rangle$ is the unit normal vector closest to the initial framing of $\mathfrak{s}_1$; the terminal framing of $\langle\mathfrak{s}_1\mathfrak{s}_2\cdots\mathfrak{s}_m\rangle$ is the unit normal vector closet to the terminal framing of $\mathfrak{s}_m$.

\item For a $\de$-consecutive cycle $\mathfrak{s}_1\cdots\mathfrak{s}_m$, we define its \textit{reduced\ cyclic\ concatenation}, denoted by
\begin{equation*}
    [\mathfrak{s}_1\mathfrak{s}_2\cdots\mathfrak{s}_m],
\end{equation*}
to be the unframed oriented closed geodesic curve freely homotopic to the concatenation of the carrier segments of each $\mathfrak{s}_i$, assuming the result nontrivial. 

\item A \textit{continuous} chain of oriented framed segments is a consecutive chain $\mathfrak{s}_1\cdots\mathfrak{s}_m$ with all bending angle $\pi/2$ and
\begin{equation*}
    \vec{n}_{ini}(\mathfrak{s}_{i+1})=\vec{n}_{ter}(\mathfrak{s}_{i})=\vec{t}_{ter}(\mathfrak{s}_i)\times\vec{t}_{ini}(\mathfrak{s}_{i+1}),
\end{equation*}
for $i=1,\cdots,m-1$. It is a \textit{continuous} cycle if furthermore $\mathfrak{s}_m$ is consecutive to $\mathfrak{s}_1$ with bending angle $\pi/2$ and
\begin{equation*}
    \vec{n}_{ini}(\mathfrak{s}_{1})=\vec{n}_{ter}(\mathfrak{s}_{m})=\vec{t}_{ter}(\mathfrak{s}_m)\times\vec{t}_{ini}(\mathfrak{s}_{1}).
\end{equation*}

\item A continuous cycle $\mathfrak{s}_1\cdots\mathfrak{s}_{2m}$ is called $(L,d,\De)$-\textit{tame}, if the length of each $\mathfrak{s}_{2i-1}$ is greater than $2L$, the length of each $\mathfrak{s}_{2i}$ is no greater than $2d$ (which is allowed to be 0) and $|\ln|\sinh(\lle(\mathfrak{s}_{2i})/2)||$ is bounded by $\De$, here the geometric meaning of $|\ln|\sinh(\lle(\mathfrak{s}_{2i})/2)||$ will be explained later.

\end{enumerate}
\end{defn}

\subsection{Inefficiency of framed segments} 

We first recall that for any bending angle $\theta$, we have the \textit{inefficiency} of $\theta$ defined as 
\begin{equation*}
    I(\theta)=2\ln(\sec(\theta/2)).
\end{equation*}
The next lemma interpret the geometric meaning of the inefficiency of angles.

\begin{lem}[Lemma 4.10 in \cite{LM15}]
For any $0<\theta<\pi$ and $\ep>0$, there exists $L>0$ such that the following holds. Suppose that $\De ABC$ is a geodesic triangle in $\HH^3$ with $|CA|,|CB|>L$ and $\angle C=\pi-\theta$, then
\begin{enumerate}

\item $\angle A+\angle B<\ep$;

\item $I(\theta)-\ep<|CA|+|CB|-|AB|<I(\theta)$.
\end{enumerate}
\end{lem}

Now we want to define the inefficiency of a framed segment with complex length $d\in\R_{\geqslant0}\times\R/2\pi\Z$, which can be viewed as a subset of $\C/2\pi i\Z$, by the following lemma. 

\begin{lem}\label{iod}
For any $\De >0$ and $\ep>0$, there exists $L>0$ such that the following holds. Suppose $\mathfrak{s}_1\mathfrak{s}_2\mathfrak{s}_3$ is a continuous chain of framed segments satisfying $\lle(\mathfrak{s}_2)=2d$ with $d\neq0\in\R_{\geqslant0}\times\R/2\pi\Z$ and $|\ln|\sinh(d)||<\De$, $l(\mathfrak{s}_1),l(\mathfrak{s}_3)>2L$. Then
\begin{enumerate}
\item $\angle(\mathfrak{s}_3,\overline{\langle\mathfrak{s}_1\mathfrak{s}_2\mathfrak{s}_3\rangle}),\angle(\overline{\langle\mathfrak{s}_1\mathfrak{s}_2\mathfrak{s}_3\rangle},\mathfrak{s}_1)>\pi-\ep$;

\item $\left|l(\mathfrak{s}_1)+l(\mathfrak{s}_2)+l(\mathfrak{s}_3)-l(\langle\mathfrak{s}_1\mathfrak{s}_2\mathfrak{s}_3\rangle)
    -(\re(2d)-2\ln\left|\sinh\left(d\right)\right|)\right|<\ep$;

\item $\left|\phi(\mathfrak{s}_1)+\phi(\mathfrak{s}_2)+\phi(\mathfrak{s}_3)-\phi(\langle\mathfrak{s}_1\mathfrak{s}_2\mathfrak{s}_3\rangle)-(\im(2d)-2\Arg\left(\sinh\left(d\right)\right))\right|<\ep$,
\end{enumerate}
where $|\cdot|$ on $\R/2\pi\Z$ is understood as the distance from 0 valued in $[0,\pi]$ and $\Arg$ is the principle value of the argument.
\end{lem}

\begin{proof}
Let $\mathfrak{s}_4=\overline{\langle\mathfrak{s}_1\mathfrak{s}_2\mathfrak{s}_3\rangle}$, so $\mathfrak{s}_4$ and $\langle\mathfrak{s}_1\mathfrak{s}_2\mathfrak{s}_3\rangle$ have the same length and phase. Let $A_i$ be the joint point of $\mathfrak{s}_{i-1}$ and $\mathfrak{s}_{i}$ for $i\in\Z/4\Z$, so $A_1A_2A_3A_4$ is a hyperbolic quadrilateral with two right angles, which also can be regarded as a degenerate right-angles hexagon. Let $\vec{n}_1=\vec{t}_{ter}(\mathfrak{s}_4)\times\vec{t}_{ini}(\mathfrak{s}_1)$ and $\vec{n}_4=\vec{t}_{ter}(\mathfrak{s}_3)\times\vec{t}_{ini}(\mathfrak{s}_4)$ be the common normal vectors at $A_1$ and $A_4$ respectively, and then we can define the complex length of $A_iA_{i+1}$, denoted by $l_i$, in this right-angles hexagon, where $l_2=2d$. By the hyperbolic Cosine Law for right-angled hexagon, we have
\begin{equation*}
    \cosh(2d)=\frac{\cosh(l_1)\cosh(l_3)+\cosh(l_4)}{\sinh(l_1)\sinh(l_2)}.
\end{equation*}
Therefore
\begin{equation*}
    \frac{\cosh(l_4)}{e^{(l_1+l_3)}}=\cosh(2d)\frac{\sinh(l_1)}{e^{l_1}}\frac{\sinh(l_3)}{e^{l_3}}-\frac{\cosh(l_1)}{e^{l_1}}\frac{\cosh(l_3)}{e^{l_3}}.
\end{equation*}
When $\re(l_1),\re(l_3)\to+\infty$,
\begin{equation*}
    \frac{\cosh(l_4)}{e^{l_1+l_3}}\to \frac{1}{4}(\cosh(2d)-1).
\end{equation*}
Since $d\neq0$, we have $\re(l_4)\to+\infty$. So when $\re(l_1),\re(l_3)\to+\infty$
\begin{equation}\label{cl}
    e^{l_4-l_1-l_3}\to\sinh^2\left(d\right),
\end{equation}
where $|\sinh(d)|$ is bounded above and below. By $\re(l_i)=l(\mathfrak{s}_i)$ and (\ref{cl}), we know there exists $L_1>0$ such that when $l(\mathfrak{s}_1),l(\mathfrak{s}_3)>L_1$,
\begin{equation}\label{liod}
    \left|l(\mathfrak{s}_1)+l(\mathfrak{s}_3)-l(\mathfrak{s}_4)
    +2\ln\left|\sinh\left(d\right)\right|\right|
    <\ep.
\end{equation}
By $l(\mathfrak{s}_2)=\re(d)$ and (\ref{liod}), we get the inequality (2). 

By the hyperbolic Cosine Law again, we have
\begin{equation*}
    \cos(\angle A_1)=\frac{\cosh(l_1)\cosh(l_4)+\cosh(l_3)}{\sinh(l_1)\sinh(l_4)}.
\end{equation*}
When $\re(l_1),\re(l_3)\to+\infty$, by (\ref{cl})
\begin{equation*}
    \cos(\angle A_1)\to\frac{e^{l_1+l_4}+e^{l_3}}{e^{l_1+l_4}}\to 1.
\end{equation*}
Therefore there exists $L_2>0$ such that when $l(\mathfrak{s}_1),l(\mathfrak{s}_3)>L_2$, $\angle A_1<\ep$. Then by $\angle(\vec{t}_{ini}(\mathfrak{s}_1),\vec{t}_{ter}(\langle\mathfrak{s}_1\mathfrak{s}_2\mathfrak{s}_3\rangle))=\pi-\angle A_1$,
\begin{equation*}
    \angle(\vec{t}_{ini}(\mathfrak{s}_1),\vec{t}_{ter}(\langle\mathfrak{s}_1\mathfrak{s}_2\mathfrak{s}_3\rangle))>\pi-\ep.
\end{equation*}
And the second part of (1) is true for the similar reason.

For (3), we notice that when $\re(l_1),\re(l_3)>L_2$, $\angle A_1,\angle A_4<\ep$. Thus $\mathfrak{s}_3$ is $\ep$-consecutive to $\mathfrak{s}_4$ and $\mathfrak{s}_4$ is $\ep$-consecutive to $\mathfrak{s}_1$. Hence we have
\begin{equation}\label{piod}
    |(\phi(\mathfrak{s}_1)+\phi(\mathfrak{s}_3)-\phi(\mathfrak{s}_4))-(\im(l_1)+\im(l_3)-\im(l_4))|<2\ep.
\end{equation}
By (\ref{cl}), we know there exists $L_3>0$, such that when $l(\mathfrak{s}_1),l(\mathfrak{s}_3)>L_3$,
\begin{equation*}
    |\im(l_1)+\im(l_3)-\im(l_4)+2\Arg(\sinh(d))|<\ep.
\end{equation*}
Together with (\ref{piod}), we have
\begin{align*}
    |(\phi(\mathfrak{s}_1)+\phi(\mathfrak{s}_3)-\phi(\mathfrak{s}_4))+2\Arg(\sinh(d))|<3\ep.
\end{align*}
Hence the proof is completed since $\phi(\mathfrak{s}_2)=\im(d)$.
\end{proof}

We define the length inefficiency of $d$ as
\begin{equation*}
    I_l(d)=\re(2d)-2\ln\left|\sinh\left(d\right)\right|
\end{equation*}
and the phase inefficiency of $d$ as 
\begin{equation*}
    I_\phi(d)
    =\im(2d)-2\Arg\left(\sinh\left(d\right)\right),
\end{equation*}
and the geometric meaning is explained by the above lemma.

When $d=i\theta$ is pure imaginary, we have
\begin{equation*}
    I_l(i\theta)=-2\ln|\sinh(i\theta)|=2\ln\csc(\theta),
\end{equation*}
which implies $I_l(i\theta)$ coincides with the inefficiency of angle $(\pi-2\theta)$ since the bending angle is defined as the exterior angle.

\subsection{Inefficiency of framed segment cycles}

For sufficiently tame approximately consecutive framed segment cycles, Liu and Markovic estimated its inefficiency in sense of length and phase. We restate the results as follows.

\begin{lem}[Lemma 4.8 in \cite{LM15}]\label{Ltheta}
Given any $\delta\geq0$, $\pi>\theta>0$ and $\ep>0$, there exists $L>0$ such that the following holds. If $\mathfrak{s}_1\cdots\mathfrak{s}_m$ is an $(L,\theta)$-tame $\de$-consecutive cycle of oriented framed segments, let $\theta_i\in[0,\pi-\theta]$ be the bending angle between $\mathfrak{s}_i$ and $\mathfrak{s}_{i+1}$ with $\mathfrak{s}_{m+1}$ equal to $\mathfrak{s}_1$, then
\begin{equation*}
    |l([\mathfrak{s}_1\cdots\mathfrak{s}_m])-\sum\limits_{i=1}^{m}l(\mathfrak{s}_i)+\sum\limits_{i=1}^{m}I(\theta_i)|<\ep,
\end{equation*}
and
\begin{equation*}
    |\phi([\mathfrak{s}_1\cdots\mathfrak{s}_m])-\sum\limits_{i=1}^{m}\phi(\mathfrak{s}_i)|<m\de+\ep,
\end{equation*}
where $|\cdot|$ on $\R/2\pi\Z$ is understood as the distance from 0 valued in $[0,\pi]$.
\end{lem}

As a matter of fact, the above lemma, which works for $(L,\theta)$-tame $\de$-consecutive cycles, can be generalized to the following lemma for $(L,d,\theta)$-tame continuous cycles.

\begin{lem}[Sum of Inefficiencies Lemma]\label{LDe}
Given any $m\in\Z_+$, $\De>0$, $d>0$ and $1/4>\ep>0$, there exists $L>0$ such that the following holds. If $\mathfrak{s}_1\cdots\mathfrak{s}_{2m}$ is an $(L,d,\De)$-tame continuous cycle of oriented framed segments, let $2d_i$ be the complex length of $\mathfrak{s}_{2i}$, then
\begin{equation}\label{dl}
    |l([\mathfrak{s}_1\cdots\mathfrak{s}_{2m}])-\sum\limits_{i=1}^{2m}l(\mathfrak{s}_{i})+\sum\limits_{i=1}^{m}I_l(d_i)|<\ep,
\end{equation}
and
\begin{equation}\label{dp}
    |\phi([\mathfrak{s}_1\cdots\mathfrak{s}_{2m}])-\sum\limits_{i=1}^{2m}\phi(\mathfrak{s}_{i})+\sum\limits_{i=1}^{m}I_\phi(d_i)|<\ep,
\end{equation}
where $|\cdot|$ on $\R/2\pi\Z$ is understood as the distance from 0 valued in $[0,\pi]$.
\end{lem}

\begin{proof}
Take $\ep_1=\ep/(2m+1)$. We write
\begin{equation*}
    \De l(\mathfrak{s}_1\cdots\mathfrak{s}_k)=\left\{
    \begin{aligned}
    \sum\limits_{i=1}^{k}l(\mathfrak{s}_i)-l(\langle\mathfrak{s}_1\cdots\mathfrak{s}_k\rangle),&\ \mathrm{if}\ \mathfrak{s}_1,\cdots,\mathfrak{s}_k\ \mathrm{is\ a\ chain;}\\
    \sum\limits_{i=1}^{k}l(\mathfrak{s}_i)-l([\mathfrak{s}_1\cdots\mathfrak{s}_k]),&\ \mathrm{if}\ \mathfrak{s}_1,\cdots,\mathfrak{s}_k\ \mathrm{is\ a\ cycle.}
    \end{aligned}\right.
\end{equation*}
and
\begin{equation*}
    \De \phi(\mathfrak{s}_1\cdots\mathfrak{s}_k)=\left\{
    \begin{aligned}
    \sum\limits_{i=1}^{k}\phi(\mathfrak{s}_i)-\phi(\langle\mathfrak{s}_1\cdots\mathfrak{s}_k\rangle),&\ \mathrm{if}\ \mathfrak{s}_1,\cdots,\mathfrak{s}_k\ \mathrm{is\ a\ chain;}\\
    \sum\limits_{i=1}^{k}\phi(\mathfrak{s}_i)-\phi([\mathfrak{s}_1\cdots\mathfrak{s}_k]),&\ \mathrm{if}\ \mathfrak{s}_1,\cdots,\mathfrak{s}_k\ \mathrm{is\ a\ cycle.}
    \end{aligned}\right.
\end{equation*}
for convinience.

For each $1\leq i\leq m$, write $\mathfrak{s}_{2i-1}$ as the concatenation of two consecutive oriented framed segments $\mathfrak{s}_{2i-1}^-$ and $\mathfrak{s}_{2i-1}^+$ of equal length and phase. Let $\tilde{\mathfrak{s}}_{i}=\langle\mathfrak{s}_{2i-1}^+\mathfrak{s}_{2i}\mathfrak{s}_{2i+1}^-\rangle$ with $\mathfrak{s}_{2m+1}=\mathfrak{s}_1$. Since the phase of each $\mathfrak{s}_{2i}$ is at least $\theta$ away from $0$, then by Lemma \ref{iod}, there exists $L_1>0$ such that
\begin{equation}\label{dl1}
    |\De l(\mathfrak{s}_{2i-1}^+\mathfrak{s}_{2i}\mathfrak{s}_{2i+1}^-)-I_l(d_i)|<\ep_1
\end{equation}
and
\begin{equation}\label{dp1}
    |\De \phi(\mathfrak{s}_{2i-1}^+\mathfrak{s}_{2i}\mathfrak{s}_{2i+1}^-)-I_\phi(d_i)|<\ep_1,
\end{equation}
for $i=1,\cdots,m$, and $\tilde{\mathfrak{s}}_i$ is $\ep_1$-consecutive to $\tilde{\mathfrak{s}}_{i+1}$. 

Now $\tilde{\mathfrak{s}}_1\cdots\tilde{\mathfrak{s}}_m$ is a $\ep_1$-consecutive cycle of framed segments. Let $\theta_i=\angle(\tilde{\mathfrak{s}}_i,\tilde{\mathfrak{s}}_{i+1})$ be the bending angle between $\tilde{\mathfrak{s}}_i$ and $\tilde{\mathfrak{s}}_{i+1}$. By Lemma \ref{iod}, the unsigned angle between $\vec{t}_{ini}(\tilde{\mathfrak{s}}_i)$ and $\vec{t}_{ini}(\mathfrak{s}_{2i-1}^+)$ is less than $\ep_1$, and the same for the unsigned angle between $\vec{t}_{ter}(\tilde{\mathfrak{s}}_{i-1})$ and $\vec{t}_{ter}(\mathfrak{s}_{2i-1}^-)$. Thus $\theta_i$ is less than $2\ep_1$. By Lemma \ref{Ltheta}, there exists $L_2>0$ such that when $l(\mathfrak{s}_i)>L_2$, we have
\begin{equation}\label{dl2}
    |\De l(\tilde{\mathfrak{s}}_1\cdots\tilde{\mathfrak{s}}_m)-\sum\limits_{i=1}^{m}I(\theta_i)|<\ep_1
\end{equation}
and
\begin{equation}\label{dp2}
    |\De \phi(\tilde{\mathfrak{s}}_1\cdots\tilde{\mathfrak{s}}_m)|<m\ep_1+\ep_1.
\end{equation}
Since $\ep_1<\ep<1/4$, we have
\begin{equation*}
    I(\theta_i)<I(2\ep_1)=2\ln(\sec(\ep_1))<\ep_1.
\end{equation*}
Together (\ref{dl2}), we get
\begin{equation}\label{dl3}
    |\De l(\tilde{\mathfrak{s}}_1\cdots\tilde{\mathfrak{s}}_m)|<(m+1)\ep_1.
\end{equation}
By (\ref{dl1}) and (\ref{dl3}),
\begin{equation*}
    \begin{aligned}
    |\De l(\mathfrak{s}_1\cdots\mathfrak{s}_{2m})-\sum\limits_{i=1}^{m}I_l(d_i)|
    &=|\sum\limits_{i=1}^{m}\left(\De l(\mathfrak{s}_{2i-1}^+\mathfrak{s}_{2i}\mathfrak{s}_{2i+1}^-)-I_l(d_i)\right)
    +\De l(\tilde{\mathfrak{s}}_1\cdots\tilde{\mathfrak{s}}_m)|\\
    &\leq (2m+1)\ep_1=\ep,
    \end{aligned}
\end{equation*}
which proves (\ref{dl}). By (\ref{dp1}) and (\ref{dp2}), we have
\begin{equation*}
    \begin{aligned}
    |\De \phi(\mathfrak{s}_1\cdots\mathfrak{s}_{2m})-\sum\limits_{i=1}^{m}I_\phi(d_i)|
    &=|\sum\limits_{i=1}^{m}\left(\De \phi(\mathfrak{s}_{2i-1}^+\mathfrak{s}_{2i}\mathfrak{s}_{2i+1}^-)-I_\phi(d_i)\right)
    +\De \phi(\tilde{\mathfrak{s}}_1\cdots\tilde{\mathfrak{s}}_m)|\\
    &\leq (2m+1)\ep_1=\ep,
    \end{aligned}
\end{equation*}
which proves (\ref{dp}).
\end{proof}

\subsection{Zigzag geodeiscs}

Now we want to study the cycles of 4 framed segments, which will appear several times in subsequent sections.

\begin{defn}
For $L,\ep>0$, a continuous cycle of 4 framed segments $\mathfrak{s}_1\mathfrak{s}_2\mathfrak{s}_3\mathfrak{s}_{4}$ is called $(L,\ep)$-\textit{zigzag}, if the followings are satisfied:
\begin{enumerate}
\item $l(\mathfrak{s}_i)>2L$, for $i=1,3$;

\item $|\lle(\mathfrak{s}_i)-i\pi|<2\ep$, for $i=2,4$.
\end{enumerate}
\end{defn}

In Lemma \ref{LDe}, we estimated the length of the reduced cyclic concatenation of tame enough framed segment cycles. It turns out that we can also detect its location if furthermore it is zigzag, which is also the last lemma of this section.

\begin{lem}\label{zglem}
For any $\de>0$, there exists constants $L_0,\ep_0$ such that the following holds. Suppose $M$ is an oriented hyperbolic 3-manifold and $\mathfrak{s}_1\mathfrak{s}_2\mathfrak{s}_3\mathfrak{s}_{4}$ is $(L,\ep)$-zigzag in $M$, for some $L>L_0$ and $0<\ep<\ep_0$. Then the Hausdorff distance between the carriers of $\mathfrak{s}_1\mathfrak{s}_2\mathfrak{s}_3\mathfrak{s}_{4}$ and $[\mathfrak{s}_1\mathfrak{s}_2\mathfrak{s}_3\mathfrak{s}_{4}]$ is less than $\de$.
\end{lem}

\begin{proof}
For convenience, let $\mathfrak{s}=\mathfrak{s}_1\mathfrak{s}_2\mathfrak{s}_3\mathfrak{s}_{4}$, $l_{2i-1}=\hl(\mathfrak{s}_{2i-1})$ and $l_{2i}=\hl(\mathfrak{s}_{2i})-\frac{i\pi}{2}$ for $i=1,2$, with $|l_{2i}|<\ep$. Let $A_i$ be the intersection of $\mathfrak{s}_i$ and $\mathfrak{s}_{i+1}$ for $i=1,2,3,4$ where $\mathfrak{s}_5=\mathfrak{s}_1$. Without loss of generality, we assume that $p_0=(0,1)\in\HH^3$ in the upper half space model, also denoted by $\tilde{A}_1$, is a lift of $A_1$, the initial direction $\vec{t}_{ini}(\mathfrak{s}_1)$ at $A_1$ lifts to $\mathbf{u}_0=(0,1)\in T_{p_0}\HH^3$ and the initial framing $\vec{n}_{ini}(\mathfrak{s}_1)$ lifts to $\mathbf{v}_0=(1,0)\in T_{p_0}\HH^3$. Let $\tilde{\mathfrak{s}}$ be the lift of $\mathfrak{s}$ which passes through $\tilde{A}_1$, and now we want to prove that there is a lift of $[\mathfrak{s}]$ which is close to $\tilde{\mathfrak{s}}$.

Since $[\mathfrak{s}]$ is the geodesic representative of $\mathfrak{s}$, we know that the right action along $\mathfrak{s}$ is conjugate to the right action along $[\mathfrak{s}]$ in $\PSL(2,\C)$. That is,
\begin{align}\label{zigzag}
    X=A\left(l_1\right)B\left(\frac{\pi}{2}\right)A\left(l_2+\frac{i\pi}{2}\right)B\left(\frac{\pi}{2}\right)A(l_3)B\left(\frac{\pi}{2}\right)A\left(l_4+\frac{i\pi}{2}\right)B\left(\frac{\pi}{2}\right)\sim A(l) \in \PSL(2,\C),
\end{align}
where $l=\hl([\mathfrak{s}])$. Let $Y\in\PSL(2,\C)$ such that $X=YA(l)Y^{-1}$, then $\Psi_0\cdot Y$ is a frame on a lift of $[\mathfrak{s}]$ with $\mathbf{u}_0\cdot Y=Y\cdot\mathbf{u}_0$ a tangent vector. Since the right action and the left action have the same results when acting on $\Psi_0$, hence $Y\cdot \tilde{A}_1$ is a point on a lift of $[\mathfrak{s}]$. Moreover, the collection of all $Y\cdot \tilde{A}_1$ for $Y$ conjugating $X$ to $A(l)$ is a lift of $[\mathfrak{s}]$, and this lift is also the fixed geodesic of $X$ as left action in $\HH^3$. We denote this lift by $[\tilde{\mathfrak{s}}]$, then the distance from $A_1$ to $\mathfrak{s}$ is exactly the distance from $\tilde{A}_1$ to $[\tilde{\mathfrak{s}}]$. 

The right action along $\mathfrak{s}$ is given by
\begin{align}
\begin{split}
    X &= A\left(l_1\right)B\left(\frac{\pi}{2}\right)A\left(l_2+\frac{i\pi}{2}\right)B\left(\frac{\pi}{2}\right)A(l_3)B\left(\frac{\pi}{2}\right)A\left(l_4+\frac{i\pi}{2}\right)B\left(\frac{\pi}{2}\right)\\
    &= \begin{pmatrix} e^{l_1} & 0 \\ 0 & e^{-l_1}\end{pmatrix}\begin{pmatrix} \sqrt{2}/2 & \sqrt{2}/2 \\ -\sqrt{2}/2 & \sqrt{2}/2 \end{pmatrix} \begin{pmatrix} ie^{l_2} & 0 \\ 0 & -ie^{-l_2}\end{pmatrix}\begin{pmatrix} \sqrt{2}/2 & \sqrt{2}/2 \\ -\sqrt{2}/2 & \sqrt{2}/2 \end{pmatrix}\\
    &\ \ \ \begin{pmatrix} e^{l_3} & 0 \\ 0 & e^{-l_3}\end{pmatrix}\begin{pmatrix} \sqrt{2}/2 & \sqrt{2}/2 \\ -\sqrt{2}/2 & \sqrt{2}/2 \end{pmatrix}\begin{pmatrix} ie^{l_4} & 0 \\ 0 & -ie^{-l_4}\end{pmatrix}\begin{pmatrix} \sqrt{2}/2 & \sqrt{2}/2 \\ -\sqrt{2}/2 & \sqrt{2}/2 \end{pmatrix}\\
    &= -\begin{pmatrix} e^{l_1} & 0 \\ 0 & e^{-l_1}\end{pmatrix}\begin{pmatrix} \cosh{(l_2)} & \sinh{(l_2)} \\ -\sinh{(l_2)} & -\cosh{(l_2)} \end{pmatrix}\begin{pmatrix} e^{l_3} & 0 \\ 0 & e^{-l_3}\end{pmatrix}\begin{pmatrix} \cosh{(l_4)} & \sinh{(l_4)} \\ -\sinh{(l_4)} & -\cosh{(l_4)} \end{pmatrix}\\
    &\xlongequal{\PSL(2,\C)} \begin{pmatrix} e^{l_1} & 0 \\ 0 & e^{-l_1}\end{pmatrix}\begin{pmatrix} \cosh{(l_2)} & \sinh{(l_2)} \\ -\sinh{(l_2)} & -\cosh{(l_2)} \end{pmatrix}\begin{pmatrix} e^{l_3} & 0 \\ 0 & e^{-l_3}\end{pmatrix}\begin{pmatrix} \cosh{(l_4)} & \sinh{(l_4)} \\ -\sinh{(l_4)} & -\cosh{(l_4)} \end{pmatrix}\\
    &=\begin{pmatrix} a & b \\
    c & d \end{pmatrix},
\end{split}
\end{align}
where
\begin{equation}\label{abcd}
    \begin{aligned}
    a(l_1,l_2,l_3,l_4)&=e^{l_1}(e^{l_3}\cosh({l_2})\cosh({l_4})-e^{-l_3}\sinh({l_2})\sinh({l_4})),\\
    b(l_1,l_2,l_3,l_4)&=e^{l_1}(e^{l_3}\cosh({l_2})\sinh({l_4})-e^{-l_3}\sinh({l_2})\cosh({l_4})),\\
    c(l_1,l_2,l_3,l_4)&=e^{-l_1}(e^{-l_3}\cosh({l_2})\sinh({l_4})-e^{l_3}\sinh({l_2})\cosh({l_4})),\\
    d(l_1,l_2,l_3,l_4)&=e^{-l_1}(e^{-l_3}\cosh({l_2})\cosh({l_4})-e^{l_3}\sinh({l_2})\sinh({l_4}))
    \end{aligned}
\end{equation}
are functions of $l_1,l_2,l_3$ and $l_4$. Let $\lambda=e^l$, and then
\begin{equation}\label{z1z2}
    z_1=\frac{b}{\lambda-a}=\frac{\lambda-d}{c}, z_2=\frac{b\lambda}{1-a\lambda}=\frac{1-d\lambda}{c\lambda}
\end{equation}
are the two fixed points of $X$, as left action, on $\C=\partial\HH^3$. Hence we want to determine the distance $D$ from $\tilde{A}_1=(0,1)$ to the geodesic connecting $z_1$ and $z_2$, which is $[\tilde{\gamma}]$. Consider M\"obius transformation
\begin{equation*}
    T:z\mapsto\frac{z-z_1}{z-z_2}
\end{equation*}
which sends $[\mathfrak{s}]$ to the positive $t$-axis, and we have
\begin{equation*}
    T((0,1))=\left(\frac{1+z_1\overline{z_2}}{1+|z_2|^2},\frac{|z_1-z_2|}{1+|z_2|^2}\right).
\end{equation*}
Hence $D$ is also the distance from $T((0,1))$ to the positive $t$-axis, and it is given by
\begin{equation*}
    \begin{aligned}
    e^D&=\frac{\left|\frac{1+z_1\overline{z_2}}{1+|z_2|^2}\right|
    +\sqrt{\left|\frac{1+z_1\overline{z_2}}{1+|z_2|^2}\right|^2+\left|\frac{|z_1-z_2|}{1+|z_2|^2}\right|^2}}
    {\frac{|z_1-z_2|}{1+|z_2|^2}}\\
    &=\frac{|1+z_1\overline{z_2}|
    +\sqrt{|1+z_1\overline{z_2}|^2+|z_1-z_2|^2}}{|z_1-z_2|}.
    \end{aligned}
\end{equation*}
Therefore we have 
\begin{equation}\label{sinhD}
    \sinh(D)=\left|\frac{1+z_1\overline{z_2}}{z_1-z_2}\right|.
\end{equation}
By (\ref{z1z2}),
\begin{equation*}
    \begin{aligned}
    |z_1-z_2|&=\left|\frac{\lambda-d}{c}-\frac{1-d\lambda}{c\lambda}\right|=\left|\frac{\lambda^2-1}{c\lambda}\right|
    =\left|\frac{1}{c}\left(\lambda-\frac{1}{\lambda}\right)\right|\\
    &=\left|\frac{\sqrt{(a+d)^2-4}}{c}\right|=\frac{\sqrt{|(a+d)^2-4|}}{|c|}
    \end{aligned}
\end{equation*}
and
\begin{equation*}
    \begin{aligned}
    |1-z_1\overline{z_2}|^2&=(1-z_1\overline{z_2})(1-\overline{z_1}z_2)=1+|z_1z_2|^2-z_1\overline{z_2}-\overline{z_1}z_2\\
    &=1+\left|\frac{(\lambda-d)(1-d\lambda)}{c^2\lambda}\right|^2
    -\frac{\lambda-d}{c}\cdot\frac{1-\overline{d}\overline{\lambda}}{\overline{c}\overline{\lambda}}
    -\frac{\overline{\lambda}-\overline{d}}{\overline{c}}\cdot\frac{1-d\lambda}{c\lambda}\\
    &=1+\left|\frac{(1-ad)\lambda}{c^2\lambda}\right|
    -\frac{1}{|c|^2}\left(
    \frac{\lambda}{\overline{\lambda}}+\frac{\overline{\lambda}}{\lambda}+2|d|^2-d\left(\overline{\lambda}+\frac{1}{\overline{\lambda}}\right)-\overline{d}\left(\lambda+\frac{1}{\lambda}\right)
    \right)\\
    &=1+\frac{|b|^2}{|c|^2}+\frac{1}{|c|^2}\left(a\overline{d}+\overline{a}d-\frac{\lambda}{\overline{\lambda}}-\frac{\overline{\lambda}}{\lambda}\right)\\
    &\leq 1+\frac{|b|^2}{|c|^2}+\frac{2|ad|+2}{|c|^2}
    =1+\frac{|b|^2}{|c|^2}+\frac{2|bc+1|+2}{|c|^2}
    \leq\frac{(|b|+|c|)^2+4}{|c|^2}\\
    &\leq\left|\frac{|b|+|c|+2}{c}\right|^2.
    \end{aligned}
\end{equation*}
Together with (\ref{sinhD}) we have
\begin{equation}\label{sinhD2}
    \sinh(D)\leq\frac{|b|+|c|+2}{\sqrt{|(a+d)^2-4|}}.
\end{equation}
Take $L_1=10^{10}$ and $\ep_1=10^{-10}$, then by (\ref{abcd}), $L>L_1$ and $0<\ep<\ep_1$, we have
\begin{equation}\label{bc}
    \begin{aligned}
    |b|+|c|+2=&|e^{l_1}(e^{l_3}\cosh({l_2})\sinh({l_4})-e^{-l_3}\sinh({l_2})\cosh({l_4}))|\\
    &+|e^{-l_1}(e^{-l_3}\cosh({l_2})\sinh({l_4})-e^{l_3}\sinh({l_2})\cosh({l_4}))|+2\\
    \leq& 4(3\ep)e^{\re(l_1+l_3)}+2
    \end{aligned}
\end{equation}
and
\begin{equation}
    \begin{aligned}\label{ad}
    |(a+d)^2-4|=&|(\cosh(l_2)\cosh(l_4)\cosh(l_1+l_3)+\sinh(l_2)\sinh(l_4)\cosh(l_1-l_3))^2-4|\\
    \geq& (|\cosh(l_2)\cosh(l_4)\cosh(l_1+l_3)|-|\sinh(l_2)\sinh(l_4)\cosh(l_1-l_3)|)^2-4\\
    \geq& \frac{1}{4}e^{\re(l_1+l_3)}.
    \end{aligned}
\end{equation}
Thus by (\ref{sinhD2}), (\ref{bc}) and (\ref{ad}),
\begin{equation}\label{sinhD3}
    \sinh(D)\leq\frac{{12\ep}e^{\re(l_1+l_3)}+2}{\frac{1}{4}e^{\re(l_1+l_3)}}={48\ep}+8e^{-\re(l_1+l_3)}\leq {48\ep}+8e^{-2L}.
\end{equation}
Hence for any $\de>0$, we can find $L_2,\ep_2$ such that when $L>L_2$ and $0<\ep<\ep_2$, we have
\begin{equation}\label{de}
    {48\ep}+8e^{-2L}<\de.
\end{equation}
So let $L_0=\max\{L_1,L_2\}$ and $\ep_0=\min\{\ep_1,\ep_2\}$, then by (\ref{sinhD3}) and (\ref{de}), we get $D<\sinh(D)<\de$ when $L>L_0$ and $0<\ep<\ep_0$.

Similarly, we have the same result for $\tilde{A}_3$. For $\tilde{A}_2$ and $\tilde{A}_4$, we only need to reverse the orientation of each $\mathfrak{s}_i$ and use the same technique for $\overline{\mathfrak{s}}=(\overline{\mathfrak{s}_4})(\overline{\mathfrak{s}_3})(\overline{\mathfrak{s}_2})(\overline{\mathfrak{s}_1})$. And then by the property of distance function between geodesics in $\HH^3$, we know that for each point $\tilde{A}$ on $\tilde{\mathfrak{s}}$, the distance from $\tilde{A}$ to $[\tilde{\mathfrak{s}}]$ is less than $\de$.

On the other hand, we consider the projection map $\rho$ from $\tilde{\mathfrak{s}}$ to $[\tilde{\mathfrak{s}}]$ by sending each point on $\tilde{\mathfrak{s}}$ to its closest point on $[\tilde{\mathfrak{s}}]$. Then $\rho$ is a continuous map, so it is onto since $\tilde{\mathfrak{s}}$ and $[\tilde{\mathfrak{s}}]$ have same end point on boundary at infinity of $\HH^3$. For each $\tilde{B}$ on $[\tilde{\mathfrak{s}}]$, let $\tilde{B}'$ be a preimage of $\tilde{B}$ under $\rho$. Since $\tilde{B}$ is the closest point of $\tilde{B}'$ on $[\tilde{\mathfrak{s}}]$, thus $d(\tilde{B},\tilde{B}')<\de$. Hence the distance from $\tilde{B}$ to $\tilde{\mathfrak{s}}$ is less than $\de$.

To conclude, when $L>L_0$ and $0<\ep<\ep_0$, the Hausdorff distance between $\tilde{\mathfrak{s}}$ and $[\tilde{\mathfrak{s}}]$ is less than $\de$ in $\HH^3$, so the same result holds for $\mathfrak{s}$ and $[\mathfrak{s}]$ in $M$.

\end{proof}

Actually by (\ref{de}), we can have more accurate estimate on $\delta$ from the above lemma in some situations, and we state it as the following lemma.

\begin{lem}

For any $C>0$, there exists constants $R_0$ and $B>0$ such that the following holds. Suppose $M$ is an oriented hyperbolic 3-manifold and $\mathfrak{s}_1\mathfrak{s}_2\mathfrak{s}_3\mathfrak{s}_{4}$ is $(R/2,C e^{-R/2})$-zigzag in $M$, for some $R>R_0$. Then the Hausdorff distance between the carriers of $\mathfrak{s}_1\mathfrak{s}_2\mathfrak{s}_3\mathfrak{s}_{4}$ and $[\mathfrak{s}_1\mathfrak{s}_2\mathfrak{s}_3\mathfrak{s}_{4}]$ is less than $B e^{-R/2}$.

\end{lem}

\section{Pants decomposition of genus-2 quasi-Fuchsian groups}\label{S5}

In this section, we will introduce how to construct a pants decomposition of a genus-2 quasi-Fuchisian group with long cuffs and bounded shears. From now on, assume $\Gamma$ is a genus-2 $K$-quasi-Fuchsian group for some given $K>1$, and $CM(\Gamma)$ be the convex core of $\Gamma$ which is homotopy equivalent to a topological genus-2 oriented surface $S$. Let $f:S\to CM(\Gamma)$ be the homotopy equivalent map. 

\subsection{Good pants decomposition}

Suppose $\{C'_1,C'_2,C'_3\}$, a set of 3 disjoint oriented simple closed curves on $S$, is a non-separating pants decomposition of $S$, then we take the oriented closed geodesic $C_i$ in $CM(\Gamma)$ freely homotopic to $f(C'_i)$, for $i\in\Z/3\Z$. Then $\{C_1,C_2,C_3\}$ is a non-separating pants decomposition of $\Gamma$, so there are two immersed pairs of pants $P_1,P_2$ in $CM(\Gamma)$ whose boundary components are both $\{C'_1,C'_2,C'_3\}$. In each $P_k$, let $\gamma_{i,k}$ be the short orthogonal geodesic between $C_{i+1}$ and $C_{i+2}$, oriented from $C_{i+1}$ to $C_{i+2}$ for $i\in\Z/3\Z$ and $k\in\Z/2\Z$. Then let $\eta_{i,i+1}$ and $\eta_{i,i+2}$ be the orthogonal geodesic between $\gamma_{i,1}$ and $\gamma_{i,2}$ along $C_{i+1}$ and $C_{i+2}$, respectively, and the orientation of $\eta_{i,i+1}$ follows the orientation of $C_{i+1}$. Notice that $\eta_{i+1,i}$ and $\eta_{i-1,i}$ are the short shear between $P_1$ and $P_2$ along $C_{i}$. 

Right now we want to assign frames at endpoints of these oriented geodesic arcs to make them framed segments and use them to construct new pants decompositions of $\Gamma$. For a tuple of positive integers $(n_1,n_2,n_3)\in\Z_{+}^3$, consider the following sets of geodesic arcs, $A_i(n_1,n_2,n_3)=\{\eta_{i,i+1}C_{i+1}^{n_{i+1}},\gamma_{i,2},C_{i-1}^{-n_{i-1}}\eta_{i,i-1}^{-1},\gamma_{i,1}^{-1}\}$, for $i\in\Z/3\Z$. Then the 4 geodesic arcs in each $A_i(n_1,n_2,n_3)$ can form a piecewise geodesic curve, because the concatenation of $\eta_{i,i+1}$ and $C_{i+1}^{n_{i+1}}$ can be regarded as a smooth geodesic. And we can find unique framing at each joint point, such that there are framed segments $\mathfrak{a}_{i}^{(n_{i+1})},\mathfrak{b}_{i},\mathfrak{c}_{i}^{(n_{i-1})}$ and $\mathfrak{d}_{i}$, whose carriers are $\eta_{i,i+1}C_{i+1}^{n_{i+1}},\gamma_{i,2},C_{i-1}^{-n_{i-1}}\eta_{i,i-1}^{-1}$ and $\gamma_{i,1}^{-1}$, forming a zigzag continuous cycle. Let $\mathfrak{s}_{i}^{(n_{i+1},n_{i+2})}=\mathfrak{s}_{i}^{(n_{i+1},n_{i-1})}=\mathfrak{a}_{i}^{(n_{i+1})}\mathfrak{b}_{i}\mathfrak{c}_{i}^{(n_{i-1})}\mathfrak{d}_{i}$. Since we have
\begin{equation*}
    \begin{aligned}
    \lle(\mathfrak{a}_{i}^{(n_{i+1})})=\lle(\eta_{i,i+1})+n_{i+1}\lle(C_{i+1}),\\
    \lle(\mathfrak{c}_{i}^{(n_{i-1})})=\lle(\eta_{i,i-1})+n_{i-1}\lle(C_{i-1}),
    \end{aligned}
\end{equation*}
so by Lemma \ref{LDe}, for any $\ep>0$, there exists $N\in\Z_+$ such that when $n_1,n_2,n_3>N$, we have
\begin{equation*}
    |l([\mathfrak{s}_{i}^{(n_{i+1},n_{i+2})}])-l(\mathfrak{a}_{i}^{(n_{i+1})})-l(\mathfrak{b}_{i})-l(\mathfrak{c}_{i}^{(n_{i-1})})-l(\mathfrak{d}_{i})+I_l(\lle(\mathfrak{b}_{i})/2)+I_l(\lle(\mathfrak{d}_{i})/2)|<\ep
\end{equation*}
and
\begin{equation*}
    |\phi([\mathfrak{s}_{i}^{(n_{i+1},n_{i+2})}])-\phi(\mathfrak{a}_{i}^{(n_{i+1})})-\phi(\mathfrak{b}_{i})-\phi(\mathfrak{c}_{i}^{(n_{i-1})})-\phi(\mathfrak{d}_{i})+I_\phi(\lle(\mathfrak{b}_{i})/2)+I_\phi(\lle(\mathfrak{d}_{i})/2)|<\ep.
\end{equation*}
for $i\in\Z/3\Z$. To simplify the above equations, we let
\begin{equation*}
    \sigma_{i}=\lle(\eta_{i,i+1})+\lle(\eta_{i,i-1})-(I_l(\lle(\mathfrak{b}_{i})/2)+I_l(\lle(\mathfrak{d}_{i})/2))-i(I_\phi(\lle(\mathfrak{b}_{i})/2)+I_\phi(\lle(\mathfrak{d}_{i})/2)).
\end{equation*}
Then when $n_1,n_2,n_3>N$,
\begin{equation}\label{n1n2n3}
    |\lle([\mathfrak{s}_{i}^{(n_{i+1},n_{i+2})}])-n_{i+1}\lle(C_{i+1})-n_{i-1}\lle(C_{i-1})-\sigma_i|<\sqrt{2}\ep,
\end{equation}
for $i\in\Z/3\Z$.

To see that $\{[\mathfrak{s}_{1}^{(n_{2},n_{3})}],[\mathfrak{s}_{2}^{(n_{3},n_{1})}],[\mathfrak{s}_{3}^{(n_{1},n_{2})}]\}$ gives us a pants decomposition of $\Gamma$, we only need to prove the pullback homotopy classes on $S$ is a pants decomposition of $S$. Thus we only need to find disjoint representatives in these three homotopy classes on $S$.

\begin{figure}
    \centering
    \begin{tikzpicture}

\draw[dashed, black, thick] (0,0)arc (180:0:2 and 0.3);
\draw[black, thick] (0,0)arc (-180:0:2 and 0.3);
\draw[thick] (0,0) -- (0,4);
\draw[thick] (4,0) -- (4,4);
\draw[black, thick] (0,4)arc (180:0:2 and 0.3);
\draw[black, thick] (0,4)arc (-180:0:2 and 0.3);
\draw[black, thick, dashed] (0,2)arc (180:0:2 and 0.3);
\draw[black, thick] (0,2)arc (-180:0:2 and 0.3);

\draw[blue] (2,-.3) -- (4,.5);
\draw[blue, dashed] (4,.5) -- (0,1.5);
\draw[blue] (4,2.5) -- (0,1.5);
\draw[blue, dashed] (4,2.5) -- (0,3.5);
\draw[blue] (0,3.5) -- (2,3.7);

\draw[red, dashed] (2,.3) -- (0,.5);
\draw[red] (4,1.5) -- (0,.5);
\draw[red, dashed] (4,1.5) -- (0,2.5);
\draw[red] (0,2.5) -- (4,3.5);
\draw[red, dashed] (4,3.5) -- (2,4.3);

\coordinate [label=180:{$C_3$}] (a) at (0,2);
\coordinate [label=0:{$\gamma_1$}] (b) at (4,.5);
\coordinate [label=0:{$\gamma_2$}] (c) at (4,3.5);



\end{tikzpicture}
    \caption{A neighborhood of $C_3$.}
    \label{localpicture}
\end{figure}
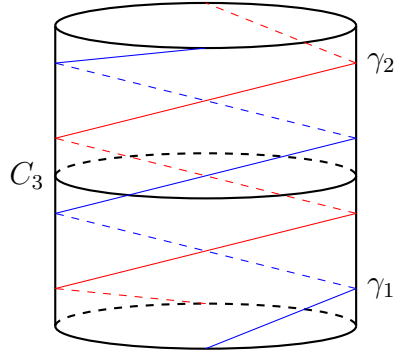

We notice that $\mathfrak{s}_{1}^{(n_{2},n_{3})}$ and $\mathfrak{2}_{i}^{(n_{3},n_{1})}$ both spin around $C_3$ by $n_3$ times, so we perturb these two curves in a tubular neighbourhood of $C_3$ as indicated in Figure \ref{localpicture}.

Similarly, we can use the same technique on $C_1$ and $C_2$, so we have three disjoint simple closed curves in the three free homotopy classes. Therefore the geodesic representatives of these three homotopy classes are disjoint simple closed geodesics on $S$, which helps us prove that $\{[\mathfrak{s}_{1}^{(n_{2},n_{3})}],[\mathfrak{s}_{2}^{(n_{3},n_{1})}],[\mathfrak{s}_{3}^{(n_{1},n_{2})}]\}$ is indeed us a pants decomposition of $\Gamma$.

Next we want to make the real lengths of $\{[\mathfrak{s}_{1}^{(n_{2},n_{3})}],[\mathfrak{s}_{2}^{(n_{3},n_{1})}],[\mathfrak{s}_{3}^{(n_{1},n_{2})}]\}$ be almost at the same size by adjusting the value of $(n_1,n_2,n_3)$. Consider the lattice $g:\Z^3\rightarrow\R^3$ by 
\begin{align*}
    g(n_1,n_2,n_3)=
    (\re(n_1\lle(C_1)-\sigma_{1}),
    \re(n_2\lle(C_2)-\sigma_{2}),
    \re(n_2\lle(C_2)-\sigma_{2})),
\end{align*}
then we know there is a constant $m_1>0$, such that for any point $P(x,y,z)\in\R^3$, we can find $n_1,n_2,n_3\in\Z$ satisfying 
\begin{align*}
    d(g(n_1,n_2,n_3),P)<m_1,
\end{align*}
here $d$ is the standard distance function in $\R^3$. Especially, for $R'$ big enough, we can find $n_1,n_2,n_3$ such that the distance between $g(n_1,n_2,n_3)$ and $(R',R',R')$ is less than $m_1$. Since $R'$ is big enough, we can assume $n_i$'s are positive and $n_i>N$, here $N$ comes from (\ref{n1n2n3}). Therefore there is $R_1>0$, such that for any $R'>R_1$, there exist $n_i\in\Z$ and $n_i>N$ satisfying that 
\begin{align}\label{m1}
    |\re(n_i\lle(C_i)-\sigma_{i})-R'|<m_1,
\end{align}
for $i\in\Z/3\Z$.

Let $m_2=\frac{1}{2}(m_1+\sqrt2\ep)$, then for any $R_0>0$, let $R'=\max\{R_1,2R_0+2\Sigma_{i=1}^{3}|\sigma_{i}|,3m_1\}+1$. Thus we can find $n_i$'s satisfying (\ref{m1}). Let $2R=\re(n_1\lle(c_1)+n_2\lle(c_2)+n_3\lle(c_3))-R'$, and we have
\begin{equation}\label{m1,2}
    |\re(n_{i+1}\lle(C_{i+1})+n_{i+2}\lle(C_{i+2})+\sigma_{i})-2R|=|-\re(n_{i}\lle(C_{i}))+\re(\sigma_{i})+R'|<m_1.
\end{equation}
Combine (\ref{n1n2n3}) and (\ref{m1,2}), we have
\[    
    |\re(\lle([\mathfrak{s}_{i}^{(n_{i+1},n_{i+2})}]))-2R|<\sqrt2\epsilon+m_1=2m_2,
\]
so
\[    
    |\re(\hl([\mathfrak{s}_{i}^{(n_{i+1},n_{i+2})}]))-R|<m_2.
\]
And we also have
\begin{align*}
    2R&=\re(n_1\lle(c_1)+n_2\lle(c_2)+n_3\lle(c_3))-R'\\
    &>(\re(\sigma_{1})+R'-m_1)+(\re(\sigma_{2})+R'-m_1)+(\re(\sigma_{3})+R'-m_1)-R'\\
    &=\sum_{i=1}^{3}\re(\sigma_{i})+2R'-3m_1\\
    &=(R'+\sum_{i=1}^{3}\re(\sigma_{i}))+(R'-3m_1)>2R_0.
\end{align*}

Let $m=m_2+\pi$. To summarize the above discussion, we actually prove the following theorem.

\begin{thm}\label{pd1} 
There is a constant $m>0$, such that for any positive real number $R_0$, there exist $R>R_0$ and $n_1,n_2,n_3\in \Z_+$ satisfying that $\{[\mathfrak{s}_{1}^{(n_{2},n_{3})}],[\mathfrak{s}_{2}^{(n_{3},n_{1})}],[\mathfrak{s}_{3}^{(n_{1},n_{2})}]\}$ forms a non-separating $(R,m)$-good pants decomposition of $\Gamma$.
\end{thm}

\subsection{Ideal triangulation and shears}

By Theorem \ref{pd1}, we can suppose $\{P_1,P_2\}$ is a $(R,m)$-good non-separating pants decomposition of $\Gamma$ right now, where $m$ is a constant and $R$ is big enough and will be determined later. We repeat the process in Section $5.1$ for this $(R,m)$-good pants decomposition instead of an arbitrary non-separating pants decomposition, and also inherit the notation. That means there is a constant $m'$, and for any $R_0^{'}>0$, there exist $R'>R_0^{'}$ and $n_1,n_2,n_3$ such that $\{[\mathfrak{s}_{1}^{(n_{2},n_{3})}],[\mathfrak{s}_{2}^{(n_{3},n_{1})}],[\mathfrak{s}_{3}^{(n_{1},n_{2})}]\}$ is a $(R',m')$-good non-separating pants decomposition of $\Gamma$. Now we want to use ideal triangulation of $\{P_1,P_2\}$ to estimate the short shears of the new $(R',m')$-good pants decomposition.

\begin{figure}
    \centering
    \begin{tikzpicture}[scale=0.2]

\draw[decoration={markings, mark=at position 0.25 with {\arrow{>}}}, postaction={decorate}] (0,0) ellipse (2 and 1) node[font=\fontsize{6}{6}\selectfont, below=.15] {$C_1$};
\draw[shift={(120:24)}, rotate=60, decoration={markings, mark=at position 0.75 with {\arrow{>}}}, postaction={decorate}] (0,0) ellipse (2 and 1)
node[font=\fontsize{6}{6}\selectfont, above left=.1] {$C_3$};
\draw[shift={(60:24)}, rotate=120, decoration={markings, mark=at position 0.25 with {\arrow{>}}}, postaction={decorate}] (0,0) ellipse (2 and 1)
node[font=\fontsize{6}{6}\selectfont, above right=.1] {$C_2$};


\draw[name path=A] (-2,0) .. controls ($(0,{sqrt(192)})+(210:3)$) .. ($(120:24)+(240:2)$)
\foreach \p in {0,1,...,20} {node[pos=\p*0.05, inner sep=0] (a\p) {}};
\draw[name path=B] (2,0) .. controls ($(0,{sqrt(192)})+(330:3)$) .. ($(60:24)+(-60:2)$)
\foreach \p in {0,1,...,20} {node[pos=\p*0.05, inner sep=0] (b\p) {}};
\draw[name path=C] ($(120:24)-(240:2)$) .. controls ($(0,{sqrt(192)})+(85:3)$) .. ($(60:24)-(-60:2)$)\foreach \p in {0,1,...,20} {node[pos=\p*0.05, inner sep=0] (c\p) {}};

\draw[green!85!blue] (c1) -- (a18);
\draw[green!85!blue] (c3) -- (a16);
\draw[green!85!blue] (c5) -- (a14);
\draw[green!85!blue] (c8) -- (a11);
\draw[green!85!blue] (b8) -- (a6);
\draw[green!85!blue] (b5) -- (a4);
\draw[green!85!blue] (b3) -- (a2);
\draw[green!85!blue, dashed] (c3) -- (a18);
\draw[green!85!blue, dashed] (c5) -- (a16);
\draw[green!85!blue, dashed] (c8) -- (a14);
\draw[green!85!blue, dashed] (b8) -- (a11);
\draw[green!85!blue, dashed] (b5) -- (a6);
\draw[green!85!blue, dashed] (b3) -- (a4);
\draw[green!85!blue, dashed] (b1) -- (a2);

\draw[red] (a1) -- (b2);
\draw[red] (a3) -- (b4);
\draw[red] (a5) -- (b6);
\draw[red] (a8) -- (b9);
\draw[red] (c12) -- (b14);
\draw[red] (c15) -- (b16);
\draw[red] (c17) -- (b18);
\draw[red, dashed] (a3) -- (b2);
\draw[red, dashed] (a5) -- (b4);
\draw[red, dashed] (a8) -- (b6);
\draw[red, dashed] (c12) -- (b9);
\draw[red, dashed] (c15) -- (b14);
\draw[red, dashed] (c17) -- (b16);
\draw[red, dashed] (c19) -- (b18);

\draw[blue] (b19) -- (c18);
\draw[blue] (b17) -- (c16);
\draw[blue] (b15) -- (c14);
\draw[blue] (b12) -- (c11);
\draw[blue] (a12) -- (c6);
\draw[blue] (a15) -- (c4);
\draw[blue] (a17) -- (c2);
\draw[blue, dashed] (b17) -- (c18);
\draw[blue, dashed] (b15) -- (c16);
\draw[blue, dashed] (b12) -- (c14);
\draw[blue, dashed] (a12) -- (c11);
\draw[blue, dashed] (a15) -- (c6);
\draw[blue, dashed] (a17) -- (c4);
\draw[blue, dashed] (a19) -- (c2);

\draw (30,11) circle (11);
\coordinate (u) at ($(30,11)+(-60:11)$);
\coordinate (v) at ($(30,11)+(-120:11)$);
\draw[decoration={markings, mark=at position 0.5 with {\arrow{>}}}, postaction={decorate}] (u) -- (v) node[above, pos=0.2,font=\fontsize{8}{8}\selectfont] {$\tilde{C_1}$};
\coordinate (w) at ($(30,11)+(60:11)$);
\coordinate (x) at ($(30,11)+(0:11)$);
\draw[decoration={markings, mark=at position 0.5 with {\arrow{>}}}, postaction={decorate}] (w) -- (x) node[left, pos=0.3,font=\fontsize{8}{8}\selectfont] {$\tilde{C_2}$};
\coordinate (y) at ($(30,11)+(180:11)$);
\coordinate (z) at ($(30,11)+(120:11)$);
\draw[decoration={markings, mark=at position 0.5 with {\arrow{>}}}, postaction={decorate}] (y) -- (z) node[right, pos=0.1,font=\fontsize{8}{8}\selectfont] {$\tilde{C_3}$};
\draw[red] (v) -- (x);
\draw[blue] (z) -- (x);
\draw[green!85!blue] (v) -- (z);

\end{tikzpicture}
    \caption{An ideal triangulation of $P_1$ and its corresponding picture in the universal cover $\HH^3$, where $\tilde{C_i}$ is a lift of $C_i$ for $i=1,2,3$ and all these lifts may not be on the same hyperbolic plane.}
    \label{idealtriangulation}
\end{figure}
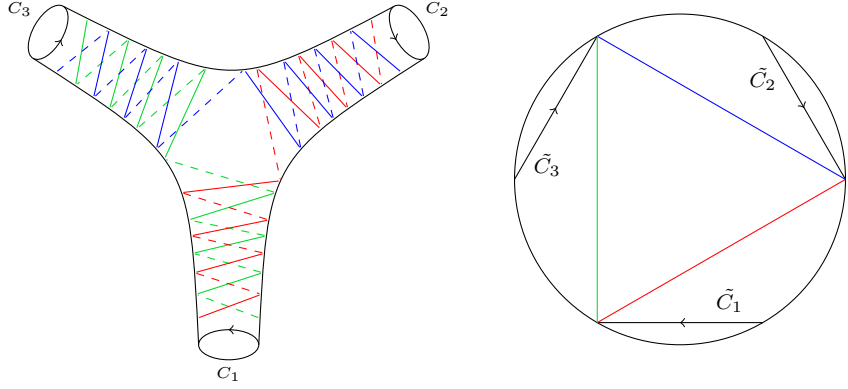

We fix the orientations of $C_i$'s, $i=1,2,3$, then there is a unique ideal triangulation of $P_1$ such that the sides of ideal triangles asymptotic to each $C_j$ in the same direction as the orientation of this $C_j$. This ideal triangulation can be also determined in the hyperbolic 3-space $\HH^3$, where one can connect the attracting points of the lifts of $C_i$'s as shown in Figure \ref{idealtriangulation}. Similarly, there is a unique ideal triangulation of $P_2$ such that the sides of the ideal triangles are asymptotic to each $C_j$ in the reverse direction of the orientation of this $C_j$. In each $P_k$, two ideal triangles share three sides, and we let $p_{k}^{i}$ be the common side approaching $C_{i+1}$ and $C_{i+2}$ in $P_k$ with the orientation spinning from $C_{i+1}$ to $C_{i+2}$. 

For each pair $\{p_{1}^{i},p_{2}^{i}\}$, consider all orthogeodesics between them, and we can label them in a unique way as $\mu_{i+1,k}^{i}$ and $\mu_{i+2,k}^{i}$ for all $k\in\Z$ such that $\delta^{i}_{n_{i+1},n_{i+2}}$ is freely homotopic to $[\mathfrak{s}_{i}^{(n_{i+1},n_{i+2})}]$, where $\delta^{i}_{n_{i+1},n_{i+2}}$ is defined as follows. For any $m,n$, $p_1^{i}$ and $p_2^{i}$ are orthogeodesics between $\mu_{i+1,m}^{i}$ and $\mu_{i+2,n}^{i}$. We denote the orthogonal segment between $\mu_{i+1,m}^{i}$ and $\mu_{i+2,n}^{i}$ on $p_1^{i}$ by $p^{i}_{1,m,n}$, and the one on $p_2^{i}$ by $p^{i}_{2,m,n}$, with the same orientation as $p_1^{i}$ and $p_2^{i}$. For a tuple of positive integers $(n_1,n_2,n_3)$, we define
\begin{equation*}
    \delta^{i}_{n_{i+1},n_{i+2}}=\mu^{i}_{i+1,n_{i+1}}p_{2,n_{i+1},n_{i+2}}^{i}(\mu^{i}_{i+2,n_{i+2}})^{-1}(p_{1,n_{i+1},n_{i+2}}^{i})^{-1}, i\in\Z/3\Z,
\end{equation*}
which are piecewise geodesic curves. Fix the orientation of $\mu_{j,k}^{i}$ from $p_{1}^{i}$ to $p_{2}^{i}$, then each piece of $\delta^{i}_{n_{i+1},n_{i+2}}$ is oriented and orthogonal to the two intersecting pieces. So there exists a unique framing at each joint point, such that there are framed segments $\mathfrak{w}_{i}^{(n_{i+1})},\mathfrak{x}_{i}^{(n_{i+1},n_{i+2})},\mathfrak{y}_{i}^{(n_{i+2})}$ and $\mathfrak{z}_{i}^{(n_{i+1},n_{i+2})}$, whose carriers are $\mu^{i}_{i+1,n_{i+1}},p_{2,n_{i+1},n_{i+2}}^{i},(\mu^{i}_{i+2,n_{i+2}})^{-1}$ and $(p_{1,n_{i+1},n_{i+2}}^{i})^{-1}$,
forming a zigzag continuous cycle. Let $\mathfrak{t}_{i}^{(n_{i+1},n_{i+2})}=\mathfrak{w}_{i}^{(n_{i+1})}\mathfrak{x}_{i}^{(n_{i+1},n_{i+2})}\mathfrak{y}_{i}^{(n_{i+2})}\mathfrak{z}_{i}^{(n_{i+1},n_{i+2})}$.

Next we want to estimate the length of each segment of $\mathfrak{t}_{i}^{(n_{i+1},n_{i+2})}$'s. Without loss of generality, we assume $i=1$. Let $\lle(C_1)=\lambda$ and $\lle(\mathfrak{w}_{1}^{(k)})=d_{k}$. By composition with M\"obius transformations, we can assume the geodesic $\rho$ connecting $-1$ and $0$ in $\HH^3$ is a lift of $p_1^1$. Let $z_k\in\C$ satisfying that the geodesic $\psi_k$ connecting $\infty$ and $z_k$ is a lift of $p_2^1$ whose distance to $\rho$ is $d_k$. Then we have
\begin{equation*}
    z_k=e^{-k\lambda}z_0,
\end{equation*}
for all $k\in\Z$. By elementary hyperbolic geometry, we know \begin{equation*}
    e^{d_k}=-(1+2z_k+2\sqrt{z_k^2+z_k}).
\end{equation*}
Then when $k\to+\infty$, we have $|z_k|\to 0$. Therefore
\begin{equation*}
    \lim\limits_{k\to+\infty}e^{d_k}=-1,
\end{equation*}
and 
\begin{equation*}
    |d_k-i\pi|\sim |e^{d_k-i\pi}-1|=\left|2z_k+2\sqrt{z_k^2+z_k}\right|\sim 2e^{-k\re(\lambda)/2}.
\end{equation*}
We can also have similar estimates for all other $\mathfrak{w}_{i}^{(k)}$'s and $\mathfrak{y}_{i}^{(k)}$'s. Hence for any $\ep>0$, there exists 
$N_\ep$ such that when $k>N_\ep$, we have
\begin{equation}\label{mu}
    |\lle(\mathfrak{w}_{i}^{(k)})-i\pi|<\ep, |\lle(\mathfrak{y}_{i}^{(k)})-i\pi|<\ep
\end{equation}
for any $i\in\Z/3\Z$. On the other hand, by the symmetry of $P_1$ and $P_2$, we have
\begin{equation*}
    l(\mathfrak{x}_{1}^{(n_{2},n_{3})})=l(\mathfrak{z}_{1}^{(n_{2},n_{3})}).
\end{equation*}
Therefore
\begin{align}
\begin{split}
    l(\mathfrak{x}_{1}^{(n_{2},n_{3})})=l(\mathfrak{z}_{1}^{(n_{2},n_{3})})
    &=\frac{1}{2}\left(l(\mathfrak{t}_{1}^{(n_{2},n_{3})})-l(\mathfrak{w}_{1}^{(n_2)})-l(\mathfrak{y}_{1}^{(n_3)})\right)\\
    &\geq\frac{1}{2}\left(l([\mathfrak{t}_{1}^{(n_{2},n_{3})}])-l(\mathfrak{w}_{1}^{(n_2)})-l(\mathfrak{y}_{1}^{(n_3)})\right)\\
    &=\frac{1}{2}\left(l([\mathfrak{s}_{1}^{(n_{2},n_{3})}])-l(\mathfrak{w}_{1}^{(n_2)})-l(\mathfrak{y}_{1}^{(n_3)})\right).
\end{split}
\end{align}
Thus for any $\epsilon>0$, we can choose $R_0^{'}$ big enough such that each $n_i$ is big, $|\lle(\mathfrak{w}_{i}^{(k)})-i\pi|<\ep$ and $   |\lle(\mathfrak{y}_{i}^{(k)})-i\pi|<\ep$, and then $\{[\mathfrak{t}_{1}^{(n_{2},n_{3})}],[\mathfrak{t}_{2}^{(n_{3},n_{1})}],[\mathfrak{t}_{3}^{(n_{1},n_{2})}]\}$ is a $(R',m')$-good non-separating pants decomposition for some $R'>R'_0$ and $\mathfrak{t}_{i}^{(n_{i+1},n_{i+2})}$'s are $(R'-m'-\epsilon)/2,\epsilon/2)$-zigzag continuous cycles, for $i\in\Z/3\Z$.

Now we want to use the ideal triangle shears along $p_{i}^{j}$ in $P_i$ to estimate the long shears of pants decomposition $\{[\mathfrak{t}_{1}^{(n_{2},n_{3})}],[\mathfrak{t}_{2}^{(n_{3},n_{1})}],[\mathfrak{t}_{3}^{(n_{1},n_{2})}]\}$ along $[\mathfrak{t}_j^{n_{j+1},n_{j+2}}]$. Let $T_i^1$ and $T_i^2$ be the two ideal triangles in $P_i$. $i\in\Z/2\Z$. In each ideal triangle $T_i^k$, let $q_i^{k,j}$ be the altitude to the side $p_i^j$ with orientation as $\left[0,\infty\right)$, for $j\in\Z/3\Z$ and $i\in\Z/2\Z$, $k\in\Z/2\Z$. For pants decomposition $\{[\mathfrak{t}_{1}^{(n_{2},n_{3})}],[\mathfrak{t}_{2}^{(n_{3},n_{1})}],[\mathfrak{t}_{3}^{(n_{1},n_{2})}]\}$, let $P'_1$ and $P'_2$ be two pairs of pants and $\xi_{k,n_j}^j$ be the third connection of $[\mathfrak{t}_j^{n_{j+1},n_{j+2}}]$ in $P'_{k}$ for $k\in\Z/2\Z$ and $j\in\Z/3\Z$. 

We also want to label those orthogeodesics between $q_1^{k,j}$ and $q_2^{k,j}$ as $\tau^{k,j}_{n_j}$ for any $k$, $j$ and $n_j$ such that $q_1^{k,j,n_j}\tau^{k,j}_{n_j}(q_2^{k,j,n_j})^{-1}$ is in the same free homotopy class as $\xi^j_{k,n_j}$, for any $n_j$, $k\in/Z/2\Z$ and $j\in\Z/3\Z$. Here $q_i^{k,j,n_j}$ is the segment on $q_i^{k,j}$ between $\tau^{k,j}_{n_j}$ and $p_i^j$. Similar with (\ref{mu}), we have
\begin{equation}\label{taulength}
    \lle(\tau^{k,j}_{n_j})-i\pi\sim e^{-n_j \hl(C_j)}\to 0, \mathrm{as}\ n_j\rightarrow\infty, k=1,2.
\end{equation}

Let $X_i^{k,j}$ be the intersection between $q_i^{k,j}$ and $p_i^j$, $Y_i^{k,j}$ be the intersection between $\xi^{j}_{k,n_j}$ and $[\mathfrak{t}_j^{n_{j+1},n_{j+2}}]$ in $P_i$, then we will show that $X_i^{k,j}$ and $Y_i^{k,j}$ are close to each other by the following lemma when $n_j$'s are big enough. For two points $A,B\in\HH^3$, we denote by $d(A,B)$ the distance between $A$ and $B$. And the estimate of $d(X_i^{k,j},Y_i^{k,j})$ is given by the following lemma.

\begin{lem}\label{lem5.8}
There exist $W_0,\omega_0>0$ and a function $\kappa(\omega)$ satisfying
\begin{equation*}
    \lim\limits_{\omega\rightarrow0}\kappa(\omega)=0,
\end{equation*}
such that when $0<\omega<\omega_0$, the following statement holds: Suppose $\delta_1$ and $\delta_2$ are two non-intersecting geodesics in $\HH^3$, $\xi$ is the orthogeodesic between them with $Y_i$ the intersection of $\xi$ and $\delta_i$, $i=1,2$ and $d(Y_1,Y_2)>W_0$. Suppose $A_1,A_2,B_1,B_2$ are four points satisfying $d(A_i,B_i)>W_0$, $d(A_i,\delta_i)<\omega$ and $d(B_i,\delta_i)<\omega$ for $i=1,2$. Let $X_i$ be a point on the geodesic segment $A_i B_i$ and $\zeta_i$ be the perpendicular geodesic to $A_i B_i$ which passes through $X_i$. If the following are satisfied:
\begin{enumerate}

\item $\zeta_1$ and $\zeta_2$ are disjoint;

\item $|d_{\C}(\zeta_1,\zeta_2)-i\pi|<\omega$,

\end{enumerate}
then $d(X_i,Y_i),d(Y_i,Z_i),d(Z_i,X_i)<\kappa(\omega)$.
\end{lem}

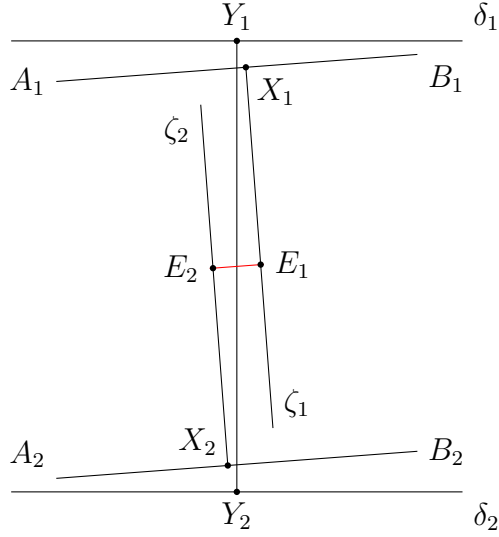
\begin{figure}
    \centering
    \begin{tikzpicture}[scale=0.6]

\draw (-5,5) -- (5,5);
\draw (-5,-5) -- (5,-5);
\draw (0,-5) -- (0,5);

\draw (-4,4.1) -- (4,4.7);
\draw (-4,-4.7) -- (4,-4.1);

\draw [name path=zeta1] (0.2,4.415) -- (0.8,-3.585);
\draw [name path=zeta2] (-0.2,-4.415) -- (-0.8,3.585);

\draw[red,name path=tau] (0.6,0.045) -- (-0.6,-0.045);

\coordinate [label=60:{$\delta_1$}] (a) at (5,5);
\coordinate [label=90:{$Y_1$}] (b) at (0,5);
\coordinate [label=-60:{$\delta_2$}] (a) at (5,-5);
\coordinate [label=-90:{$Y_2$}] (b) at (0,-5);

\coordinate [label=180:{$A_1$}] () at (-4,4.1);
\coordinate [label=170:{$A_2$}] () at (-4,-4.7);
\coordinate [label=-10:{$B_1$}] () at (4,4.7);
\coordinate [label=0:{$B_2$}] () at (4,-4.1);

\coordinate [label=-45:{$X_1$}] () at (0.2,4.415);
\coordinate [label=135:{$X_2$}] () at (-0.2,-4.415);
\coordinate [label=0:{$E_1$}] () at (0.6,0.045);
\coordinate [label=180:{$E_2$}] () at (-0.6,-0.045);

\coordinate [label=45:{$\zeta_1$}] () at (0.8,-3.585);
\coordinate [label=225:{$\zeta_2$}] () at (-0.8,3.585);

\fill [black] (0,-5) circle (2pt);
\fill [black] (0,5) circle (2pt);
\fill [black] (0.2,4.415) circle (2pt);
\fill [black] (-0.2,-4.415) circle (2pt);
\fill [name intersections={of=zeta1 and tau}] (intersection-1) circle (2pt);
\fill [name intersections={of=zeta2 and tau}] (intersection-1) circle (2pt);

\end{tikzpicture}
    \caption{A picture for Lemma \ref{lem5.8}}
    \label{lemma5.8pic}
\end{figure}

\begin{proof}
We will prove this lemma by contradiction and limit process. First we take $W_0$ big enough and $\omega_0$ small enough, for example $W_0>10^{10}$ and $\omega_0<10^{-10}$. And we can also assume $d(A_i,B_i)=W_0$ by taking a shorter geodesic segment with the position of $X_i$ unchanged. Let $\tau$ denote the common orthogonal geodesic between $\zeta_1$ and $\zeta_2$ and $E_i=\tau\cap\zeta_i$ for $i=1,2$ and $Z_1Z_2$ be the common orthogonal geodesic between $A_1B_1$ and $A_2B_2$ with $Z_i\in A_iB_i$. A set of geodesics and geodesic segments $\{\delta_1,\delta_2,\xi,A_1 B_1,A_2 B_2,Z_1Z_2,\zeta_1,\zeta_2,\tau\}$ is called a \emph{picture of $\omega$}, denoted by $\mathcal{P}(\omega)$, if all the requirements in the above statement are satisfied, i.e. $d(Y_1,Y_2)>W_0$, $d(A_i,B_i)=W_0$, $d(A_i,\delta_i)<\omega$, $d(B_i,\delta_i)<\omega$ for $i=1,2$, $\zeta_1$ and $\zeta_2$ are disjoint, $|d_{\C}(\zeta_1,\zeta_2)-i\pi|<\omega$.

We first study $d(X_i,Y_i)$. Suppose the inequality is not correct, then there exist a decreasing sequence $\{\omega^{(n)}\}$ which converges to 0, a sequence of pictures $\{\mathcal{P}^{(n)}(\omega^{(n)})\}$ and $\kappa_0>0$ such that in each $\mathcal{P}^{(n)}(\omega^{(n)})$, $d(X_1^{(n)},Y_1^{(n)})<\kappa_0$ and $d(X_2^{(n)},Y_2^{(n)})<\kappa_0$ cannot hold together. Without loss of generality, we can assume that $d(X_1^{(n)},Y_1^{(n)})\geq\kappa_0>0$ for all $n$ by passing to a subsequence. 

Now we take the compactified hyperbolic plane $\overline{\mathbb{H}^3}$ and study $\{\mathcal{P}^{(n)}(\omega^{(n)})\}$ on $\overline{\mathbb{H}^3}$. Since $\overline{\mathbb{H}}^3$ is compact, a sequence of points in $\overline{\mathbb{H}}^3$ will have a converging subsequence and its limit can be a point in $(\overline{\mathbb{H}^3})^{\circ}$ or a point on $\partial\overline{\mathbb{H}^3}$ at infinity. Besides, a sequence of geodesics also has a converging subsequence, and its limit can be a geodesic in $(\overline{\mathbb{H}^3})^{\circ}$ or a point on $\partial\overline{\mathbb{H}^3}$ at infinity as a degenerate geodesic.

Because each object in the picture $\mathcal{P}^{(n)}(\omega^{(n)})$ is either a point or a geodesic, we can assume $\{\mathcal{P}^{(n)}(\omega^{(n)})\}$ has a limit $\{\mathcal{P}^{(\infty)}\}$ by passing to subsequences. By applying isometries of $\mathbb{H}^3$, we can fix $E_1^{(i)}$ and $\zeta_1^{(i)}$, i.e. 
\begin{equation*}
    E_1^{(1)}=E_1^{(2)}=\cdots=E_1^{(n)}=\cdots=E_1^{(\infty)}
\end{equation*}
and
\begin{equation*}
    \zeta_1^{(1)}=\zeta_1^{(2)}=\cdots=\zeta_1^{(n)}=\cdots=\zeta_1^{(\infty)}.
\end{equation*}
Then by
\begin{equation*}
    d(E_1^{(\infty)},E_2^{(n)})=d(E_1^{(n)},E_2^{(n)})=l(\tau^{(n)})<\omega^{(n)},
\end{equation*}
we have
\begin{equation*}
    E_2^{(\infty)}=\lim\limits_{n\to+\infty}E_2^{(n)}=E_1^{(\infty)}.
\end{equation*}
Moreover, since $\tau=E_1 E_2$ is perpendicular to both $\zeta_1$ and $\zeta_2$, and $|\lle(\tau^{(n))}|<\omega^{(n)}$, so we know
\begin{equation*}
    \zeta_2^{(\infty)}=\lim\limits_{n\to+\infty}\zeta_2^{(n)}=\zeta_1^{(\infty)},
\end{equation*}
and we let $\zeta^{(\infty)}=\zeta_1^{(\infty)}=\zeta_2^{(\infty)}$. Then since $X_1^{(n)}=\zeta_1^{(n)}\cap\delta_1^{(n)}$ and $X_2^{(n)}=\zeta_2^{(n)}\cap\delta_2^{(n)}$, $X_1^{(\infty)}$ and $X_2^{(\infty)}$ are on $\zeta^{(\infty)}$.

By $d(A_1^{(n)},\delta_1^{(n)})<\omega^{(n)}$ and $\{\omega^{(n)}\}$ decreasing to 0, we know that $A_1^{(\infty)}\in\delta_1^{(\infty)}$. Similarly, we have $X_1^{(\infty)},B_1^{(\infty)}\in\delta_1^{(\infty)}$ and $A_2^{(\infty)},X_2^{(\infty)},B_2^{(\infty)}\in\delta_2^{(\infty)}$. Therefore $\zeta^{(\infty)}$ is orthogonal to $\delta_1^{(\infty)}$ and $\delta_2^{(\infty)}$, because $\zeta_i^{(n)}$ is perpendicular to $A_i^{(n)} B_i^{(n)}$. Hence $\zeta^{(\infty)}$ is the orthogeodesic between $\delta_1^{(\infty)}$ and $\delta_2^{(\infty)}$, i.e $\zeta^{(\infty)}$ coincides with $\xi^{(\infty)}$. Next we will show that the distance between $X_1^{(\infty)}$ and $Y_1^{(\infty)}$ is 0.

\textbf{Case I.} If $X_1^{(\infty)}$ is in $(\overline{\mathbb{H}^3})^{\circ}$, then we know
\begin{equation*}
    Y_1^{(\infty)}=\xi^{(\infty)}\cap\delta_1^{(\infty)}=\zeta^{(\infty)}\cap A_1^{(\infty)} B_1^{(\infty)}=X_1^{(\infty)}.
\end{equation*}

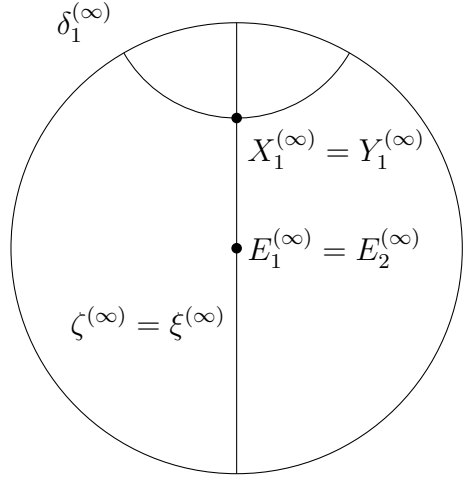
\begin{figure}
    \centering
    \begin{tikzpicture}[scale=1]

\draw (0,0) circle [radius=3];

\coordinate [label=225:{}] (a1) at (120:3);
\coordinate [label=225:{}] (a2) at (60:3);

\draw [name path=zeta] (120:3) arc (-150:-30:{sqrt(3)});
\draw [name path=xi] (0,3) -- (0,-3);
\fill [name intersections={of=zeta and xi}] (intersection-1) circle (2pt) node[black,below right] {$X_1^{(\infty)}=Y_1^{(\infty)}$};
\fill [black] (0,0) circle (2pt) node[right] {$E_1^{(\infty)}=E_2^{(\infty)}$};

\coordinate [label=180:{$\zeta^{(\infty)}=\xi^{(\infty)}$}] () at (0,-1);
\coordinate [label=135:{$\delta_1^{(\infty)}$}] () at (120:3);

\end{tikzpicture}
    \caption{$\mathcal{P}^{(\infty)}$ in \textbf{Case I}.}
    \label{case1}
\end{figure}

\textbf{Case II.} If $X_1^{(\infty)}$ is on $\partial\overline{\mathbb{H}^3}$, i.e. the boundary at infinity of $\HH^3$. Then since $\delta_1^{(\infty)}$ is perpendicular to $\zeta^{(\infty)}$ and passes through $X_1^{(\infty)}$, we know $\delta_1^{(\infty)}$ is a point on $\partial\overline{\mathbb{H}^3}$ which coincides with $X_1^{(\infty)}$. Moreover, $X_1^{(\infty)}=Y_1^{(\infty)}$ on $\partial\overline{\mathbb{H}^3}$. Now we let each $\{\mathcal{P}^{(n)}(\omega^{(n)})\}$ be acted by an isometry, which translates along $\zeta_1^{(n)}$, and denoted by $\{\mathcal{P'}^{(n)}(\omega^{(n)})\}$, such that 
\begin{equation*}
    {X'}_1^{(1)}={X'}_1^{(2)}=\cdots={X'}_1^{(n)}=\cdots={X'}_1^{(\infty)}.
\end{equation*}
Then $\{\mathcal{P'}^{(n)}(\omega^{(n)})\}$ has a limit picture $\mathcal{P'}^{(\infty)}$ in $\overline{\mathbb{H}^3}$ by passing to subsequences.

\begin{figure}
    \centering
    \begin{tikzpicture}[scale=1]

\draw (0,0) circle [radius=3];

\draw [name path=xi] (0,3) -- (0,-3);
\fill [black] (0,0) circle (2pt) node[right] {$E_1^{(\infty)}=E_2^{(\infty)}$};

\fill [black] (0,3) circle (2pt) node[above] {$\delta_1^{(\infty)},X_1^{(\infty)},Y_1^{(\infty)}$};

\coordinate [label=180:{$\zeta^{(\infty)}=\xi^{(\infty)}$}] () at (0,-1);

\draw (8,0) circle [radius=3];

\draw (5,0) -- (11,0);
\draw [name path=xi'] (8,3) -- (8,-3);
\fill [black] (8,0) circle (2pt) node[above right] {${X'}_1^{(\infty)}={Y'}_1^{(\infty)}$};
\coordinate [label=90:{${\delta'}_1^{(\infty)}$}] () at (6,0);

\fill [black] (8,-3) circle (2pt) node[below] {${E'}_1^{(\infty)},{E'}_2^{(\infty)}$};

\coordinate [label=180:{${\zeta'}^{(\infty)}={\xi'}^{(\infty)}$}] () at (8,-1);
\end{tikzpicture}
    \caption{$\mathcal{P}^{(\infty)}$(the left one) and $\mathcal{P'}^{(\infty)}$ in \textbf{Case II}.}
    \label{case2}
\end{figure}
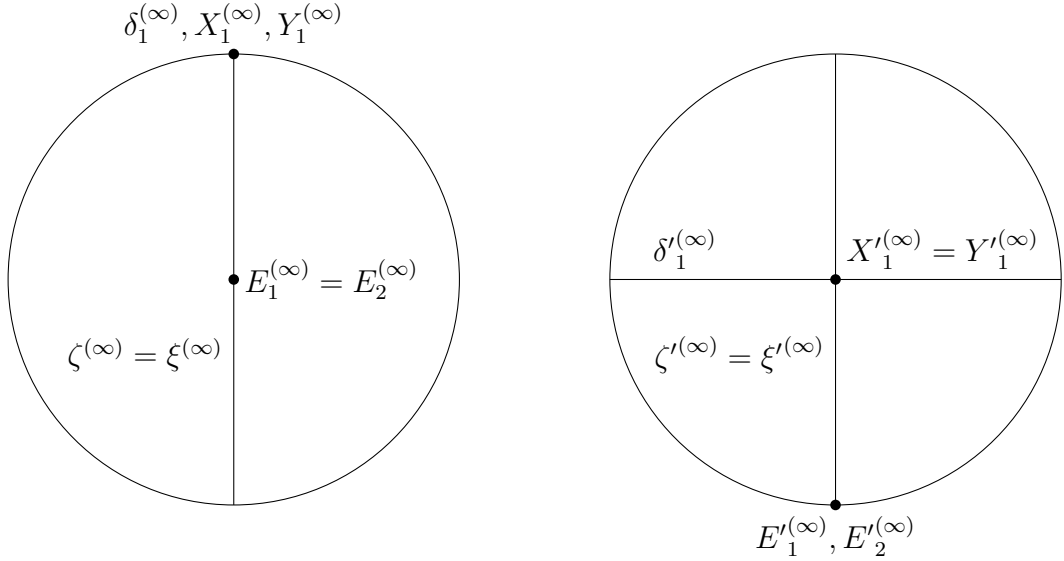

Since $X_1^{(\infty)}\in\mathcal{P}^{(n)}(\omega^{(n)})$ is on $\partial\overline{\mathbb{H}^3}$ and ${E}_1^{(\infty)}\in\mathcal{P}^{(n)}(\omega^{(n)})$ is in $(\overline{\mathbb{H}^3})^{\circ}$, we know ${E'}_1^{(\infty)}\in\mathcal{P'}^{(n)}(\omega^{(n)})$ is on $\partial\overline{\mathbb{H}^3}$ after translations. Therefore ${E'}_1^{(\infty)}={E'}_2^{(\infty)}$ on $\partial\overline{\mathbb{H}^3}$. Since isometries will keep the distance between any pair of points, so if two sequences of points have the same limits in $\mathcal{P}^{(\infty)}$, then their limits are also the same point in $\mathcal{P'}^{(\infty)}$. Therefore we still have ${\zeta'}^{(\infty)}={\zeta'}_1^{(\infty)}={\zeta'}_2^{(\infty)}$. By ${\delta'}_2^{(\infty)}$ perpendicular to ${\zeta'}_2^{(\infty)}={\zeta'}^{(\infty)}$, we know ${\zeta'}^{(\infty)}$ is orthogonal to both ${\delta'}_1^{(\infty)}$ and ${\delta'}_2^{(\infty)}$. Hence ${\zeta'}^{(\infty)}={\xi'}^{(\infty)}$ and we have
\begin{equation*}
    {Y'}_1^{(\infty)}={\xi'}^{(\infty)}\cap{\delta'}_1^{(\infty)}={\zeta'}^{(\infty)}\cap {\delta'}_1^{(\infty)}=X_1^{(\infty)}.
\end{equation*}
Thus $d({X}_1^{(\infty)},{Y}_1^{(\infty)})=0$.

Now we put these two cases together, and we always have $d({X}_1^{(\infty)},{Y}_1^{(\infty)})=0$, which contradicts to the assumption that $d(X_1^{(n)},Y_1^{(n)})\geq\kappa_0>0$ for all $n$. Contradiction!

And we can have similar results for $d(Y_i,Z_i)$ and $d(Y_i,X_i)$, so the lemma is proved.

\end{proof}

Now we require $R>100m$, then for $\de>0$ with $\kappa(\de)<m/2$, let $\ep_0$ and $L_0$ in Lemma \ref{zglem} for $\de$. By (\ref{mu}) and (\ref{taulength}) we can find $n_j$'s big enough such that
\begin{equation*}
    |\lle(\mu_{j,n_j}^{i})-i\pi|<\ep_0,
\end{equation*}
\begin{equation*}
    \frac{R'-m'-\ep_0}{2}>L_0
\end{equation*}
and
\begin{equation*}
    |\lle(\tau_{n_j}^{k,j})-i\pi|<\de.
\end{equation*}
Thus by Lemma \ref{zglem} and Lemma \ref{lem5.8}, we have
\begin{equation}\label{shearerror}
    d(X_i^{k,j},Y_i^{k,j})<\kappa(\de)<m/2,
\end{equation}
for $i\in\Z/2\Z$, $k\in\Z/2\Z$ and $j\in\Z/3\Z$.

On $p_i^j$, $X_i^{1,j}X_i^{2,j}$ is the ideal triangle shear between $T_i^1$ and $T_i^2$ along $p_i^j$. Therefore we know
\begin{equation*}
    \lle(X_i^{1,j}X_i^{2,j})=\hl(C_{j+1})+\hl(C_{j+2})-\hl(C_j),
\end{equation*}
for $j\in\Z/3\Z$ and $i\in\Z/2\Z$. By $\{P_1,P_2\}$ being a $(R,m)$-good pants decomposition, we have $R-m<c_t<R+m$, for $t=1,2,3$. Hence we have
\begin{equation}\label{trishear}
    R-3m<l(X_i^{1,j}X_i^{2,j})<R+3m.
\end{equation}
Thus by (\ref{shearerror}) and (\ref{trishear}), we have
\begin{equation}\label{longshear}
    R-4m<l(Y_i^{1,j}Y_i^{2,j})<R+4m.
\end{equation}
By definition, $Y_i^{1,j}Y_i^{2,j}$, $i=1,2$ are the long shears between $P'_1$ and $P'_2$ along $[\mathfrak{t}_j^{n_{j+1},n_{j+2}}]$. Since the short shear $s_j$ between $P'_1$ and $P'_2$ along $[\mathfrak{t}_j^{n_{j+1},n_{j+2}}]$ is the average of $\lle(Y_1^{1,j}Y_1^{2,j})$ and $\lle(Y_2^{1,j}Y_2^{2,j})$, we have
\begin{equation}\label{shortshear}
    R-4m<\re(s_j)<R+4m,
\end{equation}
which gives us upper and lower bounds of real parts of short shears of the pants decomposition $\{P'_1,P'_2\}$. And we can state the above result as the following theorem.

\begin{thm}\label{gpdwbs}
Suppose $\Gamma$ is a genus-2 quasi-Fuchsian group, then there exist $B^+>B^->0$ and $\delta>0$ such that for any $R_0>0$, there exists $R>R_0$ such that $\Gamma$ admits a nonseparating $(R,\delta)$-good pants decomposition with the real parts of twists in the interval $(B^-,B^+)$.
\end{thm}

At the end of this section, we want to use Lemma \ref{lem5.8} to refine Lemma \ref{lem5.8} itself for some future purpose. 

\begin{lem}\label{lem5.14}
For any $B,m>0$, there exist $R_0>0$ and $B'>0$ such that the following statement holds: Let $\delta_1$ and $\delta_2$ be two non-intersecting geodesics in $\HH^3$, $\xi$ be the orthogeodesic between them with $Y_i$ the intersection of $\xi$ and $\delta_i$, $i=1,2$ , and $A_1,A_2,B_1,B_2$ be four points satisfying $d(A_i,\delta_i)<B e^{-R/2}$ and $d(B_i,\delta_i)<B e^{-R/2}$ for $i=1,2$. Suppose $d_{\C}(\de_1,\de_2)$, $d(A_1,B_1)$ and $d(A_2,B_2)$ are within distance $m$ from $R$. Let $X_i$ be a point on the geodesic segment $A_i B_i$, $\zeta_i$ be the perpendicular geodesic to $A_i B_i$ which passes through $X_i$ and $\tau$ be the common orthogonal between $\zeta_1$ and $\zeta_2$, and the orientation of $\zeta_i$ is from $A_i B_i$ to $\tau$. If the following are satisfied:
\begin{enumerate}

\item $\zeta_1$ and $\zeta_2$ are disjoint;

\item $|d_{\C}(\zeta_1.\zeta_2)-i\pi|<B e^{-R/2}$,

\end{enumerate}
then $d(X_i,Y_i)<B' e^{-R/4}$.
\end{lem}

\begin{proof}

Let $Z_1 Z_2$ be the common orthogonal between $A_1 B_1$ and $A_2 B_2$ with $Z_i$ on geodesic $A_i B_i$. Then we will estimate $d(X_i, Z_i)$ and $d(Z_i, Y_i)$ separately. By Lemma \ref{lem5.8}, we can take $R_1>0$ such that when $R>R_1$, we have $d(X_i,Y_i),d(Y_i,Z_i),d(Z_i,X_i)<1$.

\begin{enumerate}[(a)]

\item Let $a=d_{\C}(\zeta_1,Z_1 Z_2)$, $b=d_{\C}(\zeta_2,Z_1 Z_2)$, $c=d_{\C}(\zeta_1,\zeta_2)$, $x=d_{\C}(A_1 B_1,\tau)$, $y=d_{\C}(A_2 B_2,\tau)$ and $z=d_{\C}(A_1 B_1, A_2 B_2)$, then we know 
\begin{equation}\label{5.15}
    |c-i\pi|<B e^{-R/2}.
\end{equation}

By Hyperbolic Cosine Rule for right-angled hexagon, we have
\begin{equation*}
    \cosh(c)=\frac{\cosh(x)\cosh(y)+\cosh(z)}{\sinh(x)\sinh(y)}.
\end{equation*}
Therefore
\begin{equation}\label{5.14cosh}
\begin{aligned}
    \cosh(z-i\pi)&=\cosh(x)\cosh(y)+\sinh(x)\sinh(y)\cosh(c-i\pi)\\
    &=\cosh(x+y)+\sinh(x)\sinh(y)(\cosh(c-i\pi)-1).
\end{aligned}
\end{equation}
And by Hyperbolic Sine Rule, we have 
\begin{equation}\label{5.14sinh}
    \frac{\sinh(a)}{\sinh(y)}=\frac{\sinh(b)}{\sinh(x)}=\frac{\sinh(c)}{\sinh(z)}=\frac{\sinh(c-i\pi)}{\sinh(z-i\pi)}.
\end{equation}

Since $d(Y_i,Z_i)<1$ and $d(Y_1,Y_2)>R-m$, we have $d(Z_1,Z_2)>R-m-2$. So there exists $R_2>0$ such that when $R>R_2$, we have
\begin{equation}\label{5.18}
    |\cosh(z-i\pi)|<2|\sinh(z-i\pi)|,
\end{equation}
\begin{equation}\label{5.19}
    \re(x+y)\geq\re(z-c)\geq R-m-2-Be^{-R/2}\geq R-m-3,
\end{equation}
\begin{equation}\label{5.20}
    |\cosh(x+y)|>\frac{|e^{x+y}|}{4}
\end{equation}
and
\begin{equation}\label{5.21}
    |\sinh(c-i\pi)|<2|c-i\pi|<2B e^{-R/2}<1/8.
\end{equation}

Because $|\sinh(t)|\leq |e^t|$ for any $\re(t)\geq0$, then by (\ref{5.15}), (\ref{5.14cosh}), (\ref{5.14sinh}), (\ref{5.18}), (\ref{5.19}), (\ref{5.20}) and (\ref{5.21}), we get
\begin{equation*}
    \begin{aligned}
        |\sinh(a)|&=|\sinh(c-i\pi)|\cdot\frac{|\sinh(y)|}{|\sinh(z-i\pi)|}\\
        &<2Be^{-R/2}\frac{|\sinh(y)|}{|\sinh(z-i\pi)|/2}\\
        &=4Be^{-R/2}\frac{|\sinh(y)|}{|\cosh(x+y)+\sinh(x)\sinh(y)(\cosh(c-i\pi)-1)|}\\
        &\leq4Be^{-R/2}\frac{|\sinh(y)|}{|\cosh(x+y)|-|\sinh(x)\sinh(y)(\cosh(c-i\pi)-1)|}\\
        &\leq 4Be^{-R/2}\frac{|e^y|}{|e^{x+y}|/4-|e^x e^y|\cdot|\cosh(c-i\pi)-1|}\\
        &=4Be^{-R/2}\frac{1}{|e^x|/4-|\cosh(c-i\pi)-1|}\\
        &\leq 4Be^{-R/2}\frac{1}{\frac{1}{4}-\frac{1}{8}}=32Be^{-R/2}.
    \end{aligned}
\end{equation*}

Hence $d(X_i,Z_i)=\re(a)<B_1 e^{-R/2}$ for some $B_1>0$.

\item Let $r=d(Y_1, Y_2)$ and $s=d(Z_1, Z_2)$. Since $d(X_i,Y_i),d(Z_i,X_i)<1$, therefore there exists $B_2>0$ such that 
\begin{equation}\label{ZY}
    d(Z_i,\de_i),d(Y_i,A_i B_i)<B_2 e^{-R/2}.
\end{equation}
Then $|r-s|<2B_2 e^{-R/2}$. Let $M_i$ be the projection from $Z_i$ on $\de_i$ and $N_1$ be the projection from $M_2$ to $\de_1$. Then
\begin{equation}
    d(M_2,N_1)\leq d(M_2, M_1)\leq d(M_2, Z_2) + d(Z_1, Z_2) + d(Z_1, M_1)<r+4B_2 e^{-R/2}.
\end{equation}

In the hyperbolic quadrilateral $Y_1Y_2M_2N_1$, let $u=d_{\C}(Y_1, Y_2)$, $v=d_{\C}(N_1, M_2)$ and $w=d_{\C}(Y_2, M_2)$. Then
\begin{equation}\label{uv}
    \re(u)=r\leq\re(v)<r+4B_2 e^{-R/2},
\end{equation}
and we have
\begin{equation*}
    \cosh(u+i\frac{\pi}{2})=\frac{\cosh(i\frac{\pi}{2})\cosh(w-i\frac{\pi}{2})+\cosh(v+i\frac{\pi}{2})}{\sinh(i\frac{\pi}{2})\sinh(w-i\frac{\pi}{2})},
\end{equation*}
which is simplified as 
\begin{equation*}
    \cosh(w)=\frac{\sinh(v)}{\sinh(u)}.
\end{equation*}
So by (\ref{uv}) and $r>R-m$, we know there exist $R_3>0$ and $B_3>0$ such that when $R>R_2$, we have
\begin{equation}\label{w}
    |w|<B_3 e^{-R/4}.
\end{equation}
Thus by (\ref{ZY}) and (\ref{w}), we have 
\begin{equation}
    d(Z_2, Y_2)<B_2 e^{-R/2} + B_3 e^{-R/4},
\end{equation}
and we can have a similar result for $d(Z_1, Y_1)$.
\end{enumerate}

Now take $R_0>\max\{R_1,R_2,R_3\}$ and combine the results from part (a) and part (b), we can find $B'>0$ such that when $R>R_0$, we have
\begin{equation*}
    d(X_i, Y_i)<B' e^{-R/4}.
\end{equation*}

\end{proof}

\section{Counting and matching pants in the 3-manifold}\label{S6}

In this section, we will first count good curves and good pants in compact hyperbolic 3-manifolds, and then match good pants along each good curve by Hall marriage Theorem to get a good assembly.

By Theorem \ref{gpdwbs}, we can take a nonseparating $(\overline{R},\delta)$-good pants decomposition of $\Gamma$ with real parts of short shears bounded by $(B^-,B^+)$, where $\overline{R}>0$ is big enough and will be determined later and $\de$ is a constant. Let $C_i$'s be the cuffs of this pants decomposition with $\hl(C_i)=R_i\in\C/2\pi i\Z$ and $s_i$ the short shear along $C_i$, $i\in\Z/3\Z$, then we have
\begin{equation}
    |R_i-\overline{R}|<\delta
\end{equation}
and
\begin{equation}
    \re(s_i)\in(B^-,B^+).
\end{equation}

For a pair of pants $P$ with cuffs $C'_i$, we call $P$ is $(R_i,\ep)_{i=1}^3$-good, if
\begin{equation*}
    |\hl(C'_i)-R_i|<\ep,
\end{equation*}
for $i=1,2,3$. One can check Section 2.2 in \cite{KW21} for $assemblies$ and in this paper an assembly will be always constructed only by pants. And we also want to change the definition of good assemblies a little bit to adapt our situation. We say an assembly $\mathcal{A}$ is $(R_i,s_i,\ep)_{i=1}^3$-good, if
\begin{enumerate}

\item Each pair of pants in $\mathcal{A}$ is $(R_i,\ep)_{i=1}^3$-good;

\item If two pants of $\mathcal{A}$ are glued along a curve $\gamma$ which is $(R_i,\ep)$-good for some $i$, then let $\al_1,\al_2\in N^1(\sqrt{\gamma})$ be two feet of these two pants, and we have
\begin{equation*}
    |\al_1-\al_2-(s_i+i\pi)|<\ep/\overline{R}.
\end{equation*}
Here these two pants should be oriented and induce opposite orientation on $\gamma$.
\end{enumerate}

\subsection{Counting good curves, geodesic connections and pants}
Now we will follow results in Section 3 and Section 5 of \cite{KW21} to help us count $(R_i,\ep)_{i=1}^3$-good pants in $\HH^3/G$, and we do not need to consider cusps since $G$ is cocompact.

We first want to count $(R,\ep)$-good curves in $\HH^3/G$ for any $R\in\C/2\pi i\Z$ and $\ep>0$, which follows from the Margulis argument. Let $\Gamma_{\ep,R}$ be the set of $(R,\ep)$-good curves in $\HH^3/G$, then by (3.1.1) in \cite{KW21}, we have
\begin{equation}
    \lim\limits_{R\rightarrow\infty}\frac{\#(\Gamma_{\ep,R})}{\ep^2 e^{4\re(R)}/\re(R)}=c_{\ep},
\end{equation}
where $c_{\ep}$ is a non-zero constant depending on $\ep$.

Our next step is to count geodesic connections between two geodesics. Let $\gamma_0,\gamma_1$ be two oriented closed geodesics. Then for a connection $\al$ between $\gamma_0$ and $\gamma_1$, we let $n_i(\al)$ be the unit vector that points in toward $\al$ at the point where $\al$ meets $\gamma_i$ and $\theta(\al)$ be the angle between the tangent
vector to $\gamma_1$ where it meets $\al$ and the parallel transport along $\al$ of the tangent vector to $\gamma_0$ where it meets $\al$, and define $w(\al) = l(\al)+i\theta(\al)$. Let
\begin{equation*}
    \mathbb{I}(\gamma_0,\gamma_1)=N^1(\gamma_0)\times N^1(\gamma_1)\times \C/2\pi i\Z,
\end{equation*}
then for each $\al$, we have a triple
\begin{equation*}
    \mathbf{I}(\al)=(n_0(\al),n_1(\al),w(\al))\in\mathbb{I}(\gamma_0,\gamma_1).
\end{equation*}
We fix the measure on $\mathbb{I}(\gamma_0,\gamma_1)$ by regarding $\C/2\pi i\Z$ as $S^1\times\R$ and then taking the product of Lebesgue measures on the first three coordinates times $e^{2t}dt$ on $\R$. And we also have a metric on $\mathbb{I}(\gamma_0,\gamma_1)$ that is the $L^2$ norm of the distances in each coordinate. Let $\mathcal{N}_{\eta}(A)$ be the set of points with distance less than $\eta$ to the set $A$ and $\mathcal{N}_{-\eta}(A)$ be the set of points with distance greater than $\eta$ to the complement of $A$. Then the following theorem (Theorem 3.2 in \cite{KW21}) holds:

\begin{thm}\label{countconn}
There exists $q>0$ depending on $G$ such that the following holds when $R^-$ is sufficiently large. Suppose $A\subset \mathbb{I}(\gamma_0,\gamma_1)$, and let $R^-$ be the infimum of the fourth coordinate of values in $A$. Let $\eta=e^{-q R^-}$, then the numbers of connections $\mathbf{n}(A)$ for $\al$ between $\gamma_0$ and $\gamma_1$ that have $\mathbf{I}(\al)\in A$ satisfies
\begin{equation*}
    (1-\eta)|\mathcal{N}_{-\eta}(A)|
    \leq 32\pi^2\mathbf{n}(A)|\HH^3/G|
    \leq (1+\eta)|\mathcal{N}_{\eta}(A)|.
\end{equation*}
\end{thm}

By letting $\gamma_0=\gamma_1=\gamma$ be a $(R_1,\ep)$-good curve in Theorem \ref{countconn}, we can have an estimate on the number of $(R_i,\ep)_{i=1}^3$-good pants which have $\gamma$ as a boundary. Actually we can almost follow the proof in Section 3.3 of \cite{KW21} with only changing the lengths of cuffs.

Let $\Pi_{\ep,R_i}$ and $\Pi^*_{\ep,R_i}$ be the set of all unoriented and oriented $(R_i,\ep)_{i=1}^3$-good pants in $\HH^3/G$ respectively, and $\Pi^*_{\ep,R_i}(\gamma)$ be the set of pants in $\Pi^*_{\ep,R_i}$ for which $\gamma$ is a cuff. We know that the unit normal bundle $N^1(\gamma)$ is a torsor for $\C/(\lle(\gamma)\Z+2\pi i\Z)$, then the map $n\mapsto n+\hl(\gamma)$ is an involution on $N^1(\gamma)$. We use $N^1(\sqrt{\gamma})$ to denote the quotient of $N^1(\gamma)$ by this involution, then $N^1(\sqrt{\gamma})$ is a torsor for $\C/(\hl(\gamma)\Z+2\pi i\Z)$. For $P\in\Pi^*_{\ep,R_i}(\gamma)$, let $n_0$ and $n_1$ be the first two coordinates of $\mathbf{I}(\al)$ where $\al$ is the third connection for $P$, then it turns out that we have a well-defined map $u:\Pi^*_{\ep,R_i}(\gamma)\to N^1(\sqrt{\gamma})$ by
\begin{equation*}
    u(P)=(n_0+n_1)/2.
\end{equation*}
So we can now state the following theorem and the proof is the same as the proof of Theorem 3.3 in \cite{KW21} with slight modification.

\begin{thm}\label{countpants}
There exists positive constant $q$ depending on $G$ such that for any $\ep>0$ the following holds when $R_i$'s are sufficiently large. Let $\gamma$ be an $(R_1,\ep)$-good curve. If $B\subset N^1(\sqrt{\gamma})$, then
\begin{equation*}
    (1-\xi)\mathrm{Vol}(\mathcal{N}_{-\xi}(B))
    \leq \frac{\#\{P\in\Pi^*_{\ep,R_i}(\gamma)\ |\ u(P)\in B\}}
    {C_{count}(\ep)\ep^4 e^{2\re(R_2)+2\re(R_3)-l(\gamma)}/|\HH^3/G|}
    \leq (1+\xi)\mathrm{Vol}(\mathcal{N}_{\xi}(B)),
\end{equation*}
where $\xi=e^{-q\overline{R}/2}$ and $C_{count}(\ep)\to 1$ as $\ep\to 0$.   
\end{thm}

\subsection{Matching pants}

At each $(R_i,\ep)$-good curve $\gamma$, we want to match each oriented $(R_i,\ep)_{i=1}^3$-good pants with $\gamma$ as a cuff with another such good pants which has the opposite orientation on $\gamma$. For $i=1,2,3$ and $\gamma\in \Gamma^*_{\ep,R_i}$, we define $\tau:N^1(\sqrt{\gamma})\to N^1(\sqrt{\gamma})$ by $\tau(v)=v+i\pi+s_i$. Now we will follow the idea in Section 5.2 of \cite{KW21} to prove the following theorem by using the Hall Marriage theorem:

\begin{thm}\label{matchpants}
For all $\ep$ there exists $R_0>0$ such that for all $\overline{R}>R_0$: Let $\gamma$ be an oriented $(R_i,\ep)$-good curve for some $i\in\{1,2,3\}$, where $R_i$'s are determined at the beginning of Section \ref{S6}, then there exists a permutation $\sigma_{\gamma}:\Pi_{\ep,R_i}(\gamma)\to\Pi_{\ep,R_i}(\gamma)$ such that
\begin{equation*}
    |\mathbf{foot}_\gamma(\sigma_{\gamma}(\pi))-\tau(\mathbf{foot}_\gamma(\pi))|<\ep/\overline{R},
\end{equation*}
for all $\pi\in\Pi_{\ep,R_i}(\gamma)$.
\end{thm}

Before we prove this theorem, we introduce the following notation. For $A\subset N^1(\sqrt{\gamma})$, let
\begin{equation*}
    \#A:=|\{\pi\in\Pi_{\ep,R_i}(\gamma):\mathbf{foot}_\gamma(\pi)\in A|.
\end{equation*}
And we recall the following proposition from Corollary 5.6 in \cite{KW21}:
\begin{prop}\label{prop}
If $A\subset N^1(\sqrt{\gamma})$ and $|\mathcal{N}_\eta(A)|\leq1/2|N^1(\sqrt{\gamma})|$, then 
\begin{equation*}
    \frac{|\mathcal{N}_\eta(A)|}{|A|}>1+\frac{\eta}{\re(R_i)}.
\end{equation*}
Moreover, 
\begin{equation*}
    \frac{|\mathcal{N}_\eta(A)|}{|A|}>1+\frac{\eta}{2\overline{R}},
\end{equation*}
when $\overline{R}$ is sufficiently large.
\end{prop}

\begin{proof}[Proof of Theorem \ref{matchpants}]
Let $\eta=\ep/\overline{R}$, $C=C_{count}(\ep)\ep^4 e^{2\re(R_2)+2\re(R_3)-l(\gamma)}/|\HH^3/G|$ and $\xi$ is as it appears in Theorem \ref{countpants}. By the Hall marriage theorem, Theorem \ref{matchpants} follows from the statement that $\#\mathcal{N}_\eta(A)\geq\#\tau(A)$ for every finite set $A$. Actually when $\overline{R}$ is large, we have 
\begin{equation*}
    \eta/2=\ep/2\overline{R}>5\overline{R}\xi>(4\overline{R}+1+\ep/2\overline{R})\xi=(4\overline{R}+1+\eta/2)\xi,
\end{equation*}
therefore by Proposition \ref{prop},
\begin{equation}\label{cp1}
    \frac{|\mathcal{N}_{\eta/2}(A)|}{|\mathcal{N}_\xi(A)|}>1+\frac{\eta/2-\xi}{2\overline{R}}>\frac{1+\xi}{1-\xi}.
\end{equation}
And the remaining of the proof is the same as the proof of Theorem 5.7 in \cite{KW21}.
\end{proof}

\section{Existence of quasiconformal map}\label{S7}

The goal of this section is to show that the gluing of good pants by Theorem \ref{matchpants} forms an assembly which is close to its perfect model. The idea follows from the appendix of \cite{KW21}, where the perfect model is considered as an assembly consisting of all identical pants.

\subsection{An estimate for matrix multiplication} We first quote an important theorem in \cite{KW21}, which helps us estimate the error of matrix multiplication. 

For any element $U$ of a Lie algebra (for a given Lie group), we let $U(t)$ be a shorthand for $\mathrm{exp}(tU)$. Let $X\in \mathfrak{sl}_2(\R)$ be $\begin{pmatrix}\frac{1}{2} & 0 \\ 0 & -\frac{1}{2} \end{pmatrix}$, $\theta=\begin{pmatrix}0 & -\frac{1}{2} \\ \frac{1}{2} & 0 \end{pmatrix}$ and $Y=e^{(\pi/2)ad_{\theta}}X$, so that $Y(t)=\theta(\pi/2)X(t)\theta(-\pi/2)$. We consider $\SL_2(\R)$ as a subset of $M_2(\R)$, then we can add or subtract elements of $\SL_2(\R)$ from each other, and also take the matrix operator norm. And here is the theorem. 

\begin{thm}\label{mm}
Suppose $(a_i)_{i=1}^{n},(b_i)_{i=1}^{n},(a'_i)_{i=1}^{n},(b'_i)_{i=1}^{n}$ are sequences of complex numbers, and $A,B,\ep$ are positive real numbers such that $\ep<\min(1/A,1/\ep)$,
\begin{equation*}
    \sum\limits_{i=1}^{n}|a_i|e^{|b_i|}\leq B,
\end{equation*}
and for all $i$,
\begin{enumerate}

\item $2|a_i|e^{|\re(b_i)|+1}\leq A$;

\item $|b_i-b'_i|<\ep$, and 

\item $|a_i-a'_i|<\ep|a_i|$.
\end{enumerate}
Then
\begin{equation*}
    \left\Vert\sum\limits_{i=1}^{n}Y(b_i)X(a_i)Y(-b_i)-\sum\limits_{i=1}^{n}Y(b'_i)X(a'_i)Y(-b'_i)\right\Vert\leq12e^{A+2B}B\ep.
\end{equation*}
\end{thm}

\subsection{Frames and distortion} 

We first want to introduce the distance and distortion in $\mathcal{F}\HH^3$. Fix a left-invariant metric $d$ on $\mathrm{Isom}(\HH^3)$, and for any $g\in\mathrm{Isom}(\HH^3)$, let $d(g)=d(1,g)$.

For $u,v\in\mathcal{F}\HH^3$, there is a unique $u\to v\in\mathrm{Isom}(\HH^3)$ such that $u\cdot(u\to v)=v$. Provided $X\subset\mathcal{F}\HH^3$ and a map $\tilde{e}:X\to\mathcal{F}\HH^3$, we say $\tilde{e}$ has $\ep$-$bounded$ $distortion$ $to$ $distance$ $D$ if
\begin{equation*}
    d(u\to v,\tilde{e}(u)\to\tilde{e}(v))<\ep,
\end{equation*}
for all $u,v\in X$ with $d(u,v)<D$ (where $d(u,v)=d(u\to v)$).

Given $\tilde{e}_1,\tilde{e}_2:X\to\mathcal{F}\HH^3$, we say that $\tilde{e}_1$ and $\tilde{e}_2$ are $\ep$-$related$ if
\begin{equation*}
    d(\tilde{e}_1(u),\tilde{e}_2(u))<\ep
\end{equation*}
for all $u\in X$. Here are some observations.

\begin{lem}
For all $D$, there exists $D'$ such that if $U\in\SL(2,\C)$ and $\Vert U\Vert<D$, then $d(U)<D'$. This is also true if $\Vert\cdot\Vert$ and $(\cdot)$ are interchanged.
\end{lem}

\begin{lem}
For all $D,\ep$ there exists $\de$ such that if $U,V\in\SL(2,\C)$, $\Vert U\Vert,\Vert V\Vert<D$ and $\Vert U-V\Vert<\de$, then $d(U-V)<\ep$.
\end{lem}

\begin{lem}
For all $\ep,D,k$ there exists $\de$ such that if $u_0,\cdots,u_k,v_0\cdots,v_k\in\mathcal{F}]\HH^3$, and
\begin{equation*}
    d(u_i\to u_{i+1},v_i\to v_{i+1})<\de,
\end{equation*}
and
\begin{equation*}
    d(u_i\to u_{i+1})<D,
\end{equation*}
then
\begin{equation*}
    d(u_0\to u_k,v_0\to v_k)<\ep.
\end{equation*}
\end{lem}

Here are some notations about frames in a 3-manifold $M^3$. Suppose $\gamma$ is an oriented geodesic in $M$, and then any unit normal vector $v$ to $\gamma$ determines a unique 3-frame $q,w,v$ in $M$, where $w$ is the unit tangent vector to $\gamma$ and is positive oriented. And we call this the associated 3-frame for $v$ with respect to $\gamma$. We denote by $\mathcal{F}(\gamma)$ the set of all such frames. If $\gamma$ is unoriented, then we let $\mathcal{F}(\gamma)=\mathcal{F}(\gamma^+)\cup\mathcal{F}(\gamma^-)$, where $\gamma^+$ and $\gamma^-$ are two possible oriented versions of $\gamma$. For a pair of pants $Q$ in $M$, let $\partial^\mathcal{F}Q$ be the union of $\mathcal{F}(\gamma)$ for all oriented $\gamma\in\partial Q$. And we also let $\hat{\partial}^{\mathcal{F}}(Q)$ denote the union of the associated 3-frames for the unique slow and constant vector field on each boundary of $Q$, which is determined by the feet of the short orthogeodesics.

To briefly introduce the logic of this section, suppose that we have $(R_i,\ep_i)_{i=1}^3$-good pants $Q_i$ glued along boundaries, satisfying that each boundary component of one $C_i$ is geometrically identified with some boundary component of another $C_j$, and then to form an assembly $\mathcal{A}$. A $perfect$ $pants$ of an $(R_i,\ep_i)_{i=1}^3$-good pants is a pants with cuff length $R_1,R_2$ and $R_3$, and we will construct a perfect model $\mathcal{\hat{A}}$ for $\mathcal{A}$ later which provide a perfect pants $\hat{Q}_i$ for each $Q_i$ in $\mathcal{A}$ along with a map $h_i:\hat{Q}_i\to M$ sending $\partial\hat{Q}_i$ to $\partial Q_i$ up to homotopy through such maps. We can also get the gluings of $\hat{Q}_i$'s and the geometric and isometric identifications of the boundary components of $\hat{Q}_i$'s. Therefore we obtain a surface $S_{\hat{\mathcal{A}}}$ with a quasi-Fuchsian structure and a homotopy class of maps $h:S_{\hat{\mathcal{A}}}\to M$ which send each $\hat{Q}_i$ to $Q_i$ in $M$. Then we lift this map to $\tilde{h}:\HH^3\to\HH^3$, where the homotopy class of $h$ determines the relationship of each lift of $\hat{Q}_i$ and the corresponding lift of $Q_i$ and the same for boundary geodesics. Now suppose we have maps $e:\partial^\mathcal{F}(\hat{Q}_i)\to\partial^\mathcal{F}(Q_i)$ that maps frames over each boundary geodesic of $\hat{Q}_i$ to frames over the corresponding boundary geodesic of $Q_i$. Then we can use $\hat{h}$ to get a canonical lift $\tilde{e}$ of $e$ to $\partial^\mathcal{F}(\tilde{\hat{\mathcal{A}}})$. And we say $e:\partial^\mathcal{F}(\hat{\mathcal{A}})\to\partial^\mathcal{F}(\mathcal{A})$ has $\ep$-distortion at distance $D$ if and only if $\tilde{e}$ does.

\subsection{Sequences of geodesics}\label{7.3}

In this subsection we will prove a theorem which is based on Theorem \ref{mm}. We first recall some definitions from \cite{KW21}.

A $linear$ $sequence$ $of$ $geodesics$ in $\HH^2$ is a sequence $(\gamma_i)_{i=0}^n$ of disjoint geodesics in $\HH^2$ such that each one of them separates those before it and those after it. The geodesics are oriented such that those following a given geodesic are to the left of that geodesic. For a linear sequence $(\gamma_i)$ of geodesics and given $x_0\in\gamma_0$, we inductively define $x_{i}\in\gamma_i$ such that $x_{i+1}$ and $x_i$ are related by the unique orientation-preserving isometry sending $\gamma_i$ to $\gamma_{i+1}$. And we say that the $x_i$'s form a $homologous$ $sequence$ of points on the $\gamma_i$'s. Similarly we can define a homologous sequence of associated frames since the associated 2-frame to $\gamma_0$ at each point $x_0\in\gamma_0$ is unique determined. A $semi$-$linear$ $sequence$ $of$ $geodesics$ in $\HH^3$ is a sequence $(\gamma_i)_{i=0}^n$ satisfying that each pair of geodesics are disjoint and have a common orthogonal. We can also define a homologous sequence of points and associated frames over a semi-linear sequence.

For two (semi-)linear sequences $(\gamma_i)_{i=0}^n$ and $(\gamma'_i)_{i=0}^n$, let $\eta_i$ be the common orthogonal to $\gamma_i$ and $\gamma_{i+1}$ with orientation from $\gamma_i$ to $\gamma_{i+1}$, $u_i$ be the signed complex distance from $\gamma_i$ to $\gamma_{i+1}$ and $v_i$ be the signed complex distance along $\gamma_i$ from $\eta_{i-1}$ to $\eta_i$. And we likewise define $\eta'_i$, $u'_i$ and $v'_i$ for $(\gamma'_i)$. Furthermore we can define a map $e:\mathcal{F}(\gamma_0)\cup\mathcal{F}(\gamma_n)\to\mathcal{F}(\gamma'_0)\cup\mathcal{F}(\gamma'_n)$ such that $e:\mathcal{F}(\gamma_0)\to\mathcal{F}(\gamma'_0)$ and $e:\mathcal{F}(\gamma_n)\to\mathcal{F}(\gamma'_n)$ are isometric embeddings and $e$ maps the foot of $\eta_0$ on $\gamma_0$ to the foot of $\eta'_0$ on $\gamma'_0$, and the same holds for the foot of $\eta_{n-1}$ on $\gamma_n$.

We say that two sequences $(\gamma_i)_{i=0}^n$ and $(\gamma'_i)_{i=0}^n$ are $(R,B,\ep,B^-,B^+)$-$well$-$matched$, if the following properties hold for each $i$:
\begin{enumerate}
\item $B^-<\re(v_i)<B^+$;

\item $|v'_i-v_i|<\frac{B^-\ep}{2R}$;

\item $B^{-1}<|u_i|e^{R/2}<B$;

\item $|u'_i-u_i|<\ep|u_i|$.
\end{enumerate}

And we say that a sequence $(\gamma_i)_{i=0}^{n}$ is $(R,B,B^-,B^+,K)$-related to another sequence $(\gamma'_i)_{i=0}^{n}$, if there exists a $K$-quasiconformal map from $\C$ to itself sending the attracting and repelling endpoints of $\gamma_i$ to the corresponding points of $\gamma'_i$ for each $i$, and:

\begin{enumerate}
\item $B^{-1}<|u_i|e^{R/2},|u'_i|e^{R/2}<B$;

\item $B^-<\re(v_i),\re(v'_i)<B^+$;

\item $|\re(v_i)-\re(v'_i)|<1/R$.

\end{enumerate}

Here is our theorem.

\begin{thm}\label{7.5}
For any $B>0,D>0,B^+>B^->0$ and $K>1$, there exists $R_0>0$, $\ep_0$ and $C$ such that when $R>R_0$ and $0<\ep<\ep_0$, the following holds. Suppose $(\gamma_i)_{i=0}^n$ and $(\gamma'_i)_{i=0}^n$ are $(R,B,\ep,B^-,B^+)$-well-matched, and $(\gamma_i)_{i=0}^n$ is $(R,B,B^-,B^+,K)$-related to a linear sequence $(\gamma''_i)_{i=0}^{n}$. Then the map $e:\mathcal{F}(\gamma_0)\cup\mathcal{F}(\gamma_n)\to\mathcal{F}(\gamma'_0)\cup\mathcal{F}(\gamma'_n)$ has $C\ep$-bounded distortion at distance $D$.
\end{thm}

Before proving this theorem, we prove some lemmas.

\begin{lem}\label{7.6}
For any $K>1$, there exists $t=t(K)$ such that the following holds. Suppose $f$ is a $K$-quasiconformal mapping from $\C$ to $\C$ fixing $0$. Then for any $z_1,z_2\in\C$ with $|z_1|<|z_2|$, we have $|f(z_1)|<t|f(z_2)|$.
\end{lem}

\begin{proof}
Without loss of generality, we can assume $z_2=f(z_2)=1$. Thus we only need to prove that if $f$ is $K$-quasiconformal, then for any $|z|<1$, we have $|f(z)|<t(K)$. First we fix $|z|<1$ with $z\neq0$.

Let $\mu$ be the Beltrami coefficient of $f$. Then $\|\mu\|_{\infty}=\frac{K-1}{K+1}$. Now for any $s$ in the unit disk, let 
\begin{equation*}
    \mu_s=s\cdot\frac{\mu}{\|\mu\|_{\infty}},
\end{equation*}
and $f_s$ be the quasiconformal automorphism of $\C$ fixing 0 and 1 with Beltrami coefficient $\mu_s$. Hence we know $f_0=id$ and $f_{\|\mu\|_{\infty}}=f$. Let $g(s)=f_s(z)$, then $g$ is holomorphic on the unit disk and $g(s)\neq0,1$. Therefore by Theorem VI.19 in \cite{tsu59}, there is an absolute constant $c$, such that
\begin{equation*}
    |g(s)|\leq \exp(\frac{c\ln(|g(0)|+2)}{1-|s|}).
\end{equation*}
By $|g(0)|=|f_0(z)|=|z|<1$, we have
\begin{equation*}
    |f(z)|=|f_{\|\mu\|_{\infty}}(z)|=|g(\|\mu\|_{\infty})|<\exp(\frac{c\ln3}{1-\|\mu\|_{\infty}}).
\end{equation*}
So $|f(z)|$ is bounded by a constant only depending on $K$, and the lemma is proved.
\end{proof}

For a oriented geodesic $\gamma\in\HH^3$ and a map $f$ from $\partial\HH^3$ to itself, we denote by $[f(\gamma)]$ the geodesic in $\HH^3$ determined by the image of the endpoints of $\gamma$ under $f$, and the orientation of $[f(\gamma)]$ is determined by $f$ and the orientation of $\gamma$.

\begin{lem}\label{7.7}
For any $K>1,D>0$, there exists $C(K,D)>0$ such that the following holds: Suppose $\gamma_1$ and $\gamma_2$ are two disjoint geodesics in $\HH^2$ and $f$ is a $K$-quasiconformal mapping from $\C$ to itself, and we let $\gamma'_i=[f(\gamma_i)]$ for $i=1,2$. Let $\eta$ be the common orthogonal of $\gamma'_1$ and $\gamma'_2$ and $A_i\in\gamma'_i$ be two points. If $d(A_1, A_2)<D$, then 
\begin{equation*}
    |\mathbf{d}(A_1, \eta)-\mathbf{d}(A_2, \eta)|<C(K,D),
\end{equation*}
where $\mathbf{d}(A_i,\eta)$ is the signed distance along $\gamma'_i$.
\end{lem}

\begin{proof}

We first normalize $f$ such that $f$ fixes $0,1$ and $\infty$ and we assume that the endpoints of $\gamma_1$ are $0$ and $\infty$ and the endpoints of $\gamma_2$ are $1$ and $x$ for some $1\neq x>0$. Without loss of generality, we assume $0<x<1$. Let $d=d(\gamma_1,\gamma_2)$, which is actually real and we have
\begin{equation}\label{7.7.0}
    d=\ln\left(\frac{1+\sqrt{x}}{1-\sqrt{x}}\right).
\end{equation}
Similarly, we have
\begin{equation}\label{7.7.01}
    d_{\C}(\gamma'_1,\gamma'_2)=\ln\left(\frac{1+\sqrt{f(x)}}{1-\sqrt{f(x)}}\right)\ \mathrm{or}\ -\ln\left(\frac{1+\sqrt{f(x)}}{1-\sqrt{f(x)}}\right). 
\end{equation}
So there exists $\ep_1$ such that when $|z|<\ep_1$, we have
\begin{equation}\label{7.7.02}
    |\sqrt{z}|<\left|\ln\left(\frac{1+\sqrt{z}}{1-\sqrt{z}}\right)\right|<4|\sqrt{z}|.
\end{equation}

It's also well-known that given $K>1$, for a normalized $K$-quasiconformal mapping $f$ there exist $\ep_2>$ and $M>0$ such that if $|x|<\ep_2$, then 
\begin{equation}\label{7.7.1}
    |f(x)|<M|x|^{1/K}.
\end{equation}

Therefore by (\ref{7.7.0}), (\ref{7.7.01}), (\ref{7.7.02}) and (\ref{7.7.1}), there exists $\ep_0>0$ such that when $0<d>\ep_0$, we have
\begin{equation}\label{7.7.2}
    \begin{aligned}
        |d_{\C}(\gamma'_1,\gamma'_2)|&=\left|\ln\left(\frac{1+\sqrt{|f(x)|}}{1-\sqrt{|f(x)|}}\right)\right|<4\sqrt{|f(x)|}\\
        &<4\sqrt{M x^{1/K}}< 4\sqrt{M}\left(\ln\left(\frac{1+\sqrt{x}}{1-\sqrt{x}}\right)\right)^{1/K}\\
        &=4\sqrt{M}d^{1/K}\leq 4\sqrt{M}\ep_0^{1/K},
    \end{aligned}
\end{equation}
which is a constant depending only on $K$. If $\mathbf{d}(A_1,\eta)\mathbf{d}(A_2,\eta)>0$, then by triangle inequality, we have
\begin{equation*}
    |\mathbf{d}(A_1,\eta)-\mathbf{d}(A_2,\eta)|\leq d(A_1, A_2)<D.
\end{equation*}
If $\mathbf{d}(A_1,\eta)\mathbf{d}(A_2,\eta)<0$, then by Lemma \ref{iod} and (\ref{7.7.2}), there is a constant $C_1(K,D)$ depending on $K$ and $D$ such that 
\begin{equation*}
    |\mathbf{d}(A_1,\eta)-\mathbf{d}(A_2,\eta)|=|\mathbf{d}(A_1,\eta)|+|\mathbf{d}(A_2,\eta)|<d(A_1, A_2)+C_1(K,D)<D+C_1(K,D).
\end{equation*}
Thus we always have
\begin{equation}\label{7.7.8}
    |\mathbf{d}(A_1,\eta)-\mathbf{d}(A_2,\eta)|<D+C_1(K,D).
\end{equation}

When $d>\ep_0$, we first prove a statement that there exists $\delta$ depending only on $K$, such that
\begin{equation*}
    |d_{\C}(\gamma'_1,\gamma'_2)|>\de
\end{equation*}
and
\begin{equation*}
    |d_{\C}(\gamma'_1,\gamma'_2)-i\pi|>\de.
\end{equation*}
By $d>\ep_0$ and (\ref{7.7.0}), there exists $\de_1$ such that $x>\de_1$. Then we apply (\ref{7.7.1}) to $f^{-1}$, we have
\begin{equation*}
    \de_1<|x|<M|f(x)|^{1/K}\ \mathrm{or}\ |f(x)|>\ep_1.
\end{equation*}
Thus $|f(x)|>\min\{(\de_1/M)^K,\ep_1\}$. Since $0<x<1$, $|f(x)|$ is also bounded above by Lemma \ref{7.6}. By $$d_{\C}(\gamma'_1, \gamma'_2)=\ln\left(1+\frac{2\sqrt{f(x)}}{1-\sqrt{f(x)}}\right),$$
we know there exists $\de$ such that $|d_{\C}(\gamma'_1,\gamma'_2)|>\de$. To prove the second inequality, we only need to reverse the orientation of $\gamma_2$ and then it follows from the first inequality.

Now if $d_{\R}(\gamma'_1,\gamma'_2)>\de/2$, then by Lemma \ref{iod}, there exists a constant $C_2(K,\de)$ such that
\begin{equation}\label{7.7.10}
    |\mathbf{d}(A_1,\eta)|+|\mathbf{d}(A_2,\eta)|<d(A_1, A_2)+C_2(K,\de).
\end{equation}
If $d_{\R}(\gamma'_1,\gamma'_2)\leq\de/2$, then $\de/2<\im(d_{\C}(\gamma'_1,\gamma'_2))<\pi-\de/2$. Then by Lemma \ref{iod} again, there exists $C_3(K,\de)$ such that 
\begin{equation}\label{7.7.11}
    |\mathbf{d}(A_1,\eta)|+|\mathbf{d}(A_2,\eta)|<d(A_1, A_2)+C_3(K,\de).
\end{equation}

Since $\de$ only depends on $K$, then by (\ref{7.7.8}), (\ref{7.7.10}) and (\ref{7.7.11}), there exists a constant $C(K,D)$ such that
\begin{equation*}
    |\mathbf{d}(A_1, \eta)-\mathbf{d}(A_2, \eta)|<C(K,D).\qedhere
\end{equation*}

\end{proof}

The next two lemma are Lemma A.10 and Lemma A.13 in \cite{KW21}.

\begin{lem}\label{7.16}
For all $\ep,D,k$, there exists $\de$: When $u_0,\dots,u_k,v_0,\dots,v_k\in\mathcal{F}{\HH^3}$, and
\begin{equation*}
    d(u_i\to u_{i+1}, v_i\to v_{i+1})<\de,
\end{equation*}
and 
\begin{equation*}
    d(u_i\to u_{i+1})<D,
\end{equation*}
then
\begin{equation*}
    d(u_0\to u_k,v_0\to v_k)<\ep.
\end{equation*}
\end{lem}

\begin{lem}\label{7.17}
Suppose $(\gamma_i)_{i=0}^n$ is a linear sequence, with $(u_i)$ and $(v_i)$ defined as at the beginning of this subsection. Let $D=d(\gamma_0,\gamma_n)$, and suppose that $u_0,u_{n-1}<1$. Then
\begin{equation*}
    \left|\sum\limits_{i=1}^{n-1} v_i\right|\leq D+2\ln{D}-\ln{u_0}-\ln{u_{n-1}}+3.
\end{equation*}
\end{lem}

Now we can prove Theorem \ref{7.5}.

\begin{proof}[Proof of Theorem \ref{7.5}]

Suppose $x\in\mathcal{F}(\gamma_0)$ and $y\in\mathcal{F}(\gamma_n)$ such that $d(x,y)<D$, and then in particular, $d(\gamma_0,\gamma_n)<D$. By Theorem 1.2 in \cite{Shiga05}, there is an absolute constant $A$ and a constant $C_K$ only depending on $K$, such that $d(\gamma''_0,\gamma''_n)<AKD+C_K:=D'$. Then by Lemma \ref{7.17}, we have
\begin{equation}\label{7.5.1}
    (n-1)B^-<\sum\limits_{i=1}^{n-1}v''_i<D'+2\ln D'+R+2\ln B+3=R+C_1(B,D,K)<2R,
\end{equation}
when $R>C_1(B,D,K)$. Let $(x_i)$ be the homologous sequence of frames for $(\gamma_i)$ with $x_0=x$, and let $(x'_i)$ be the same for $(\gamma'_i)$ with $x'_0=x'=e(x)$.

We want to use Theorem \ref{mm} to control $d(x_0\to x_n, x'_0\to x'_n)$. Let $a_i=u_i$ and $b_i$ be complex numbers such that $\mathrm{foot}_{\gamma_{i+1}}\gamma_i=x_i Y(b_i)$. Then we notice that $x_{i+1}=x_i Y(b_i) X(a_i) Y(-b_i)$ and $b_{i+1}=b_i+v_{i+1}$, and we have similar results for $x'_i$ and $b'_i$. Thus $b'_0=b_0$ by the definition of the map $e$, and 
\begin{equation*}
    b_i=b_0+\sum\limits_{j=1}^{i}v_j,
\end{equation*}
and the same for $b'_i$ and $v'_i$. Hence by (\ref{7.5.1}) we have
\begin{equation*}
    |b'_i-b_i|=\left|\sum\limits_{j=1}^{i}(v'_i-v_i)\right|<\frac{(n-1)B^- \ep}{2R}<\ep,
\end{equation*}
by Condition 2 of well-matched. So we verify Condition 2 of Theorem \ref{mm}. On the other hand, Condition 3 of Theorem \ref{mm} directly follows from Condition 4 of well-matched.

Now we want to control $|a_i|e^{|\re(b_i)|}$ to satisfy the remaining conditions in Theorem \ref{mm}. Since we have $a_i$'s are about $e^{-R/2}$ in size and $B^-<\re(v_i)<B^+$, so we only need to estimate the largest $|a_i|e^{|\re(b_i)|}$, which is actually either the first or the last, in order to control the sum of all $|a_i|e^{|\re(b_i)|}$. 

To control $|a_0|e^{|\re(b_0)|}$, we only need to give an upper bound of $d(x,\gamma_1)$. Let $f$ be the $K$-quasiconformal mapping sending the endpoints of $(\gamma'_i)_{i=1}^n$ to the endpoints of $(\gamma''_i)_{i=1}^n$ and $\tilde{f}:\HH^3\to\HH^3$ be the $(K',C')$-quasi-isometric extension of $f$, where $K'$ and $C'$ are constants depend on $K$. Let $\tilde{g}:\HH^3\to\HH^3$ be the approximate inverse of $\tilde{f}$, and then $\tilde{g}$ is a quasi-isometric extension of $f^{-1}$. And we also know that $\tilde{f}(\gamma_i)$ is within distance $C'_1$ of $\gamma''_i$ and $\tilde{g}(\gamma''_i)$ is within distance $C'_1$ of $\gamma_i$ for some constant $C'_1$ depending only on $K$. Now let $x''$ and $y''$ be the projection of $\tilde{f}(x)$ and $\tilde{f}(y)$ on $\gamma''_0$ and $\gamma''_n$, respectively. And we have
\begin{equation*}
    d(x'',y'')<d(x'',\tilde{f}(x))+d(\tilde{f}(x),\tilde{f}(y))+d(\tilde{f}(y),y'')<2C'_1+K'D+C'.
\end{equation*}
Therefore
\begin{equation*}
    d(x'',\gamma''_2)\leq d(x'',\gamma''_n)\leq d(x'',y'')<2C'_1+K'D+C'.
\end{equation*}
Hence 
\begin{equation*}
    d(\tilde{f}(x),\gamma''_2)\leq d(\tilde{f}(x),x'')+d(x'',\gamma''_2)<3C'_1+K'D+C'.
\end{equation*}
Since we also have $d(x,\tilde{g}\circ\tilde{f}(x))<C'$, so
\begin{equation*}
\begin{aligned}
    d(x,\gamma_2)&<d(x,\tilde{g}(\gamma''_2))+C'_1<d(\tilde{g}\circ\tilde{f}(x),\tilde{g}(\gamma''_2))+C'+C'_1\\
    &<K'(3C'_1+K'D+C')+C'+C'_1,
\end{aligned}
\end{equation*}
which is a constant depending on $K$ and $D$.

Now we want to control $|a_{n-1}|e^{|\re(b_{n-1})|}$. By the same method as above, we can have an upper bound of $d(y,\gamma_{n_1})$, so we only need to estimate $d(x_n,y)$. Let $s=\frac{1}{2}\sum_{j=1}^{n-1}v_j$, then by (\ref{7.5.1}) and Condition (3) for $(R,B,B^-,B^+,K)$-related, we know
\begin{equation}\label{7.5.4}
\begin{aligned}
    \re(s)&=\frac{1}{2}\sum\limits_{j=1}^{n-1}\re(v_j)<\frac{1}{2}\sum\limits_{j=1}^{n-1}\left(\re(v''_j)+\frac{1}{R}\right)
    <\frac{R}{2}+\frac{1}{2}C_1(B,D,K)+\frac{n-1}{R}\\
    &<\frac{R}{2}+\frac{1}{2}C_1(B,D,K)+\frac{2}{B^-}=\frac{R}{2}+C_2(B,D,K,B^-).
\end{aligned}
\end{equation}
Let $a_0$ be the base point of the frame $\mathrm{foot}_{\gamma_{1}}\gamma_0 Y(s)$, and $a_i$'s is a homologous sequence determined by $a_0$. Then $d(a_0, \mathrm{foot}_{\gamma_{1}}\gamma_0)=d(a_n,\mathrm{foot}_{\gamma_{n-1}}\gamma_n)=\re(s)$. Then by (\ref{7.5.4}) and Condition 3 of well-matched, we know there exists a constant $C_3(B,D,K,B^-)$ such that
\begin{equation}\label{7.5.5}
    d(a_0,a_n)\leq\sum\limits_{i=0}^{n-1}d(a_i,a_{i+1})<C_3(B,D,K,B^-).
\end{equation}
Let $\eta$ be the common orthogonal between $\gamma_0$ and $\gamma_n$, then by \ref{7.5.5} and Lemma \ref{7.7}, we have
\begin{equation}\label{7.5.6}
    |\mathbf{d}(a_0, \eta)-\mathbf{d}(a_n, \eta)|<C_4(B,D,K,B^-),
\end{equation}
for some constant $C_4(B,D,K,B^-)$. Similarly, there is a constant $C_5(K,D)$ such that 
\begin{equation}\label{7.5.7}
    |\mathbf{d}(x, \eta)-\mathbf{d}(y, \eta)|<C_5(K,D).
\end{equation}
By the definition of homologous sequence, we know $\mathbf{d}(x, a_0)=\mathbf{d}(x_n, a_n)$ as singed distance. Thus
\begin{equation*}
\begin{aligned}
    d(x_n,y)&=|\mathbf{d}(x_n, \eta)-\mathbf{d}(y, \eta)|\\
    &=|(\mathbf{d}(x_n, \eta)-\mathbf{d}(x, \eta))+(\mathbf{d}(x, \eta)-\mathbf{d}(y, \eta))|\\
    &=|(\mathbf{d}(x_n, \eta)+d(a_n,x_n))-(d(a_0,x)+\mathbf{d}(x, \eta))+(\mathbf{d}(x, \eta)-\mathbf{d}(y, \eta))|\\
    &=|(\mathbf{d}(a_n, \eta)-\mathbf{d}(a_0, \eta))+(\mathbf{d}(x, \eta)-\mathbf{d}(y, \eta))|\\
    &<|\mathbf{d}(a_n, \eta)-\mathbf{d}(a_0, \eta)|+|\mathbf{d}(x, \eta)-\mathbf{d}(y, \eta)|\\
    &<C_4(B,D,K,B^-)+C_5(K,D).
\end{aligned}
\end{equation*}
And hence all the conditions in Theorem \ref{mm} are satisfied.

Next we consider $d(x_n\to y, x'_n\to y')$. Since $x_n$ and $y$ are frames on $\gamma_n$, and we also know $|b_{n-1}-b'_{n-1}|<\ep$ and $d(y,\eta_{n-1})=d(y',\eta'_{n-1})$ as signed distance, so $d(x_n\to y, x'_n\to y')<\ep$.

Now we consider $x,x_n,y$ and $x',x'_n,y'$, and by Lemma \ref{7.16}, the theorem is proved.

\end{proof}

\subsection{Good assemblies and the perfect model}

In this subsection, we construct the perfect model $\hat{\mathcal{A}}$ for a good assembly $\mathcal{A}$, then bound the distortion of the map $e: \partial^{\mathcal{F}}\hat{\mathcal{A}} \to \partial^{\mathcal{F}}\mathcal{A}$ by some geometric control.

Suppose $Q$ is a $(R_i,\ep)_{i=1}^3$-good pants, and $\hat{Q}$ is its perfect model. That means the lengths of cuffs of $\hat{Q}$ are $R_1,R_2$ and $R_3$ and we also say $\hat{Q}$ is $(R_i)_{i=1}^3$-perfect. We say $Q$ is $\ep$-$compliant$ if for every short orthogeodesic $\eta$ of $Q$ and corresponding orthogeodesic $\hat{\eta}$ in $\hat{Q}$, we have
\begin{equation*}
    |l(\eta)-l(\hat{\eta})|<\ep l(\hat{\eta}).
\end{equation*}
Then we recall and revise the Lemma A.15 in \cite{KW21} for our purpose.

\begin{lem}
For each $\de>0$, there is a universal constant $C$ such that every $(R_i,\ep)_{i=1}^3$-good pants is $C\ep$-compliant with $|R_i-\overline{R}|<\de$ for any $\overline{R}>\de$.
\end{lem}

For a $(R_i,\ep)_{i=1}^3$-good pants $Q$, there is a unique map $e:\hat{\partial}^{\mathcal{F}}\hat{Q}\mapsto\partial^{\mathcal{F}}Q$ satisfying the following three properties:
\begin{enumerate}

\item The map from $\partial\hat{Q}$ to $\partial Q$ induced by $e$ is the restriction (to $\partial \hat{Q}$) of an orientation-preserving homeomorphism from $\hat{Q}$ to $Q$.

\item The induced map is affine (linear) on each component of $\partial\hat{Q}$, and maps each component of $\hat{\partial}^{\mathcal{F}}\hat{Q}$ to the frames determined by a slow and constant turning vector field on $\partial Q$.

\item $e$ maps each foot of $\hat{Q}$ to the corresponding foot of $Q$.
\end{enumerate}
And we say $e$ is $M,\ep$-compliant if $e$ is $\ep$-bounded distortion to distance $M$.

Suppose that $\mathcal{A}$ is a good assmebly and $\hat{\mathcal{A}}$ is a perfect one. Then we say that $e:\partial^{\mathcal{F}}\hat{\mathcal{A}}\to\partial^{\mathcal{F}}\mathcal{A}$ is $M,\ep$-compliant if the following are satisfied:
\begin{enumerate}

\item For each component $\hat{Q}$ of $\hat{\mathcal{A}}$, there is a corresponding $Q$ of $\mathcal{A}$ such that the restriction of $e$ on $\hat{\partial}^{\mathcal{F}}\hat{Q}$ is a $M,\ep$-compliant map to $\partial^{\mathcal{F}}Q$.

\item If $\hat{\gamma}$ is a gluing boundary of $\hat{Q_1}$ and $\hat{Q_2}$ of $\hat{\mathcal{A}}$, and $n_i$ are the frames in $\hat{\partial}^{\mathcal{F}}\hat{Q_i}$ for $i=1,2$ with $n_1,n_2$ sharing the same base point on $\hat{\gamma}$, then
\begin{equation*}
    d(n_1\to n_2,e(n_1)\to e(n_2))<\ep.
\end{equation*}

\end{enumerate}

Before constructing the map from a good assembly to its perfect model, we want to introduce one more definition. We say two complex tuples $(R_i,s_i)_{i=1}^3$ and $(R'_i,s'_i)_{i=1}^3$ are $K$-related for some $K>1$, if there are two genus-2 quasi-Fuchsian groups $\Gamma_1, \Gamma_2$ and a $K$-quasiconformal mapping $f:\hat{\C}\to\hat{\C}$ such that:
\begin{enumerate}
\item $\Gamma$ has a non-separating pants decomposition with cuff length $R_i$ and shear $s_i$, $i=1,2,3$, and the same for $\Gamma'$ and $(R'_i,s'_i)_{i=1}^3$.
\item $f$ conjugates $\Gamma_1$ to $\Gamma_2$ and preserves corresponding the homotopy classes of each corresponding pants decomposition.
\end{enumerate}

\begin{thm}\label{7.23}
For all $M>0,\de>0,B^+>B^->0$ and $K>1$, we can find $C,R_0>0$ such that for all $\overline{R}>R_0$ and $\ep>0$: Suppose $R_i,s_i$ satisfy $|R_i-\overline{R}|<\de$ and $B^-<\re(s_i)<B^+$, and there exist $\overline{R}-\de<R'_i<\overline{R}+\de$ and $B^-<s'_i<B^+$ such that $(R_i,s_i)_{i=1}^3$ is $K$-related to $(R'_i,s'_i)_{i=1}^3$ and $|\re(s_i)-s'_i|<1/\overline{R}$. Let $\mathcal{A}$ be a $(R_i,s_i,\ep)_{i=1}^3$-good assembly, there is a $(R_i,s_i)_{i=1}^3$-perfect assembly $\hat{\mathcal{A}}$ and an $M,C\ep$-compliant map $e:\partial^{\mathcal{F}}\hat{\mathcal{A}}\to\partial^{\mathcal{F}}\mathcal{A}$.
\end{thm}

\begin{proof}
By Theorem \ref{7.5}, and $(R_i,s_i)_{i=1}^3$ being $K$-related to $(R'_i,s'_i)_{i=1}^3$, we can prove that there exists $C_1>0$ such that for each good pants $Q$ of $\mathcal{A}$, we can construct a $M,C_1\ep$-compliant map $e:\hat{\partial}^{\mathcal{F}}\hat{Q}\mapsto\partial^{\mathcal{F}}Q$ where $C_1$ does not depend on $\ep$, as the proof of Theorem A.16 in \cite{KW21}. Now given the whole assembly $\mathcal{A}$, we construct the perfect one as follows: if two good pants are glued along a $(R_i,\ep)$-good curve, then the corresponding two perfect pants are glued along the boundary curve with length $R_i$ and the two feet are joined with shear by $s_i$. Thus Condition 1 of the global map $e$ automatically holds. Condition 2 then follows from the $(R_i,s_i,\ep)_{i=1}^3$-goodness of $\mathcal{A}$, since the base points of $e(n_1)$ and $e(n_2)$ are always within $(B^++1)\ep/R$ of each other, and the difference of bending is at most $\ep$.
\end{proof}

And know we can estimate the distortion of such map.

\begin{thm}\label{7.25}
For all $D>0,B^+>B^->0,\de>0$ and $K>1$ there exist $C,R_0,\ep_0>0$ such that for all $0<\ep<\ep_0$ and $\overline{R}>R_0$: Suppose $R_i,s_i$ satisfy $|R_i-\overline{R}|<\de$ and $B^-<\re(s_i)<B^+$, and there exist $\overline{R}-\de<R'_i<\overline{R}+\de$ and $B^-<s'_i<B^+$ such that $(R_i,s_i)_{i=1}^3$ is $K$-related to $(R'_i,s'_i)_{i=1}^3$ and $|\re(s_i)-s'_i|<1/\overline{R}$, then for any $(R_i,s_i,\ep)_{i=1}^3$-good assembly $\mathcal{A}$, we can find a $(R_i,s_i)_{i=1}^3$-perfect assembly $\hat{\mathcal{A}}$ and a map $e:\partial^{\mathcal{F}}\hat{\mathcal{A}}\to\partial^{\mathcal{F}}\mathcal{A}$ which has $C\ep$-bounded distortion to distance $D$.
\end{thm}

\begin{proof}
Let $\mathcal{A}$ be as given in the statement. Then by Theorem \ref{7.23}, we have a $(R_i,s_i)_{i=1}^3$-perfect model $\hat{\mathcal{A}}$ and a $D,C\ep$-compliant map $e:\partial^{\mathcal{F}}\hat{\mathcal{A}}\to\partial^{\mathcal{F}}\mathcal{A}$, where $C$ does not depend on $\ep$. We will work in the universal cover of $\hat{\mathcal{A}}$, and apply Theorem $\ref{7.5}$ to prove that the lift of $e$ (still denoted by $e$) has bounded distortion. Thus we suppose that $p,q$ are two frames based on the boundary curves of $\hat{\mathcal{A}}$ in the universal cover with $d(p,q)<D$, and let $p$ lies on $\gamma$ and $q$ lies on $\overline{\gamma}$. 

Let $f$ be the $K$-quasiconformal mapping which relates $(R_i,s_i)_{i=1}^3$ and $(R'_i,s'_i)_{i=1}^3$, then $f$ extends to a $(K',C')$-quasi-isometry $\tilde{f}:\HH^3\to\HH^3$ and we can also assume that for each geodesic $\gamma\in\HH^3$, $\tilde{f}(\gamma)$ is always within distance $C'$ from $[\tilde{f}(\gamma)]$. Since $\hat{\mathcal{A}}$ corresponds to a finite index subgroup of a genus-2 quasi-Fuchsian group $\Gamma_1$, then we apply the conjugacy by $f$ and get a finite index subgroup of a genus-2 Fuchsian group $\Gamma_2$ which gives us an $(R'_i,s'_i)_{i=1}^3$-perfect assembly $\hat{\mathcal{A}}'$. Let $\gamma'_0=[\tilde{f}(\gamma)]$ and $\gamma'_n=[\tilde{f}(\overline{\gamma})]$ with lifts of boundary curves $\gamma'_1,\gamma'_2,\dots,\gamma'_{n-1}$ separating them in sequence. Let $\gamma_i=[\tilde{f}^{-1}(\gamma'_i)]$ for $i=0,1,2,\dots,n$, then $\gamma=\gamma_0$ and $\overline{\gamma}=\gamma_n$. And we define $(\eta_i),(u_i),(v_i)$ and $(\eta'_i),(u'_i),(v'_i)$ as in Section \ref{7.3}.

Let $p'\in\gamma'_0$ and $q'\in\gamma'_n$ such that $d(\tilde{f}(p),p')=d(\tilde{f}(p),\gamma'_0)$ and $d(\tilde{f}(q),q')=d(\tilde{f}(q),\gamma'_n)$. Also let $\tilde{g}$ be the approximate inverse of $\tilde{f}$, which is a quasi-isometric extension of $f^{-1}$. Suppose the geodesic segment $p' q'$ intersects with $\gamma'_i$ at $p'_i$ for $i=1,2,\dots,n-1$ and let $p'=p'_0$ and $q'=p'_n$. Now we know $\eta_i,\eta'_i$ is a short orthogeodesic between two cuffs of a pair of pants, so we know there is a constant $C_1$ such that $C_1^{-1}r^{-\overline{R}/2}<|u_i|,u'_i<C_1 e^{-\overline{R}/2}$. And we also know that for each $i$, there exists $j_i\in\{1,2,3\}$ such that $v_i\equiv s_{j_i}(\mathrm{mod}\ R_{j_i})$ (as complex numbers) and $v'_i\equiv s'_{j_i}(\mathrm{mod}\ R'_{j_i})$. By Lemma \ref{7.17}, there exists a constant $C_2$ such that if $d(\gamma'_{i-1},\gamma'_{i+1})<C_2$, then $v'_i=s'_{j_i}$. Moreover since $\hat{\mathcal{A}}$ and $\hat{\mathcal{A}}'$ are related by $f$, we know
\begin{equation*}
    \frac{v_i-s_{j_i}}{R_{j_i}}=\frac{v'_i-s'_{j_i}}{R'_{j_i}}.
\end{equation*}
So $v_i=s_{j_i}$ if and only if $v'_i=s'_{j_i}$.

Now we define a $run$ be an interval $\Z\cap [x,y]$ such that $v'_i=s'_{j_i}$ for all $x<i<y$. We allow the case that $y=x+1$ which is a trivial run. Then we can find integers $0=x_1<x_2<\dots<x_k=n$ such that each $[x_i,x_{i+1}]$ is a maximal run, and then the union of these intervals cover $\Z\cap[0,n]$. Here a maximal run means it cannot be extended to contain more integers. And we also know $k$ is bounded in terms of $D$ and $K$, since $d(\gamma'_{x_i-1},\gamma'_{x_i+1})$ is bounded below and $d(p_0,p_n)<K'D+2C'$ by $\tilde{f}$ quasi-isometric. Similarly we have $d(p_i,p_{i+1})<K'D+2C'$, and thus there is a constant $C_3$ such that $d(p_i,\eta_{i-1}),d(p_i,\eta_i)<\overline{R}/2+C_3$ which is the case $n=1$ of Theorem \ref{7.5}.

We notice that for each $i$, $(\gamma_i)_{i=x_i}^{x_{i+1}}$ and $(\gamma''_i)_{i=x_i}^{x_{i+1}}$ are $(R,C_1,\ep,B^-,B^+)$-well-matched; Condition 1 and Condition 3 satisfied by the definition of a run, and Condition 2 and Condition 4 hold from $(R_i,s_i,\ep)$-goodness of $\mathcal{A}$. And $(\gamma_i)_{i=x_i}^{x_{i+1}}$ is $(R,C_1,B^-,B^+,K)$-related to $(\gamma'_i)_{i=x_i}^{x_{i+1}}$ since $(R_i,s_i)_{i=1}^3$ is $K$-related to $(R'_i,s'_i)_{i=1}^3$ and $|\re(s_i)-s'_i|<1/\overline{R}$.

Let $p_{x_i}\in\gamma_{x_i}$ be the point closet to $\tilde{g}(p'_{x_i})$ for $i=1,2,\dots, n-1$, and let $p_0=p$ and $p_n=q$. Let $\al_i$ be the frame lifted from $\hat{\partial}^{\mathcal{F}}\hat{Q}$ where two boundary curves of $\hat{Q}$ lifts to $\gamma_{x_i}$ and $\gamma_{x_i+1}$, and $\be_i$ be the frame lifted from $\hat{\partial}^{\mathcal{F}}\hat{Q}$ where two boundary curves of $\hat{Q}$ lifts to $\gamma_{x_i}$ and $\gamma_{x_i-1}$.

Now we make the following claims:
\begin{enumerate}
\item $d(\al_i\to\be_{i+1},e(\al_i)\to e(\be_{i+1}))<\ep$;
\item $d(\be_i\to\al_i,e(\be_i)\to e(\al_i))<\ep$.
\end{enumerate}
The second one directly follows from the $D,C\ep$-compliance of $e$. And for the first one, we know the map in Theorem \ref{7.5} is $\ep$-related to $e$ on the relevant parts of $\mathcal{F}(\gamma_{x_i})$ and $\mathcal{F}(\gamma_{x_{i+1}})$ since $d(p_j,\eta_{j-1}),d(p_j,\eta_j)<\overline{R}/2+C_3$ for all $j$. Then Claim 1 follows from Theorem \ref{7.5}.

And then the Theorem follows by Lemma \ref{7.16}.
\end{proof}

\subsection{Extensions}

We recall that $X'\subset X$ is called $A$-$dense$ in a metric space $X$ if $\mathcal{N}_{A}(X')=X$.

\begin{thm}\label{7.26}
For all $A$ there exist $B,K$ such that for all $\de,\Omega$ there exists $\ep$: Suppose $\Lambda$ is a $K$-quasi-circle in $\hat{\C}$ and $U\subset \mathcal{F}(C(\Lambda))$ is $A$-dense, here $C(\Lambda)$ is the convex hull of $\Lambda$. If $e:U\to\mathcal{F}(\HH^3)$is a map having $\ep$-bounded distortion to distance $B$. Then $e$ is a $K$-quasi-isometric embedding, and $e$ extends to $\hat{e}:\Lambda\to\partial(\HH^3)$ to be a $(\Omega,1+\de)$-quasi-symmetric embedding.
\end{thm}

\begin{rem}
Here by $\hat{e}$ being $(\Omega,1+\de)$-quasi-symmetric, we mean that for any quadruple $(z_1,z_2,z_3,z_4)$ of four distinct points with its cross ratio in a compact set $\Omega\subset\C$, we have 
\begin{equation*}
    |[\hat{e}(z_1),\hat{e}(z_2);\hat{e}(z_3),\hat{e}(z_4)]-[z_1,z_2;z_3,z_4]|<\de.
\end{equation*}
\end{rem}

We want to quote two theorems in \cite{KW21} to prove the above theorem.

\begin{thm}
For all $K,\de$, there exist $K',D$: Suppose $X$ is a path metric space and $Y$ is $\de$-hyperbolic, and $f:X\to Y$ is such that 
\begin{equation*}
    K^{-1} d(x,x')-K<d(f(x),f(x'))<K d(x,x')+K
\end{equation*}
whenever $d(x,x')<D$. Then $f$ is a $K'$-quasi-isometric embedding.
\end{thm}

\begin{thm}
Let $X$ and $Y$ be Gromov hyperbolic, and let $f:X\to Y$ be a quasi-isometric embedding. Then $f$ extends continuously to an embedding $\hat{f}:\partial X\to \partial Y$. Moreover, $\hat{f}$ depends continuously on $f$ with the uniform topology on $\hat{f}$ and the local uniform topology on $f$.
\end{thm}

\begin{proof}[Proof of Theorem \ref{7.26}]

The proof of the first part follows from the proof of Theorem A.19 in \cite{KW21}, by applying the above two theorems. Thus we only need to show that $\hat{e}$ is $(\Omega,1+\de)$-quasisymmetric.

Since we can change $\hat{e}$ by M\"obius transformations in domain and range, we can assume that the quadruple on $\Lambda$ is $(z,-1;1,\infty)$ for some $z\in\C$ with $z\neq1,-1$ and $\hat{e}(-1)=-1$, $\hat{e}(1)=1$ and $\hat{e}(\infty)=\infty$. To show the inequality for cross ratios, we only need to prove that there exists $\ep$ such that
\begin{equation*}
    |\hat{e}(z)-z|<2\de.
\end{equation*}
Suppose to the contrary, such $\ep$ does not exist, then we can take a sequence of maps $e_n$ and $z_n\in\C-\{1,-1\}$ defined on a sequence of $A$-dense set $U_n$ with $1/n$-bounded distortion at distance $B$ and $|\hat{e}_n(z_n)-z_n|>2\de$. Since $\Omega$ is compact, we know $\{z_n\}$ is uniformly bounded. By passing to a subsequence, we have limits $e_\infty$ and $z_\infty\in\C$, and we have $\hat{e}_{\infty}(z)=z$ for any $z$; because $\hat{e}_\infty$ also preserves $-1,1$ and $\infty$ and also is a M\"obius transformation. In particular, $\hat{e}_\infty(e_\infty)=z_\infty$ Contradiction to $|\hat{e}_\infty(z_\infty)-z_\infty|>2\de$!

\end{proof}

Now our last theorem in this subsection is to extend the map on a quasi-circle to a map on the whole Riemann sphere, which is the next theorem.

\begin{thm}\label{7.30}
Given any $K\geq1$ and $M>1$, there exist $\Omega$ compact in $\C$ and $\de>0$ such that: Suppose $\gamma$ is a $K$-quasicircle in $\hat{\C}$ and $g:\gamma\to\hat{\C}$ is $(\Omega,1+\de)$-quasi-symmetric. Then $g$ extends to a $M$-quasiconformal map from $\hat{\C}$ to $\hat{\C}$. Moreover, if $g$ conjugates a group of M\"obius transformation to another such group, then the extension does as well.
\end{thm}

\begin{lem}\label{7.31}
For any $K>1$ and $\ep>0$, there exists $\de>0$ and $\Omega$ compact: Suppose $f:\C\to\C$ is a $K$-quasiconformal mapping and $g:f(S^1)\to\C$ is $(\Omega,1+\de)$-quasisymmetric. Suppose $f$ and $g$ are normalized by $f(-1)=-1=g(-1)$ and $f(1)=1=g(1)$. Then for any $z\in S^1$, we have
\begin{equation*}
    |g(f(z))-f(z)|<\ep.
\end{equation*}
\end{lem}

\begin{proof}
To prove this lemma, we only need to prove that: For any $K>1$ and $\ep>0$, there exists $\de>0$ such that if $f:\R\to\C$ is a $K$-quasi-line and $g:f(\R)\to\C$ is $(\Omega,1+\de)$-quasisymmetric with $f(i)=i=g(i)$ for $i=0,1$, then for any $z\in[0,1]$, we have
\begin{equation*}
    |g(f(z))-f(z)|<\ep.
\end{equation*}

Given $K>1$, we know that there exists $m\in\Z_+$ such that for any $K$-quasiconformal mapping $f:\C\to\C$ with $f(0)=0$ and $f(1)=1$, and $x,y\in\mathbb{D}$ (the unit disc) with $|x-y|\leq1/2^m$, we have $|f(x)-f(y)|\leq1/2$. We also choose $\Omega$ to be the close disk centered at $0$ with radius $R$, where $R$ will be determined later. Now we define
\begin{equation*}
\begin{aligned}
    T_\de(n)=\max\{&|g(f(\frac{a}{2^{nm}}))-f(\frac{a}{2^{nm}})|:0\leq a\leq 2^{nm}, 2^m\nmid a, f\ \mathrm{is}\ K\textup{-quasiconformal},\\
    &g\ \mathrm{is}\ (\Omega,1+\de)\textup{-quasisymmetric},f(0)=0=g(0),f(1)=1=g(1)\},
\end{aligned}
\end{equation*}
for $n\geq1$. We want to use recursion to prove that there exists $\de$ such that $T_\de(n)$ are universally bounded by $\ep$. We first consider $T_\de(1)$. Since $f$ is $K$-quasiconformal and normalized, we know $f(z)$ is bounded by $K$ for $0\leq z\leq 1$. Then by $m$ fixed and $g$ being $(\Omega,1+\de)$-quasisymmetric, we can find $\de$ and $R$ depending on $\ep$, $m$ and $K$ such that 
\begin{equation*}
    |g(f(\frac{a}{2^{m}}))-f(\frac{a}{2^{m}})|<\ep/2
\end{equation*}
for $a=1,2,\dots,2^m-1$, and $\de\to0$ as $\ep\to0$. Then we know $T_\de(1)\leq\ep/2$.
Now for $n$, divide $[0,1]$ into $2^{m}$ subintervals $[\frac{i}{2^{m}},\frac{i+1}{2^{m}}]$ with $i=0,1,\dots,2^{m}-1$. For each subinterval $[\frac{i}{2^{m}},\frac{i+1}{2^{m}}]$, we divide it into $2^{nm}$ pieces and renormalize $f$ and $g$ and use the result for $n$. Therefore we have
\begin{equation*}
    T_\de(n+1)\leq T_\de(1)+\frac{1}{2}T_\de(n)\leq \ep/2+\frac{1}{2}T_\de(n).
\end{equation*}
Together with $T_\de(1)<\ep/2$, we know $T_\de(n)<\ep$ for any positive integer $n$. Then by continuity of $f$ and $g$, we know for any $x\in[0,1]$,
\begin{equation*}
    |g(f(x))-f(x)|<\ep.\qedhere
\end{equation*}
\end{proof}


\begin{proof}[Proof of \ref{7.30}]

Let $f$ be a $K$-quasiconformal mapping from $\hat{\C}$ to $\hat{\C}$ that sends the unit circle $S^1$ to $\gamma$. We normalize $f$ and $g$ such that $f(-1)=-1=g(-1)$ and $f(1)=1=g(1)$. And we only need to prove that $g$ can be extended to the interior of $f(S^1)$ where the extension is compatible with M\"obius transformations as $g$, since the other side can be proved by applying the inversion along $S^1$. We take four points $A,B,C,D$ on $f(S^1)$ such that the modulus of the quadrilateral $ABCD$ is 1, which means the extremal distance from arc $\hu{AB}$ to arc $\hu{CD}$ within $f(\overline{\mathbb{D}})$ is 1. And we first want to prove that the extremal distance $b$ from arc $\hu{g(A)g(B)}$ to arc $\hu{g(C)g(D)}$ within $g\circ f(\overline{\mathbb{D}})$ is close to 1 when $\de$ is small.

For $\al>-1$, consider the circle $S_{1+\al}$ centered at the origin with radius $1+\al$ and the disc $\mathbb{D}_{1+\al}$ bounded by $S_{1+\al}$. For any $z\in S^1$, let $z'=(1+\al)z\in S_{1+\al}$. For any $K$, we know there exists $C(K,\al)>0$ such that for any $z_1,z_2\in \overline{\mathbb{D}_{1+\al}}$, we have
\begin{equation}\label{7.32}
    |f(z_1)-f(z_2)|<C(K,\al),
\end{equation}
whenever $|z_1-z_2|<\al$, and $C_1(\al)\to0$ as $\al\to0$. Applying \eqref{7.32} for $f^{-1}$ and by Lemma \ref{7.31}, for any $\al>0$ there exists $\de>0$ and compact $\Omega$ such that
\begin{equation}\label{7.33}
    |f^{-1}\circ g\circ f(z)-z|=|f^{-1}\circ g\circ f(z)-f^{-1}(f(z))|<\al,
\end{equation}
for $g$ being $(\Omega,\de)$-quasisymmetric and $z\in S^1$. Then we know $g(f(S^1))\subset f(\mathbb{D}_{1+\al}-\mathbb{D}_{1-\al})$.

\begin{figure}[H]
\centering
\includegraphics[scale=0.37]{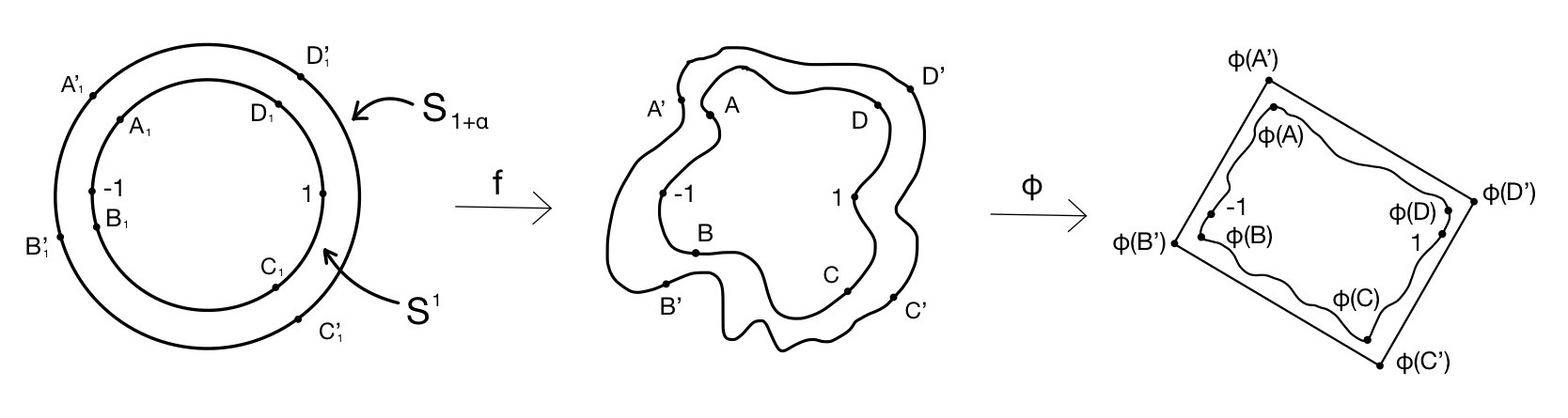} 
\caption{$S^1,S_{1+\al}$ and their images under $f$ and $\phi\circ f$}
\end{figure}

Let $A_1,B_1,C_1,D_1\in S^1$ be the preimage of $A,B,C,D$ under $f$, $A'_1,B'_1,C'_1,D'_1$ be the corresponding points on $S_{1+\al}$, and $A',B',C',D'$ be the image of $A'_1,B'_1,C'_1,D'_1$ under $f$. We first want to study the extremal distance $a$ from arc $\hu{A'B'}$ to arc $\hu{C'D'}$ within $f(\mathbb{D}_{1+\al})$. Let $\phi$ be the conformal mapping from $f(\mathbb{D}_{1+\al})$ to a rectangle $R$ in $\C$ such that $\phi(A'),\phi(B'),\phi(C'),\phi(D')$ are four vertices, $|\phi(A')\phi(D')|=a|\phi(A')\phi(B')|$ and $\phi$ preserves $-1$ and $1$. Since $\phi$ is conformal, so $\phi\circ f$ is $K$-quasiconformal on the disk $\mathbb{D}_{1+\al}$. On the other hand, the modulus of quadrilateral $ABCD$ is 1, and therefore the modulus of quadrilateral $A_1B_1C_1D_1$ is between $1/K$ and $K$. And so is the modulus of quadrilateral $A'_1B'_1C'_1D'_1$. Thus the modulus of quadrilateral $A'B'C'D'$ is between $1/K^2$ and $K^2$, and so is the modulus of $R$, which means $1/K^2<a<K^2$. Thus the rectangle $R$ admits a $\pi a$-quasiconformal reflection, which is proven in \cite{Wer97}. And then $\phi\circ f$ extends to a $\pi a K$-quasiconformal mapping $\overline{\phi\circ f}$ on $\hat{\C}$. Since $\pi a K<\pi K^3$, so $\overline{\phi\circ f}$ is $\pi K^3$-quasiconformal.

By $\overline{\phi\circ f}$ also preserves $-1$ and $1$, thus by (\ref{7.32}) for any $z\in S^1$,
\begin{equation}\label{7.34a}
    |\phi\circ f(z)-\phi\circ f(z')|<C(\pi K^3,\al),
\end{equation}
Then for any arc $\be$ inside $\phi(f(\mathbb{D}))$ connecting arc $\hu{\phi(A)\phi(B)}$ and arc $\hu{\phi(C)\phi(D)}$, we can extend $\be$ to $\be'$ which connects $\overline{\phi(A')\phi(B')}$ and $\overline{\phi(C')\phi(D')}$ and 
\begin{equation*}
    l(\be')<l(\be)+2C(\pi K^3,\al).
\end{equation*}
Let $|\phi(A')\phi(B')|=t$, then $|\phi(A')\phi(D')|=at$ and $l(\be')\geq at$. Hence $l(\be)>at-2C(\pi K^3,\al)$. Since the extremal distance from arc $\hu{AB}$ to arc $\hu{CD}$ within $f(\mathbb{D})$ is 1 and $\phi$ is conformal, so
\begin{equation}\label{7.34}
\begin{aligned}
    1&=\sup\limits_{\rho}\frac{\inf\limits_{\beta}L_\rho(\beta)}{\textup{Area}(\rho)}
    \geq \frac{\inf\limits_{\beta}l(\beta)}{\textup{Area}(\phi\circ f(\mathbb{D}))}
    \geq \frac{l(\be)^2}{at^2}\\
    &>\frac{(at-2C(\pi K^3,\al))^2}{at^2}=a-\frac{4C(\pi K^3,\al)}{t}+\frac{4C(\pi K^3,\al)^2}{at^2}\geq a-\frac{4C(\pi K^3,\al)}{t},
\end{aligned}
\end{equation}
where $\rho$ goes through all metrics and $\beta$ is among all arc connecting $\hu{\phi(A)\phi(B)}$ and $\hu{\phi(C)\phi(D)}$ within $\phi\circ f(\mathbb{D})$.
The rectangle $R$ contains -1 and 1, so its diagonal has length at least 2. Then by $1/K^2<\al<K^2$, we have
\begin{equation*}
    4\leq t^2+a^2 t^2=(1+a^2)t^2\leq (1+K^4)t^2.
\end{equation*}
Thus $t\geq\frac{2}{\sqrt{1+K^4}}$. Together with (\ref{7.34}), we know
\begin{equation*}
    1> a-\frac{4C(\pi K^3,\al)}{t}\geq a-2\sqrt{1+K^4}C(\pi K^3,\al)=:a-u(K,\al).
\end{equation*}
Hence
\begin{equation}\label{7.35}
    a<1+u(K,\al),
\end{equation}
and 
\begin{equation*}
    \lim\limits_{\al\to0}u(K,\al)=0,
\end{equation*}
for any fixed $K$. Similarly, we can consider the opposite pair of sides of the quadrilateral and we will have 
\begin{equation}\label{7.37}
    1/a<1+u(K,\al).
\end{equation}

By (\ref{7.33}), we know 
\begin{equation*}
    |f^{-1}\circ g\circ f(z)-z'|\leq|f^{-1}\circ g\circ f(z)-z|+|z'-z|<2\al.
\end{equation*}
Thus by (\ref{7.34a}),
\begin{equation*}
    |\phi(g(f(z)))-\phi(f(z'))|=|\phi(f(f^{-1}(g(f(z)))))-\phi(f(z'))|<C(\pi K^3,2\al).
\end{equation*}
Then since we know $g(f(S^1))\subset f(\mathbb{D}_{1+\al}-\mathbb{D}_{1-\al})$, we can apply the above method to the curve $g\circ f(S^1)$ and we have the following inequalities:
\begin{equation}\label{7.38}
    a<b+u(K,2\al),
\end{equation}
\begin{equation}\label{7.39}
    1/a<1/b+u(K,2\al),
\end{equation}
where $b$ is defined at the end of the first paragraph.

Now by (\ref{7.35}), (\ref{7.37}), (\ref{7.38}) and (\ref{7.39}), we have
\begin{equation}\label{7.40}
    \frac{1}{1+u(K,\al)}-u(K,2\al)<b<\left(\frac{1}{1+u(K,\al)}-u(K,2\al)\right)^{-1}.
\end{equation}
By $\lim\limits_{\al\to0}u(K,\al)=0$, we know
\begin{equation}\label{7.41}
    \lim\limits_{\al\to0^+}\left(\frac{1}{1+u(K,\al)}-u(K,2\al)\right)=1.
\end{equation}

Next we let $R_0$ be the interior of $\gamma$ and $R_1$ be the interior of $g(\gamma)$, and $\overline{R_0}$ and $\overline{R_1}$ are the corresponding closures. Then for $i=0,1$, take a homeomorphism $h_i:\overline{R_i}\to\overline{\mathbb{D}}$ such that $h_i$ is conformal on $R_i$. Here $h_i$'s are unique up to M\"obius transformations. And then we have a map $\eta:= h_2\circ g\circ h_1^{-1}:S^1\to S^1$. By $h_1$ conformal, we know the quadrilateral $h_1^{-1}(A)h_1^{-1}(B)h_1^{-1}(C)h_1^{-1}(D)$ has modulus 1. And by \eqref{7.40}, \eqref{7.41} and $h_2$ conformal, the quadrilateral $h_2(g(A))h_2(g(B))h_2(g(C))h_2(g(D))$ has modulus approaching to 1 uniformly as $\al\to1$. Thus for any $M_1>1$, there exists $\al>0$ such that $\eta$ is $M_1$-quasi-symmetric. On the other hand, for any $M>1$, there exists $M_1>1$ such that: when $\eta$ is $M_1$-quasi-symmetric, its Douady-Earle extension $E(\eta):\overline{\mathbb{D}}\to\overline{\mathbb{D}}$ is $M$-quasiconformal on $\mathbb{D}$. Moreover $\hat{g}:=h_2^{-1}\circ E(\eta)\circ h_1:\overline{R_0}\to\overline{R_1}$ is $M$-quasiconformal on $R_1$, which is the desired extension of $g$.

Finally, we only need to verify that the above extension has the natural property with M\"obius transformations when $g$ has. Suppose $\tau:G\to G'$ is an isomorphism between groups of M\"obius transformations, where $G$ preserves $\gamma$ and $G'$ preserves $g(\gamma)$, such that for any $x\in G$,
\begin{equation}\label{7.42}
    g\circ x=\tau(x)\circ g,
\end{equation}
as maps from $\gamma$ to $g(\gamma)$. Since we know $h_1\circ x\circ h_1^{-1}:\mathbb{D}\to\mathbb{D}$ is conformal, so there exists a M\"obius transformation $x_1$ such that $h_1\circ x\circ h_1^{-1}=x_1$. Thus we have
\begin{equation}\label{7.43}
\begin{aligned}
    h_1\circ x&=x_1\circ h_1,\\
    x\circ h_1^{-1}&=h_1^{-1}\circ x_1.
\end{aligned}
\end{equation}
Similarly, there is another M\"obius transformation $x_2$ such that
\begin{equation}\label{7.44}
\begin{aligned}
    h_2\circ \tau(x)&=x_2\circ h_2,\\
    \tau(x)\circ h_2^{-1}&=h_2^{-1}\circ x_2.
\end{aligned}
\end{equation}
Hence by the natural property of Douady-Earle extension with M\"obius transformations, we have
\begin{equation*}
\begin{aligned}
    \hat{g}\circ x&=h_2^{-1}\circ E(\eta)\circ h_1\circ x\\
    &\xlongequal{\eqref{7.43}}h_2^{-1}\circ E(\eta)\circ x_1\circ h_1\\
    &=h_2^{-1}\circ E(\eta\circ x_1)\circ h_1\\
    &=h_2^{-1}\circ E(h_2\circ g\circ h_1^{-1}\circ x_1)\circ h_1\\
    &\xlongequal{\eqref{7.43}}h_2^{-1}\circ E(h_2\circ g\circ x\circ h_1^{-1})\circ h_1\\
    &\xlongequal{\eqref{7.42}}h_2^{-1}\circ E(h_2\circ \tau(x)\circ g\circ h_1^{-1})\circ h_1\\
    &\xlongequal{\eqref{7.44}}h_2^{-1}\circ E(x_2\circ h_2\circ g\circ h_1^{-1})\circ h_1\\
    &=h_2^{-1}\circ E(x_2\circ \eta)\circ h_1\\
    &=h_2^{-1}\circ x_2\circ E(\eta)\circ h_1\\
    &\xlongequal{\eqref{7.44}}\tau(x)\circ h_2^{-1}\circ E(\eta)\circ h_1=\tau(x)\circ\hat{g}.\qedhere
\end{aligned}
\end{equation*}

\end{proof}

\subsection{Good is close to perfect}

We conclude this section as the following theorem.

\begin{thm}\label{7.47}
For all $B^+>B^->0,\de>0$ and $K>1$, there exists $R_0$ such that for all $M>1$, there exists $\ep>0$ such that for all $\overline{R}>R_0$: Suppose $R_i,s_i$ satisfy $|R_i-\overline{R}|<\de$, $B^-<\re(s_i)<B^+$, and suppose there exist $\overline{R}-\de<R'_i<\overline{R}+\de$ and $B^-<s'_i<B^+$ such that $(R_i,s_i)_{i=1}^3$ is $K$-related to $(R'_i,s'_i)_{i=1}^3$ and $|\re(s_i)-s'_i|<1/\overline{R}$. Let $\Gamma$ be the genus-2 quasi-Fuchsian group corresponding to the $(R_i,s_i)_{i=1}^3$-perfect assembly of 2 components. For any $(R_i,s_i,\ep)_{i=1}^3$-good assembly $\mathcal{A}$ in the hyperbolic 3-manifold $M^3$ such that $S_{\mathcal{A}}$ is connected, let $\rho_{\mathcal{A}}$ be the corresponding surface subgroup representation. Then $\rho_{\mathcal{A}}$ is $M$-quasiconformally conjugate to a finite index subgroup of $\Gamma$.
\end{thm}

\begin{proof}
Since $(R_i,s_i)_{i=1}^3$ is $K$-related to $(R'_i,s'_i)_{i=1}^3$, we know $\Gamma$ is $K$-quasi-Fuchsian. Hence the radius of its convex core is bounded by a constant $A$ only depending on $K$. Therefore all $(R_i,s_i)$-perfect assemblies are $A$-dense in the convex hull of $\Gamma$. And then the theorem follows from Theorem \ref{7.25}, \ref{7.26} and \ref{7.30}.
\end{proof}

\section{Proof of the main result}\label{S8}

\begin{proof}[Proof of Theorem \ref{mr}]

By Theorem \ref{pd1}, we take an $(R,m)$-good non-separating pants decomposition of $\Gamma$ for some $R,m>0$. Let the cuff lengths are $r_i$, $I=1,2,3$. Then there exists $K_1>0$ such that $\Gamma$ is $K_1$-quasiconformally conjugate to a Fuchsian group $\Gamma'$ with a non-separating pants decomposition of cuff lengths $\re(r_i)$, $i=1,2,3$. Then by Theorem \ref{gpdwbs}, there exists $B^+>B^->0$, such that for any $R_0>0$, there exists $\overline{R}>0,\de>0$ such that $\Gamma$ admits an $(\overline{R},\de)$-good pants decomposition with cuff lengths $R_i$ and twists $\re(s_i)\in(B^-,B^+)$, so we can choose $\overline{R}$ sufficiently large with all previous results related to $\overline{R}$ holding. Similarly $\Gamma'$ admits a corresponding pants decomposition $(R'_i,s'_i)$ by the $K_1$-quasiconformal mapping, and $s'_i\in (B^-,B^+)$. Moreover, by Lemma \ref{lem5.14}, we know 
\begin{equation*}
\begin{aligned}
    |\re(s_i)-s'_i|&\leq|\re(s_i)-\frac{\re(r_{i+1})+\re(r_{i+2})-\re(r_i)}{2}|+|s'_i-\frac{\re(r_{i+1})+\re(r_{i+2})-\re(r_i)}{2}|\\
    &<\frac{1}{\overline{R}},
\end{aligned}
\end{equation*}
when $\overline{R}$ is large enough. Then for $B^-,B^+,\de,\max\{K_1,K\},R_0$ and given $M>1$, let $\ep$ in Theorem \ref{7.47}.

Let
\begin{equation*}
    \textbf{A}=\sum_{P\in\Pi_{\ep,R_i}}P
\end{equation*}
be the formal sum of all unoriented $(R_i,\ep)_{i=1}^3$-good pants in the hyperbolic 3-manifold $M^3$. Then we want to use the \textit{doubling trick} in Section 2.5 in \cite{KW21} to match all the oriented pants together to construct a closed $(R_i,s_i,\ep)_{i=1}^3$-good assembly $\mathcal{A}$. To be more specifically, let $\gamma$ be an $(R_i,\ep)$-good curve for some $i$, then let $A_\gamma$ be the formal sum of one of each pants in $\Pi_{\ep,R_i}(\gamma)$ and $\sigma_{\gamma}$ be the permutation in Theorem \ref{matchpants}. Then we take $2A_{\gamma}$, where there are two copies with opposite orientation for each pants, and divide it into those $\pi$ with $\partial\pi$ the same orientation as $\gamma$ and those $\pi$ with $\partial\pi$ the opposite orientation. Thus each $\pi\in A_{\gamma}$ has an $\pi_+$ and an $\pi_-$, and we define the involution $\tau$ on $2A_{\gamma}$ by $\tau(\pi_+)=(\sigma_{\gamma}(\pi_+))_-$ and $\tau(\pi_-)=(\sigma^{-1}_{\gamma}(\pi_-))_+$. Then this involution gives us the way gluing pants together along each cuff and results in a closed oriented assembly.

If $S_{\mathcal{A}}$ is not connected, we can take a connected component, and we still denote it by $S_{\mathcal{A}}$ by passing to a subassembly. Then by Theorem \ref{7.47}, we know $\rho_{\mathcal{A}}$ is $K$-quasiconformally conjugate to a Fuchsian group $\rho_{\hat{\mathcal{A}}}$, where $\hat{\mathcal{A}}$ is the perfect model of $\mathcal{A}$. 

We know all pants in $\hat{\mathcal{A}}$ are identical with the pants of cuff lengths $R_1,R_2$ and $R_3$. Moreover, by Theorem \ref{matchpants}, these pants are glued by shear $s_i$ along the curve with length $R_i$. Therefore $S_{\hat{\mathcal{A}}}$ is a finite covering of our original quasi-Fuchsian surface $\HH^3/\Gamma$, which means $\rho_{\hat{\mathcal
A}}$ is a finite index subgroup of $\Gamma$.

\end{proof}

\begin{rem}\label{rem8.1}

Since we require that $\Gamma$ is a genus-2 Fuchsian group, so a non-separating pants decomposition of $\Gamma$ will have two identical pants. That is the reason why we can match those pants together by a permutation and the doubling trick, without considering the number of two different perfect pants. In general, it is possible to use good pants homology to solve this difficulty.

\end{rem}

We will end this paper by proving Theorem \ref{1.2}, which will follow from Theorem \ref{mr} and some general theorems about Hausdorff dimension.

\begin{proof}[Proof of Theorem \ref{1.2}]

By Section 7 of \cite{Ru82}, we know that the Hausdorff dimension of the limit set is a real analytic function on the deformation space of any quasi-Fuchsian group. In particular, the function is continuous. So by Theorem 1.3 in \cite{Bro03}, for any $1\leq\al<2$ and $\ep>0$, we can find a genus-2 quasi-Fuchsian group $\Gamma$ such that
\begin{equation*}
    \textup{H-dim}(\lambda(\Gamma))=\al.
\end{equation*}
By Theorem 8 and 12 in \cite{GV73} and Corollary 1.2 in \cite{Ast94}, we can find $K>1$ such that for any $K$-quasiconformal mapping $f:\hat{\C}\to\hat{C}$, we have
\begin{equation*}
    |\textup{H-dim}(f(\lambda(\Gamma)))-\al|<\ep.
\end{equation*}
Hence the result follows from Theorem \ref{mr}.

\end{proof}

\bibliographystyle{alpha}
\bibliography{reference}

\end{document}